\documentclass{amsart}
\usepackage{amsmath,amsthm,longtable}
\usepackage{graphicx}
\usepackage{float}

\newtheorem{claim}{Claim}
\newtheorem{theorem}[claim]{Theorem}
\newtheorem{proposition}[claim]{Proposition}
\newtheorem{corollary}[claim]{Corollary}
\newtheorem{lemma}[claim]{Lemma}
\newtheorem{observation}[claim]{Observation}
\newtheorem{remark}[claim]{Remark}

\theoremstyle{definition}
\newtheorem{definition}[claim]{Definition}

\newcommand{\mat}[4]{\left(\begin{smallmatrix} #1 & #2 \\ #3 & #4
\end{smallmatrix}\right)}
\newcommand{\eqs}{\overset{\text{s}}{\sim}} 
\newcommand{\eqf}{\overset{\text{F}}{\sim}} 
\newcommand{\spq}[2]{S^{\,#1}_{#2}}
\newcommand{\bspq}[2]{{\overline S}^{\,#1}_{#2}}
\newcommand{\Lc}{\text{L}}
\newcommand{\Hc}{\text{H}}

\begin{document}

\title{On pattern-avoiding partitions}
\author[V\'\i t Jel\'\i nek]{V\'\i t Jel\'\i nek$^\dag$}
\address{Department of Applied Mathematics, Charles University, Prague}
\email{jelinek@kam.mff.cuni.cz}
\thanks{$^\dag$Supported by the project MSM0021620838 of the Czech Ministry of Education, and by the grant GD201/05/H014
of the Czech Science Foundation.}

\author{Toufik Mansour}
\address{Department of Mathematics, Haifa University, 31905 Haifa, Israel}
\email{toufik@math.haifa.ac.il}

\subjclass[2000]{Primary 05A18; Secondary
05E10, 05A15, 05A17, 05A19}
\keywords{Ferrers shapes, Ordered graphs, Set partitions,
Pattern-avoidance, Polyomino shapes}
\floatstyle{boxed}
\restylefloat{figure}

\begin{abstract}
A \emph{set partition} of the set $[n]=\{1,\dotsc,n\}$ is a
collection of disjoint blocks $B_1,B_2,\dotsc, B_d$ whose union is
$[n]$. We choose the ordering of the blocks so that they satisfy
$\min B_1<\min B_2<\dotsb<\min B_d$. We represent such a set
partition by a \emph{canonical sequence} $\pi_1,\pi_2,\dotsc,\pi_n$,
with $\pi_i=j$ if $i\in B_j$.  We say that a partition $\pi$
\emph{contains} a partition $\sigma$ if the canonical sequence of
$\pi$ contains a subsequence that is order-isomorphic to the
canonical sequence of $\sigma$. Two partitions $\sigma$ and
$\sigma'$ are \emph{equivalent}, if there is a size-preserving
bijection between $\sigma$-avoiding and $\sigma'$-avoiding
partitions.

We determine several infinite families of sets of equivalent patterns;
for instance, we prove that there is a bijection between $k$-noncrossing and
$k$-nonnesting partitions, with a notion of crossing and nesting based on the
canonical sequence. We also provide new combinatorial interpretations of the
Catalan numbers and the Stirling numbers. Using a systematic computer search,
we verify that our results characterize all the pairs of equivalent
partitions of size at most seven.

We also present a correspondence between set partitions and fillings of Ferrers
shapes and stack polyominoes. This correspondence allows us to apply recent results
on polyomino fillings in the study of partitions, and conversely, some of our
results on partitions imply new results on polyomino fillings and ordered graphs.
\end{abstract}

\maketitle
\section{Introduction}
A \emph{partition} of $[n]=\{1,2,\ldots ,n\}$ is a collection
$B_{1},B_{2},\ldots ,B_{d}$ of nonempty disjoint sets, called \emph{blocks},
whose union is $[n]$. We will assume that $B_{1},B_{2},\ldots ,B_{d}$ are
listed in the increasing order of their minimum elements, that is, $\min
B_{1}<\min B_{2}<\cdots<\min B_{d}$. In this paper, we will represent a partition
of $[n]$ by its \emph{canonical sequence}, which is an integer sequence
$\pi=\pi_1\pi_2\dotsb\pi_n$ such that $\pi_i=k$ if and only if $i\in B_k$.
For instance, $1231242$ is the canonical sequence of the partition of
$\{1,2,\dotsc,7\}$ with the four blocks $\{1,4\}$, $\{2,5,7\}$, $\{3\}$ and~$\{6\}$.

Note that a sequence $\pi$ over the  alphabet $[d]$ represents a
partition of $[n]$ with $d$ blocks if and only if it has the
following properties:
\begin{itemize}
\item Each number from the set $[d]$ appears at least once in $\pi$.
\item  For each $i,j$ such that $1\le i<j\le d$, the first occurrence of $i$
precedes the first occurrence of~$j$.
\end{itemize}
We remark that sequences satisfying these properties are also known as
\emph{restricted growth functions}, and they are often encountered in the
study of set partitions \cite{Sag,WW} as well as other related topics, such
as Davenport-Schinzel sequences \cite{DS,K1,K2,MS}.

Throughout this paper, we identify a set partition with the
corresponding canonical sequence, and we use this representation to
define the notion of pattern avoidance among set partitions. Let
$\pi=\pi_1 \pi_2\dotsb\pi_n$ and
$\sigma=\sigma_1\sigma_2\dotsb\sigma_m$ be two partitions
represented by their canonical sequences. We say that $\pi$
\emph{contains} $\sigma$, if $\pi$ has a subsequence that is
order-isomorphic to $\sigma$; in other words, $\pi$ has a
subsequence $\pi_{f(1)},\pi_{f(2)},\dotsc,\pi_{f(m)}$, where $1\leq
f(1)<f(2)<\dotsb< f(m)\leq n$, and for each $i,j\in[m]$,
$\pi_{f(i)}<\pi_{f(j)}$ if and only if $\sigma_i<\sigma_j$. If $\pi$
does not contain $\sigma$, we say that $\pi$ \emph{avoids} $\sigma$.
Our aim is to study the set of all the partitions of $[n]$ that
avoid a fixed partition $\sigma$. In such context, $\sigma$ is
usually called a \emph{pattern}.

Let $P(n)$ denote the set of all the partitions of $[n]$, let $P(n;\sigma)$
denote the set of all partitions of $[n]$ that avoid $\sigma$, and let $p(n)$
and $p(n;\sigma)$ denote the cardinality of $P(n)$ and $P(n;\sigma)$,
respectively. We say that two partitions $\sigma$ and $\sigma'$ are
\emph{equivalent}, denoted by $\sigma\sim\sigma'$, if
$p(n;\sigma)=p(n;\sigma')$ for each $n$.

The concept of pattern-avoidance described above has been introduced by
Sagan~\cite{Sag} and later developed by Goyt~\cite{Goyt}. These two papers
considered, among other topics, the enumeration of partitions avoiding
patterns of size three. In our paper, we extend this study to larger
patterns. Rather than studying isolated cases of equivalent pairs of
patterns, our aim is to describe infinite families of nontrivial equivalence
classes of partitions.

For instance, we define $k$-noncrossing and $k$-nonnesting
partitions as the partitions that avoid the pattern $12\dotsb
k12\dotsb k$ and $12\dotsb k k (k-1) \dotsb 1$, respectively. We
will show that these two patterns are equivalent for every $k$, by
constructing a bijection between $k$-noncrossing and
$k$-nonnesting partitions. It is noteworthy, that a different
concept of crossings and nestings in partitions has been considered
by Chen et al.~\cite{C1,C2}, and this different notion of crossings
and nestings also admits a bijection between $k$-noncrossing and
$k$-nonnesting partitions, as has been shown in~\cite{C2}. There is,
in fact, yet another notion of crossings and nestings in partitions
that has been extensively studied by Klazar~\cite{K1,K2}.

Several of our results are proved using a correspondence between partitions
and 0-1 fillings of polyomino shapes. This correspondence allows us to
translate recent results on fillings of Ferrers shapes~\cite{DS,kra} and
stack polyominoes~\cite{rub} into the terminology of pattern-avoiding
partitions. The correspondence between fillings of shapes and
pattern-avoiding partitions works in the opposite way as well: some of our
theorems, proved in the context of partitions, imply new results about
pattern-avoiding fillings of Ferrers shapes and pattern-avoiding ordered
graphs.

Apart from these results, we also present a class of patterns
equivalent to the pattern $12\dotsb k$. This result can be viewed as
a new combinatorial interpretation of the Stirling numbers of the
second kind. Similarly, by providing patterns equivalent to $1212$,
we provide a new combinatorial interpretation of the Catalan
numbers.

To test the strength of our general equivalence theorems, we have undertaken
systematic computer enumeration of partitions avoiding single patterns of
size at most 7; the results of this enumeration are presented in the appendix
of this paper. Our methods are able to completely characterize the
equivalence of all the patterns of length at most seven. This extends earlier
results of Sagan~\cite{Sag}, who provided similar characterization of
patterns of size three.

The paper is organized as follows: in Section~\ref{sec-basic}, we present
basic facts about pattern-avoiding partitions, and we summarize previously
known results. Our main results are collected in Section~\ref{sec-main},
where we present several infinite families of classes of equivalent
patterns. In the rest of the paper, we deal with several isolated pairs of
patterns which are not covered by any of the general theorems of the previous
section. In particular, in Section~\ref{sec-1123}, we prove that the pattern
$1123$ is equivalent to the pattern $1212$, thus completing the
characterization of the patterns of size four and obtaining a new
interpretation for the Catalan numbers. In Section~\ref{sec-12112}, we prove
the equivalence $12112\sim 12212$, and discuss its implications for the
theory of pattern-avoiding ordered graphs and polyomino fillings.

\section{Basic facts and previous results}\label{sec-basic}

Let us first establish some notational conventions that will be applied
throughout this paper: for a finite sequence $S=s_1s_2\dotsb s_p$ and an
integer $k$, we use the notation $S+k$ to refer to the sequence
$(s_1+k)(s_2+k)\dotsb(s_p+k)$. For a symbol $k$ and an integer $d$, the
constant sequence $(k,k,\dotsc,k)$ of length $d$ is denoted by $k^d$.
To prevent confusion, we will use capital letters $S,T,\dotsc$ to denote
arbitrary sequences of positive integers, and we will use lowercase greek symbols
($\pi,\sigma,\tau,\dotsc$) to denote sequences in their canonical form
representing partitions.

An infinite sequence $a_0,a_1,\dotsc$ is often
conveniently represented by its \emph{exponential generating function} (or
EGF for short), which is the formal power series $F(x)=\sum_{n\ge 0}\frac{a_n
x^n}{n!}$. We are mostly interested
in the generating functions of the sequences of the form
$\left(p(n;\pi)\right)_{n\ge 0}$, where $\pi$ is a given pattern. We simply call
such a generating function \emph{the EGF of the pattern
$\pi$}.

\subsection{Simple patterns}
Several natural classes of partitions can be defined in terms of
pattern-avoidance. For instance, the partitions whose every block
is smaller than a given constant $k$ are precisely the
partitions that avoid the pattern $\pi=11\dotsb1=1^k$. Using
standard combinatorial arguments (see, e.g., \cite{fs}), we obtain
the following well-known formula for the EGF of the pattern $1^k$:
\begin{equation}\label{eq-1m}
F(x)=\exp(\exp_{<k}(x)-1),
\end{equation}
where $\exp(x)=\sum_{n\ge 0} \frac{x^n}{n!}$ and $\exp_{<k}(x)=\sum_{n=0}^{k-1}
\frac{x^n}{n!}$.

Another class of pattern-avoiding partitions with a similarly natural
combinatorial description is the set of the partitions avoiding $\pi=12\dotsb k$.
Clearly, a partition avoids $\pi$ if and only if it has less than $k$ blocks.
It is well known (see, e.g.,~\cite{fs}) that the EGF of the pattern
$12\dotsb k$ is equal to
\begin{equation}\label{eq-123k}
\exp_{<k}(\exp(x)-1).
\end{equation}
The enumeration of $\pi$-avoiding partitions is closely related to
the Stirling numbers of the second kind $S(n,m)$, defined as the
number of partitions of $[n]$ with exactly $m$ blocks (see sequence
A008277 in~\cite{oeis}).

\subsection{Patterns of size three}
Sagan~\cite{Sag} has described and enumerated the pattern-avoiding
classes $P(n;\pi)$ for the five patterns $\pi$ of length three. Here
we briefly summarize the relevant results (see Table~\ref{tab1}).
For the sake of completeness, we briefly present the arguments used
to obtain these results.
\begin{table}[h]
\begin{center}
\begin{tabular}{ l|l }
  $\tau$ & $p(n;\tau)$ \\ \hline\hline
  $111$ & sequence A000085 in~\cite{oeis}\\
  $112$, $121$, $122$, $123$ & $2^{n-1}$ \\[4pt]
\end{tabular}
\end{center}
\caption{Number of partitions in $P(n;\tau)$, where $\tau\in P(3)$.}\label{tab1}
\end{table}

Clearly, a partition avoids $111$ if and only if each of its blocks has size at
most two. This is a special case of the pattern $1^k$ discussed above. We
remark that the sequence $p(n;111)$ enumerating these partitions has several
other combinatorial interpretations, such as the number of involutions of the
set $[n]$ (see A000085 in \cite{oeis}).

The four remaining patterns of size three, namely $123$, $122$,
$121$ and $112$, are all equivalent and they satisfy
$p(n;\pi)=2^{n-1}$. To see this, observe first of all that the partitions
avoiding $123$ correspond precisely to the restricted growth
functions which start with the symbol 1 and any other symbol is
equal to $1$ or $2$.

Similarly, the partitions avoiding $122$ are encoded by sequences whose first
symbol is $1$, and each of the following symbols is either equal to
$1$ or is greater by one than the largest preceding symbol.

The case of the pattern $121$ is equally simple: each symbol of the
corresponding sequence after the first one is either equal to the
largest preceding symbol or is greater by one than the largest
preceding symbol.

For the pattern $112$, the argument is different: note that a
partition $\tau$ with $k$ blocks avoids $112$ if and only if $\tau$
consists of the sequence $12\dotsb k$ followed by a weakly
decreasing sequence, i.e., $\tau=12\dotsb kS$ where $S$ is weakly
decreasing. In particular, a $112$-avoiding partition is uniquely
determined by the ordered sequence of the sizes of its blocks. Thus,
the number of $112$-avoiding partitions of $[n]$ is equal to the
number of ordered sequences of positive integers whose sum is $n$,
and it is well known that there are precisely $2^{n-1}$ such
sequences.

\section{General classes of equivalent patterns}\label{sec-main}

In this section, we introduce the tools that will be useful in our study of
pattern-avoidance, and we prove our key results. We begin by introducing a
general relationship between pattern-avoidance in partitions and
pattern-avoidance in fillings of restricted shapes. This approach will
provide a useful tool for dealing with many pattern problems.

\subsection{Pattern-avoiding fillings of diagrams}

We will use the term \emph{diagram} to refer to any finite set of the cells
of the two-dimensional square grid. To \emph{fill} a diagram means to write a
non-negative integer into each cell.

We will number the rows of diagrams from bottom to top, so the ``first row''
of a diagram is its bottom row, and we will number the columns from left to
right. We will apply the same convention to matrices and to fillings. We
always assume that each row and each column of a diagram is nonempty; thus,
for example, when we refer to a diagram with $r$ rows, it is assumed that
each of the $r$ rows contains at least one cell of the diagram. Note
that there is a (unique) empty diagram with no rows and no columns.
Let $r(F)$ and $c(F)$ denote, respectively, the number of rows and columns of
$F$, where $F$ is a diagram, or a matrix, or a filling of a diagram.

We will mostly use diagrams of a special shape, namely Ferrers diagrams and
stack polyominoes. We begin by giving the necessary definitions:
\begin{definition}
A \emph{Ferrers diagram}, also called \emph{Ferrers shape}, is a diagram whose
cells are arranged into contiguous rows and columns satisfying the following
rules:
\begin{itemize}
\item The length of any row is greater than or equal to the length of any row
above it.
\item The rows are right-justified, i.e., the rightmost cells of the rows appear
in the same column.
\end{itemize}
\end{definition}

We admit that our convention of drawing Ferrers diagrams as right-justified
rather than left-justified shapes is different from standard practice;
however, our definition will be more intuitive in the context of our
applications.

\begin{definition}
A \emph{stack polyomino} $\Pi$ is a collection of finitely many cells of the
two-dimensional rectangular grid, arranged into contiguous rows and columns
with the property that for any $i=1,\dotsc, r(\Pi)$, every column intersecting
the $i$-th row also intersects all the rows with index smaller than~$i$.
\end{definition}

Clearly, every Ferrers shape is also a stack polyomino. On the other hand, a
stack polyomino can be regarded as a union of a Ferrers shape and a vertically
reflected copy of another Ferrers shape.

\begin{definition}
A \emph{filling} of a diagram is an assignment of non-negative integers to
the cells of the diagram. A \emph{0-1 filling} is a filling that only uses
values 0 and 1. In such filling, a \emph{0-cell} of a filling is a cell that
is filled with value 0, and a \emph{1-cell} is filled with value 1. A 0-1 filling
is called \emph{semi-standard} if each of its columns contains exactly one 1-cell. A 0-1 filling is called
\emph{sparse} if every column has at most one 1-cell. A column or row of a 0-1
filling is called \emph{zero column (or row)} if it contains no 1-cell.
\end{definition}

Among several possibilities to define pattern-avoidance in fillings, the
following approach seems to be the most useful and most common:

\begin{definition}
Let $M=(m_{ij}; i\in [r], j\in[c])$ be a matrix with $r$ rows and $c$ columns
with all entries equal to 0 or 1, and let $F$ be a filling of a
diagram. We say that $F$ \emph{contains} $M$ if $F$ contains $r$ distinct
rows $i_1<\dotsb<i_r$ and $c$ distinct columns $j_1<\dotsb<j_c$ with the
following two properties:
\begin{itemize}
\item Each of the rows $i_1,\dotsc,i_r$ intersects all columns $j_1,\dotsc,j_c$
in a cell that belongs to the underlying diagram of $F$.
\item If $m_{kl}=1$ for some $k$ and $l$, then the cell of $F$ in row $i_k$
and column $j_l$ has a nonzero value.
\end{itemize}
If $F$ does not contain $M$, we say that $F$ \emph{avoids} $M$.
We will say that two matrices $M$ and $M'$ are \emph{Ferrers-equivalent}
(denoted by $M\eqf M'$) if for every Ferrers shape $\Delta$, the number of
semi-standard fillings of $\Delta$ that avoid $M$ is equal to the number of
semi-standard fillings of $\Delta$ that avoid $M'$. We will say that $M$ and
$M'$ are \emph{stack-equivalent} (denoted by $M\eqs M'$) if the equality holds
even for semi-standard fillings of an arbitrary stack polyomino.
\end{definition}

Pattern-avoidance in the fillings of diagrams  has received considerable
attention lately. Apart from semi-standard fillings, various authors have
considered \emph{standard} fillings with exactly one 1-cell in each row and
each column (see \cite{bwx} or \cite{sw}), as well as general fillings with
non-negative integers (see \cite{adm} or~\cite{kra}). Also, nontrivial
results were obtained for fillings of more general shapes (e.g. moon
polyominoes~\cite{rub}). These results often consider the cases when the
forbidden pattern $M$ is the identity matrix (i.e., the $r$ by $r$ matrix,
with $m_{ij}=1$ if and only if $i=j$; this matrix will be denoted by $I_r$)
or the anti-identity matrix (i.e., the $r$ by $r$ matrix with $m_{ij}=1$ if
and only if $i+j=r+1$; this matrix will be denoted by $J_r$).

Since our next arguments mostly deal with semi-standard fillings, we will
drop the adjective `semi-standard' and simply use the term `filling', when
there is no risk of ambiguity.

\begin{remark}\label{rem-sparse}
In our argument, we will often encounter a mapping $f$ that transforms a
given semi-standard filling of a Ferrers diagram (or a stack polyomino) into
another semi-standard filling of the same diagram. It will be convenient to
extend such transform $f$ to act on sparse fillings as well as semi-standard
fillings; this is achieved in the following natural way: given a sparse
filling $F$ of a Ferrers diagram $\Delta$, we ignore all the columns of $F$
that contain no 1-cell and observe that the remaining columns induce a
semi-standard filling of a Ferrers diagram. We then transform the filling $F$
by letting the mapping $f$ act on the non-zero columns of $F$ (i.e. those that
contain a 1-cell), while the zero columns are left without change.

In particular, if $M$ and $M'$ are two Ferrers-equivalent 0-1 matrices with a
1-cell in every column, the argument above shows that there is a bijection
between $M$-avoiding and $M'$-avoiding sparse fillings of a given Ferrers
diagram. To see this, note that a sparse filling $F$ avoids $M$ if and only if its
subfilling induced by the nonzero columns avoids $M$, since $M$ has a 1-cell
in every column.

A completely analogous argument can be made for stack polyominoes instead of
Ferrers shapes.
\end{remark}

We now introduce some more notation, which will be useful for translating the language of
partitions to the language of fillings.
\begin{definition}\label{def-mat}
Let $S=s_1s_2\dotsb s_m$ be a sequence of positive integers, and let
$k\ge\max\{s_i;\;i\in[m]\}$ be an integer. We let $M(S,k)$ denote the 0-1
matrix with $k$ rows and $m$ columns which has a 1-cell in row $i$ and column
$j$ if and only if $s_j=i$.
\end{definition}

We now describe the correspondence between partitions and fillings
of Ferrers diagrams (recall that $\tau+k$ denotes the sequence
obtained from $\tau$ by adding $k$ to every element).

\begin{lemma}\label{lem-fil}
Let $S$ and $S'$ be two sequences over the alphabet $[k]$, let $\tau$ be an
arbitrary partition. If $M(S,k)$ is Ferrers-equivalent to $M(S',k)$ then the
partition pattern $\sigma=12\dotsb k(\tau+k)S$ is equivalent to
$\sigma'=12\dotsb k(\tau+k)S'$.
\end{lemma}
\begin{proof}
Let $\pi$ be a partition of $[n]$ with $m$ blocks. Let $M$ denote the matrix
$M(\pi,m)$. Fix a partition $\tau$ with $t$ blocks, and let $T$ denote the
matrix $M(\tau,t)$. We will color the cells of $M$ red and green in the
following way: if $\tau$ is nonempty, then a cell in row $i$ and column $j$
is colored green if and only if the submatrix of $M$ induced by the rows
$i+1,\dotsc, m$ and columns $1,\dotsc,j-1$ contains $T$. If $\tau$ is empty,
then a cell in row $i$ and column $j$ is green if and only if row $i$ has at
least one 1-cell strictly to the left of column $j$. A cell is red if it is
not green.

Note that the green cells form a Ferrers diagram, and the entries of the
matrix $M$ form a sparse filling $G$ of this diagram. Also note that the
leftmost 1-cell of each row is always red, and any 0-cell of the same row
to the left of the leftmost 1-cell is red too.

It is not difficult to see that the partition $\pi$ avoids $\sigma$ if and
only if the filling $G$ of the `green' diagram avoids $M(S,k)$, and $\pi$
avoids $\sigma'$ if and only if $G$ avoids $M(S',k)$. Since $M(S,k)\eqf
M(S',k)$, there is a bijection $f$ that maps $M(S,k)$-avoiding fillings of
Ferrers shapes onto $M(S',k)$-avoiding fillings of the same shape. By
Remark~\ref{rem-sparse}, $f$ can be extended to sparse fillings. Using this
extension of $f$, we construct the following bijection
between $P(n;\sigma)$ and $P(n;\sigma')$: for a partition $\pi\in
P(n;\sigma)$ with $m$ blocks, we take $M$ and $G$ as above. By assumption,
$G$ is $M(S,k)$-avoiding. Using the bijection $f$ and
Remark~\ref{rem-sparse}, we transform $G$ into an $M(S',k)$-avoiding sparse
filling $f(G)=G'$, while the filling of the red cells of $M$ remains the same. We
thus obtain a new matrix $M'$.

Note that if we color the cells of $M'$ red and green using the criterion
described in the first paragraph of this proof, then each cell of $M'$ will receive the
same color as the corresponding cell of $M$,
even though the occurrences of $T$ in $M'$ need not correspond exactly to the
occurrences of $T$ in $M$.

By construction, $M'$ has exactly one 1-cell in
each column, hence there is a sequence $\pi'$ over the alphabet $[m]$ such
that $M'=M(\pi',m)$. We claim that $\pi'$ is a canonical sequence of a
partition. To see this, note that for every $i\in[m]$, the leftmost 1-cell of
$M$ in row $i$ is red and the preceding 0-cells in row $i$ are red too. It
follows that the leftmost 1-cell of row $i$ in $M$ is also the leftmost
1-cell of row $i$ in $M'$, so the first occurrence of the symbol $i$ in $\pi$
appears at the same place as the leftmost occurrence of $i$ in $\pi'$, hence
$\pi'$ is indeed a partition. The green cells of $M'$ avoid $M(S',k)$, so
$\pi'$ avoids $\sigma'$. Obviously, the transform $\pi \mapsto  \pi'$ is
invertible and provides a bijection between $P(n;\sigma)$ and $P(n;\sigma')$.
\end{proof}

In general, the relation $12\dotsc k S\sim 12\dotsc k S'$ does not imply that
$M(S,k)$ and $M(S',k)$ are Ferrers equivalent: in Section~\ref{sec-12112}, we
will prove that $12112\sim 12212$, even though $M(2,112)$ is not Ferrers
equivalent to $M(2,212)$.

On the other hand, the relation $12\dotsc k S\sim 12\dotsc k S'$ allows us to
establish a somewhat weaker equivalence between pattern-avoiding fillings, using
the following lemma.

\begin{lemma}\label{lem-fil2}
Let $S$ be an arbitrary sequence over the alphabet $[k]$, and let
$\tau=12\dotsb k S$. For every $n$ and $m$, there is a bijection $f$
that maps the set of $\tau$-avoiding partitions of $[n]$ with $m$
blocks onto the set of all the $M(S,k)$-avoiding fillings $F$ of
Ferrers shapes that satisfy $c(F)=n-m$ and $r(F)\le m$.
\end{lemma}
\begin{proof}
Let $\pi$ be a $\tau$-avoiding partition of $[n]$ with $m$ blocks. Let
$M=M(\pi,m)$, and let us consider a red and green coloring of $M$ in the same
way as in the proof of Lemma~\ref{lem-fil}, i.e., the green cells of a row $i$ are precisely
the cells that are strictly to the right of the leftmost 1-cell in row $i$.

Note that $M$ has exactly $m$ red 1-cells, and each 1-cell is red if and only if
it is the leftmost 1-cell of its row. Note also that if $c_i$ is column
containing the red 1-cell in row $i$, then either $c_i$ is the rightmost column of $M$, or the
column $c_i+1$ is the leftmost column of $M$ with exactly $i$ green cells.

Let $G$ be the filling formed by the green cells. As was pointed out in the previous
proof, the filling $G$ is a sparse $M(S,k)$-avoiding filling of a Ferrers shape.
Note that for each $i=1,\dotsc m-1$, the filling $G$ has exactly one zero column
of height $i$, and this column, which corresponds to $c_{i+1}$, is the rightmost
of all the columns of $G$ with height at most $i$.

Let $G^-$ be the subfilling of $G$ induced by all the nonzero columns of $G$.
Observe that $G^-$ is a semi-standard $M(S,k)$-avoiding filling of a ferrers
shape with exactly $n-m$ columns and at most $m$ rows; we thus define
$f(\pi)=G^-$.

Let us now show that the mapping $f$ defined above can be inverted. Let $F$ be a
filling of a Ferrers shape with $n-m$ columns and at most $m$ rows. We
insert $m-1$ zero columns $c_2,c_3,\dotsc,c_m$ into the filling $F$ as follows:
each column $c_i$ has height $i-1$, and it is
inserted immediately after the rightmost column of $F\cup\{c_2,\dotsc,c_{i-1}\}$
that has height at most~$i-1$. Note that the filling obtained by this operation
corresponds to the green cells of the original matrix $M$, so let us color all its
cells green, and let us call this sparse filling $G$.

We now add a new red 1-cell on top of each zero column of $G$, and
we add a new red 1-cell in front of the bottom row, to obtain a
semistandard filling of a diagram with $n$ columns and $m$ rows.
which can be completed into a matrix $M=M(\pi,m)$, where $\pi$ is
easily seen to be a canonical sequence of a $\tau$-avoiding
partition.
\end{proof}

Lemma~\ref{lem-fil} provides a tool to deal with partition patterns of the
form $12\dotsb k(\tau+k) S $ where $S$ is a sequence over $[k]$ and $\tau$ is
a partition. We now describe a correspondence between partitions and fillings
of stack polyominoes, which is useful for dealing with patterns of the form
$12\dotsb kS(\tau+k)$. We use a similar argument as in the proof of
Lemma~\ref{lem-fil}.

\begin{lemma}\label{lem-poly} If $\tau$ is a partition, and $S$ and $S'$ are two
sequences over the alphabet $[k]$ such that $M(S,k)\eqs M(S',k)$,
then the partition $\sigma=12\dotsb kS(\tau+k)$ is
equivalent to the partition $\sigma'=12\dotsb kS'(\tau+k)$.
\end{lemma}

\begin{proof}
Fix a partition $\tau$ with $t$ blocks. Let $\pi$ be any partition of $[n]$
with $m$ blocks, let $M=M(\pi,m)$. We will color the cells of $M$
red and green as follows: a cell of $M$ in row $i$
and column $j$ is green, if it satisfies both these conditions:
\begin{itemize}
\item[(a)] The submatrix of $M$ formed by the intersection of rows
$i+1,i+2,\dotsc, m$ and columns $j+1,j+2,\dotsc,n$ contains $M(\tau,t)$.
\item[(b)] The matrix $M$ has at least one 1-cell in row $i$ appearing strictly to the
left of column~$j$.
\end{itemize}
A cell is called red, if it is not green. Note that the green cells form a stack
polyomino and the matrix $M$ induces a sparse filling $G$ of this polyomino.

Similarly to Lemma~\ref{lem-fil}, it easy to verify that the partition $\pi$
above avoids the pattern $\sigma$ if and only if the filling $G$
avoids $M(S,k)$, and $\pi$ avoids $\sigma'$ if and only if $G$ avoids
$M(S',k)$.

The rest of the argument is analogous to the proof of Lemma~\ref{lem-fil}.
Assume that $M(S,k)$ and $M(S',k)$ are stack-equivalent via a bijection $f$.
By Remark~\ref{rem-sparse}, we extend $f$ to a bijection between
$M(S,k)$-avoiding and $M(S',k)$-avoiding sparse fillings of a given stack
polyomino. Consider a partition $\pi\in P(n;\sigma)$ with $m$ blocks, and
define $M$ and $G$ as above. Apply $f$ to the filling $G$ to obtain an
$M(S',k)$-avoiding filling $G'$; the filling of the red
cells of $M$ remains the same. This yields a matrix $M'$ and a sequence $\pi'$ such that
$M'=M(\pi',k)$; we may easily check that the green cells of $M'$ are the same
as the green cells of $M$. By rule (b) above, the leftmost 1-cell of each row
of $M$ is unaffected by this transform. It follows that the first occurrence
of $i$ in $\pi'$ is at the same place as the first occurrence of $i$ in
$\pi$, and in particular, $\pi'$ is a partition. By the observation of the
previous paragraph, $\pi'$ avoids $\sigma'$ and the transform $\pi\mapsto\pi'$
is a bijection from $P(n;\sigma)$ to $P(n;\sigma')$.
\end{proof}

The following simple result about pattern-avoidance in fillings will turn out
to be useful in the analysis of pattern avoidance in partitions:

\begin{proposition}\label{pro-filshift}
If $S$ is a sequence over the alphabet $[k-1]$, then $M(S,k)$ is
stack-equivalent to $M(S+1,k)$. If $S$ and $S'$ are two sequences over $[k-1]$
such that $M(S,k-1)\eqf M(S',k-1)$ then $M(S,k)\eqf M(S',k)$, and if
$M(S,k-1)\eqs M(S',k-1)$ then $M(S,k)\eqs M(S',k)$.
\end{proposition}
\begin{proof}
To prove the first part, let us define $M=M(S,k)$, $M^-=M(S,k-1)$, and
$M'=M(S+1,k)$. Notice that a filling $F$ of a stack polyomino $\Pi$ avoids $M$
if and only if the filling obtained by erasing the topmost cell of every column
of $F$ avoids $M^-$. Similarly, $F$ avoids $M'$, if and only if the filling
obtained by erasing the bottom row of $F$ avoids $M^-$. We thus have the
following bijection between $M$-avoiding and $M'$-avoiding fillings: take an
$M$-avoiding filling $F$, and in every column, move the topmost element into the
bottom row, and move every other element into the row directly above it.
The second claim of the theorem is proved analogously.
\end{proof}

For convenience, we translate the first part of this proposition into the
language of pattern-avoiding partitions, using Lemma~\ref{lem-fil} and
Lemma~\ref{lem-poly}.

\begin{corollary}
If $S$ is a sequence over $[k-1]$ and $\tau$ is an arbitrary partition, then
$$12\dotsb k(\tau+k) S\sim 12\dotsb k(\tau+k) (S+1)\mbox{ and }12\dotsb
kS(\tau+k) \sim 12\dotsb k(S+1)(\tau+k).$$
\end{corollary}

We now state another result related to pattern-avoidance in Ferrers diagrams, which
has important consequences in  our study of partitions. Let us first fix the following
notation: for two matrices $A$ and $B$, let $\mat{A}{0}{0}{B}$ denote the matrix
with $r(A)+r(B)$ rows and $c(A)+c(B)$ columns with a copy of $A$ in the top left
corner and a copy of $B$ in the bottom right corner.

The idea of the following proposition is not new, it has already been applied by Backelin et
al.~\cite{bwx} to standard fillings of Ferrers diagrams, and later adapted
by de Mier \cite{adm} for fillings with arbitrary integers. We now apply it
to semi-standard fillings.

\begin{lemma}\label{lem-swap}
If $A$ and $A'$ are two Ferrers equivalent matrices, and if $B$ is an arbitrary
matrix, then $\mat{B}{0}{0}{A}\eqf\mat{B}{0}{0}{A'}$.
\end{lemma}
\begin{proof}
Let $F$ be an arbitrary $\mat{B}{0}{0}{A}$-avoiding filling of a Ferrers
diagram $\Delta$. We say that a cell in row $i$ and column $j$ of $F$ is
\emph{green} if the subfilling of $F$ induced by the intersection of rows
$i+1,i+2,\dotsc,r(F)$ and columns $1,2,\dotsc,j-1$ contains a copy of $B$.
Note that the green cells form a Ferrers shape $\Delta^-\subseteq
\Delta$, and that the restriction of $F$ to the cells of $\Delta^-$ is a
sparse $A$-avoiding filling $G$. By Remark~\ref{rem-sparse}, the filling $G$
can be bijectively transformed into a sparse $A'$-avoiding filling $G'$ of
$\Delta^-$, which transforms $F$ into a semi-standard $\mat{B}{0}{0}{A'}$-avoiding filling
of $\Delta$.
\end{proof}

We remark that the argument of the proof fails if the matrices
$\mat{B}{0}{0}{A}$ and $\mat{B}{0}{0}{A'}$ are replaced with $\mat{A}{0}{0}{B}$
and $\mat{A'}{0}{0}{B}$ respectively. Also, the argument fails if Ferrers shapes
are replaced with stack polyominoes.

Although Lemma~\ref{lem-swap} does not directly provide new pairs of
equivalent partition patterns, it allows us to prove the following proposition.

\begin{proposition}\label{pro-decr}
Let $s_1>s_2>\dotsb >s_m$ and $t_1>t_2>\dotsb >t_m$ be two strictly decreasing
sequences over the alphabet $[k]$, let $r_1,\dotsc, r_m$ be positive integers.
Define weakly decreasing sequences $S=s_1^{r_1} s_2^{r_2} \dotsb s_m^{r_m}$ and
$T=t_1^{r_1} t_2^{r_2} \dotsb t_m^{r_m}$. We have $M(S,k)\eqf M(T,k)$, and in
particular, if $\tau$ an arbitrary partition, then
$12\dotsb k(\tau+k)S\sim 12\dotsb k(\tau+k)T$.
\end{proposition}
\begin{proof}
We proceed by induction over minimum $j$ such that $s_i=t_i$ for each $i\le
m-j$. For $j=0$, we have $S=T$ and the result is clear. If $j>0$, assume
without loss of generality that $s_{m-j+1}-t_{m-j+1}=d>0$. Consider the
sequence $t'_1>t'_2>\dotsb>t'_m$ such that $t'_i=t_i$ for every $i\le m-j$ and
$t'_i=t_i+d$ for every $i>m-j$. The sequence $(t'_i)_{i=1}^m$ is strictly
decreasing, and its first $m-j+1$ terms are equal to $s_i$. Define
$T'=(t'_1)^{r_1} (t'_2)^{r_2} \dotsb (t'_m)^{r_m}$. By induction, $M(S,k)\eqf
M(T',k)$. To prove that $M(T,k)\eqf M(T',k)$, first write $T=T_0T_1$, where
$T_0$ is the prefix of $T$ containing all the symbols of $T$ greater than
$t_{m-j+1}$ and $T_1$ is the suffix of the remaining symbols. Notice that
$T'=T_0(T_1+d)$. We may write $M(T,k)=\mat{B}{0}{0}{A}$ and
$M(T',k)=\mat{B}{0}{0}{A'}$, where $A=M(T_1,t_{m-j}-1)$ and
$A'=M(T_1+d,t_{m-j}-1)$. By Proposition~\ref{pro-filshift}, $A\eqf A'$, and
by Lemma~\ref{lem-swap}, $M(T,k)\eqf M(T',k)$, as claimed. The last claim of
the proposition follows from Lemma~\ref{lem-fil}.
\end{proof}

\subsection{Non-crossing and non-nesting partitions}

The key application of the framework of the previous subsection is the
identity between non-crossing and non-nesting partitions. We define
non-crossing and non-nesting partitions in the following way:
\begin{definition}\label{def-crne}
A partition is \emph{$k$-noncrossing} if it avoids the pattern
$12\dotsb k12\dotsb k$, and it is \emph{$k$-nonnesting} if it
avoids the pattern $12\dotsb kk(k-1)\dotsb 1$.
\end{definition}

Let us point out that there are several different concepts of `crossings' and
`nestings' used in the literature: for example, Klazar \cite{K1} has
considered two blocks $X, Y$ of a partition to be crossing (or nesting) if
there are four elements $x_1<y_1<x_2<y_2$ (or $x_1<y_1<y_2<x_2$,
respectively) such that $x_1,x_2\in X$ and $y_1,y_2\in Y$, and similarly for
$k$-crossings and $k$-nestings. Unlike our approach, Klazar's definition
makes no assumption about the relative order of the minimal elements of $X$
and $Y$, which allows more general configurations to be considered as
crossing or nesting. Thus, Klazar's $k$-noncrossing and $k$-nonnesting
partitions are a proper subset of our $k$-noncrossing and $k$-nonnesting
partitions, (except for 2-noncrossing partitions where the two concepts
coincide).

Another approach to crossings in partitions has been pursued by Chen et al.
\cite{C1,C2}. They use the so-called \emph{linear representation}, where a
partition of $[n]$ with blocks $B_1, B_2,\dotsc, B_k$ is represented by a
graph on the vertex set $[n]$, with $a,b\in[n]$ connected by an edge if they
belong to the same block and there is no other element of this block between
them. In this terminology, a partition is $k$-crossing (or $k$-nesting) if
the representing graph contains $k$ edges which are pairwise crossing (or
nesting), where two edges $e_1=\{a<b\}$ and $e_2=\{a'<b'\}$ are crossing (or
nesting) if $a<a'<b<b'$ (or $a<a'<b'<b$ respectively). Let us call such
partitions graph-$k$-crossing and graph-$k$-nesting, to avoid confusion with
our own terminology of Definition~\ref{def-crne}. It is not difficult to see
that a partition is graph-2-noncrossing if and only if it is $2$-noncrossing,
but for nestings and for $k$-crossings with $k>2$, the two concepts are
incomparable. For instance the partition 12121 is graph-2-nonnesting but it
contains 1221, while 12112 is graph-2-nesting and avoids 1221. Similarly,
1213123 has no graph-3-crossing and contains 123123, while 1232132 has a
graph-3-crossing and avoids 123123.

Chen et al. \cite{C2} have shown that the number of graph-$k$-noncrossing and
graph-$k$-nonnesting partitions of $[n]$ is equal. Below, we prove that the
same is true for $k$-noncrossing and $k$-nonnesting partitions as well. It is
interesting to note that the proofs of both these results are based on a
reduction to theorems on pattern avoidance in the fillings of Ferrers
diagrams (this is only implicit in \cite{C2}, a direct construction is given
by Krattenthaler \cite{kra}), although the constructions employed in the
proofs of these results are quite different.

\begin{theorem}\label{thm-crne}
For every $n$ and $k$, the number of $k$-noncrossing partitions of $[n]$ is equal to
the number of $k$-nonnesting partitions of $[n]$.
\end{theorem}

By Lemma~\ref{lem-fil}, a bijection between $k$-noncrossing and
$k$-nonnesting partitions can be constructed from a bijection between
$I_k$-avoiding and $J_k$-avoiding semi-standard fillings of Ferrers
diagrams.

Krattenthaler \cite{kra} has presented a comprehensive summary of the
relationships between $I_r$-avoiding and $J_r$-avoiding fillings of a fixed
Ferrers diagram under additional constraints for row-sums and column-sums.
These relationships are based on a suitable version of the
RSK-correspondence (see \cite{ful} or \cite{S2} for a broad overview of the RSK
algorithm and related topics).

We will now state the theorem about the correspondence between $I_k$-avoiding
and $J_k$-avoiding fillings of diagrams. The result we will use is a weaker
version of Theorem~13 from \cite{kra}. Note that in the original paper, it is
not explicitly stated that the bijection between $I_k$-avoiding
and $J_k$-avoiding fillings preserves the sum of every row and every column;
however, this is an immediate consequence of the technique used in the proof.
Also, in \cite{kra}, the result is stated for arbitrary fillings with
nonnegative integers; however, the previous remark shows that the result holds
even when restricted to semi-standard fillings.

\begin{theorem}[adapted from \cite{kra}]\label{thm-kra}
For every Ferrers diagram $\Delta$ and every $k$, there is a  bijection between
the $I_k$-avoiding semi-standard fillings of $\Delta$ and the $J_k$-avoiding
semi-standard fillings of $\Delta$. The bijection preserves the number of 1-cells
in every row.
\end{theorem}

Theorem~\ref{thm-kra} and Lemma~\ref{lem-fil} give us the result we
need.
\begin{corollary}
For every $n$ and every $k$, there is a bijection between $k$-noncrossing and
$k$-nonnesting partitions of $[n]$. The bijection preserves the number of
blocks, the size of each block, and the smallest element of every block.
\end{corollary}

Applying Lemma~\ref{lem-fil} with $S=12\dotsb k$ and
$S'=k(k-1)\dotsb 1$, and translating it into the terminology of
pattern-avoiding partitions, we obtain the following result.

\begin{corollary}\label{gen-cross}
Let $\tau$ be a partition, let $k$ be an
integer. The pattern $12\dotsb k(\tau+k) 12\dotsb k$ is equivalent to
$12\dotsb k (\tau+k) k(k-1)\dotsb 1$.
\end{corollary}

Furthermore, results of Rubey, in particular~\cite[Proposition~5.3]{rub},
imply that the matrices $I_k$ and $J_k$ are in fact stack-equivalent, rather
than just Ferrers-equivalent. More precisely, Rubey's theorem deals with
fillings of moon polyominoes with prescribed row-sums. However, since a
transposed copy of a stack polyomino is a special case of a moon polyomino,
Rubey's general result applies to fillings of stack polyominoes with
prescribed column sums as well. Combining this theorem with
Lemma~\ref{lem-poly}, we obtain the following result.
\begin{corollary}
For any $k$ and any partition $\tau$, the pattern $12\dotsb k 12\dotsb
k(\tau+k)$ is equivalent to $12\dotsb k k(k-1)\dotsb 1 (\tau+k)$.
\end{corollary}

\subsection{Patterns of the form $1(\tau+1)$}

In this subsection, we will establish a general relationship between
the partitions that avoid a pattern $\tau$ and the partitions that
avoid the pattern $1(\tau+1)$. The key result is the following
theorem.

\begin{theorem}\label{thm-lift}
Let $\tau$ be an arbitrary pattern, and let $F(x)$ be its
corresponding EGF. Let $\sigma=1(\tau+1)$, and let $G(x)$ be its
EGF. For every $n\ge 1$, the following holds:
\begin{equation}\label{eq-lift}
p(n;\sigma)=\sum_{i=0}^{n-1} \binom{n-1}{i}  p(i;\tau).
\end{equation}
In terms of generating functions, this is equivalent to
\begin{equation}\label{eq-gflift}
G(x)=1+\int_0^x F(t)e^t \text{d} t.
\end{equation}
\end{theorem}

\begin{proof}
Fix $\sigma$ and $\tau$ as in the statement of the theorem. Let $\pi$ be an
arbitrary partition, and let $\pi^-$ denote the partition obtained from $\pi$ by
erasing every occurrence of the symbol 1, and decreasing every other symbol by
$1$; in other words, $\pi^-$ represents the partition obtained by removing the
first block from the partition $\pi$. Clearly, a partition $\pi$ avoids $\sigma$
if and only if $\pi^-$ avoids $\tau$. Thus, for every $\sigma$-avoiding partition
$\pi\in P(n;\sigma)$ there is a unique $\tau$-avoiding partition
$\rho\in \cup_{i=0}^{n-1} P(i;\tau)$ satisfying $\pi^-=\rho$. On the other hand,
for a fixed $\rho\in P(i;\tau)$, there are $\binom{n-1}{i}$ partitions $\pi\in
P(n;\sigma)$ such that $\pi^-=\rho$. This gives equation~\eqref{eq-lift}.

To get equation~\eqref{eq-gflift}, we multiply both sides of~\eqref{eq-lift} by
$\frac{x^n}{n!}$ and sum for all $n\ge 1$. This yields

\begin{align*}
G(x)-1&= \sum_{n\ge1}\frac{x^n}{n!} \sum_{i=0}^{n-1}\binom{n-1}{i}
p(i;\tau) =\int_0^x \sum_{n\ge1}\frac{t^{n-1}}{(n-1)!}
\sum_{i=0}^{n-1}\binom{n-1}{i}
p(i;\tau)\text{d} t\\
&=\int_0^x \sum_{n\ge 0}\frac{t^n}{n!} \sum_{i=0}^{n}\binom{n}{i}
p(i;\tau)\text{d}t=\int_0^x \sum_{n\ge
0}\sum_{i=0}^{n}\frac{t^i}{i!}
p(i;\tau)\frac{t^{n-i}}{(n-i)!}\text{d} t\\
&=\int_0^x \left(\sum_{i\ge 0}\frac{t^i}{i!}
p(i;\tau)\right)\left(\sum_{k\ge 0}\frac{t^k}{k!}\right)\text{d}
t=\int_0^x F(t)e^t \text{d}t,
\end{align*}
which is equivalent to equation~\eqref{eq-gflift}.
\end{proof}

The following result is an immediate consequence of
Theorem~\ref{thm-lift}.

\begin{corollary}\label{cor-lift}
If $\tau\sim\tau'$ then $1(\tau+1)\sim 1(\tau'+1)$, and more generally, $12\dotsb
k(\tau+k)\sim 12\dotsb k(\tau'+k)$. In particular, since $123\sim122\sim112\sim121$,
we see that for every $m\ge 2$ the patterns $12\dotsb (m-1)m(m+1)$, $12\dotsb
(m-1)mm$, $12\dotsb (m-1)(m-1)m$ and $12\dotsb (m-1)m(m-1)$ are equivalent.
Conversely, if $1(\tau+1)\sim 1(\tau'+1)$, then $\tau\sim\tau'$.
\end{corollary}
\begin{proof}
To prove the last claim, notice that equation~\eqref{eq-lift} can be inverted
to obtain
\[
p(n-1;\tau)=\sum_{i=0}^{n-1} (-1)^i\binom{n-1}{i} p(n-i;\sigma).
\]
The other claims follow directly from Theorem~\ref{thm-lift}.
\end{proof}

\subsection{Patterns equivalent to $12\dotsb m(m+1)$}

The partitions that avoid $12\dotsb m(m+1)$, or equivalently, the partitions
with at most $m$ blocks, are a very natural pattern-avoiding class of
partitions. Their number $p(n;12\dotsb (m+1))$ is equal to $\sum_{i=0}^{m}
S(n,i)$, where $S(n,i)$ is the Stirling number of the second kind, which is
equal to the number of partitions of $[n]$ with exactly $i$ blocks.

As an application of the previous results, we will now present two classes of
patterns that are equivalent to the pattern $12\dotsb (m+1)$. From this
result, we obtain an alternative combinatorial interpretation of the Stirling
numbers $S(n,i)$.

Our result is summarized in the following theorem.

\begin{theorem}\label{thm-stirling}
For every $m\ge 2$, the following patterns are equivalent:
\begin{itemize}
\item[(a)] the pattern $12\dotsb(m-1)m(m+1)$,
\item[(b)] the patterns $12\dotsb(m-1)md$, where $d$ is any number from the set
$[m]$,
\item[(c)] the patterns $12\dotsb(m-1)dm$, where $d$ is any number from the set
$[m-1]$.
\end{itemize}
\end{theorem}
\begin{proof}
From Corollary~\ref{cor-lift}, we get the following:
\[
12\dotsb m(m+1)\sim 12\dotsb(m-1)mm\sim 12\dotsb(m-1)(m-1)m.
\]
The equivalences between $$12\dotsb(m-1)mm\sim 12\dotsb(m-1)md\mbox{
and }12\dotsb(m-1)(m-1)m\sim 12\dotsb(m-1)dm$$ are obtained by a
repeated application of Proposition~\ref{pro-filshift}.
\end{proof}

\subsection{Binary patterns}
Let us now focus on the avoidance of \emph{binary} patterns, i.e.,
the patterns that only contain the symbols 1 and 2.

We will first consider the forbidden patterns of the form $1^k21^l$.
We have already seen that $112\sim 121$. The following theorem
offers a generalization.

\begin{theorem}
\label{thm-12} For any three integers $j,k,m$ satisfying $1\le j,k\le m$, the
pattern $1^j21^{m-j}$ is equivalent to the pattern $1^k21^{m-k}$.
\end{theorem}

Before we present the proof of Theorem~\ref{thm-12}, we need some
preparation. Let $\pi=\pi_1 \pi_2 \dotsb\pi_n$ be a partition.
Clearly, $\pi$ can be uniquely expressed as $1P_1 1P_2 1\dotsb
1P_{p-1}1P_p$, where the $P_i$ are (possibly empty) maximal
contiguous subsequences of $\pi$ that do not contain the symbol~$1$.
The sequence $P_i$ will be referred to as \emph{the $i$-th chunk of
$\pi$}. By concatenating the chunks into a sequence $P=P_1\dotsb
P_p$ and then subtracting 1 from every symbol of $P$, we obtain a
canonical sequence of a partition; let this partition be denoted by
$\pi^-$. The key ingredient in the proof of Theorem~\ref{thm-12} is
the following lemma.
\begin{lemma}
\label{lem-12} Let $\pi$ be a partition that has
$p$ occurrences of the symbol 1, let $P_i$ and $\pi^-$ be as above. Let $r\ge
1$ and $s\ge 0$ be two integers. The partition $\pi$ avoids $1^r21^s$ if and
only if the following two conditions hold: \begin{itemize}
\item $\pi^-$ avoids $1^r21^s$.
\item For every $i$ such that $r\le i\le p-s$, the chunk $P_i$ is empty.
\end{itemize}
\end{lemma}
\begin{proof}
Clearly, the two conditions are necessary. To see that they are sufficient, we
argue by contradiction: let $\pi$ be a partition that satisfies the
two conditions, and assume that $\pi$ has a subsequence $a^rba^s$ for two
symbols $a<b$. If $a=1$ we have a contradiction with the second condition, and
if $a>1$, then $\pi^-$ contains the sequence $(a-1)^r(b-1)(a-1)^s$,
contradicting the first condition.
\end{proof}

We are now ready prove Theorem~\ref{thm-12}.

\begin{proof}
[Proof of Theorem~\ref{thm-12}] It is enough to prove that for every $k\ge 1$
and every $m>k$ there is a bijection $f$ from $P(n;1^k21^{m-k})$ to
$P(n;1^{m}2)$. To define $f$, we will proceed by induction on the number of
blocks of $\pi$. If $\pi =1^n$ then we define $f(\pi)=\pi$. Assume that $f$
has been defined for all partitions with less than $b$ blocks, and let $\pi\in
P(n;1^k21^{m-k})$ be a partition with $b$ blocks, let $p$ be the size of the
first block of $\pi$. Let $P_1,\dotsc,P_p$ be the chunks of $\pi$ and let
$\pi^-$ be defined as above. Define $\overline\sigma=f(\pi^-)$;  this is well
defined, since $\pi^-\in P(n-p;1^k21^{m-k})$ and $\pi^-$ has $b-1$ blocks.
Let $S=\overline\sigma+1$. We express $S$ as a concatenation of the form
$S=S_1S_2\dotsb S_p$, where the length of $S_i$ is equal to the length of
$P_i$. By Lemma~\ref{lem-12}, the chunk $P_i$ (and hence also $S_i$) is empty
whenever $k\le i\le p-m+k$. We put $f(\pi)=\sigma$, where $\sigma$ is defined
as follows:
\begin{itemize}
\item If $p<m$, then $\sigma=1S_1 1S_2 1\dotsb 1S_{p-1}1S_p$.
\item If $p\ge m$, then $\sigma=1S_1 1S_2 1\dotsb 1
S_{k-1}1S_{p-m+k+1}1S_{p-m+k+2}1\dotsb 1 S_{p-1}1S_p 1^{p-m+1}$.
\end{itemize}
Using Lemma~\ref{lem-12}, we may easily see that $\sigma$ avoids $1^m2$. It is
also straightforward to check that $f$ is indeed a bijection from $P(n;1^k21^{m-k})$ to
$P(n;1^m2)$. Note that $f$ preserves not only the number of blocks of the
partition, but also the size of each block.
\end{proof}

Using our results on fillings, we can add another pattern to the
equivalence class covered by Theorem~\ref{thm-12}.

\begin{theorem}\label{thm-12m}
For every $m\ge 1$, the pattern $12^m$ is equivalent to the pattern
$121^{m-1}$.
\end{theorem}
\begin{proof}
This is just Proposition~\ref{pro-filshift} with $k=2$ and $S=1^{m-1}$.
\end{proof}

\begin{corollary}\label{co1222}
Let $m$ be a positive integer, let $\tau$ be any pattern from the set
\[T=\{1^k21^{m-k};\; 1\le
k\le m\}\cup\{12^m\}.\] The EGF $F(x)$ of a pattern $\tau\in T$ is
given by
\[
F(x)=1+\int_0^x \exp\left(t+\sum_{i=1}^{m-1}\frac{t^i}{i!}\right) \text{d}t.
\]
\end{corollary}
\begin{proof}
We have seen that all the patterns from the set $T$ are equivalent,
so let us pick $\tau=12^m$. The EGF follows from
equation~\eqref{eq-1m} on page~\pageref{eq-1m} and from Theorem~\ref{thm-lift}.
\end{proof}

We now turn to another type of binary patterns, namely the patterns of the
form $12^k12^{m-k}$ with $1\le k\le m$. It turns out that these patterns are
also equivalent, provided $m$ is fixed. In fact, we will prove several
generalizations of this fact.

Our argument will again be based on the analysis of fillings of
stack polyominoes. However, to make full use of this approach, we
will establish a stronger relation than mere stack equivalence. For
this, we need the following definition.

\begin{definition}
Let $F$ be a semi-standard filling of a stack polyomino $\Pi$ and let $t\ge 1$ be an
integer. We say that $F$ is \emph{$t$-falling} if its first $t$ rows all contain at
least one 1-cell, and the leftmost 1-cells of these rows form a decreasing
chain; formally, $F$ is $t$-falling if for every $i<j\le t$ the leftmost
1-cell in row $i$ exists and appears to the right of the leftmost 1-cell in row
$j$, which must exist as well.
\end{definition}
Notice that a $t$-falling filling of $\Pi$ only exists if the leftmost column
of $\Pi$ intersects its first $t$ rows.

In the rest of this subsection, $\spq p q$ denotes the sequence $2^p 1 2^q$
and $\bspq p q$ denotes the sequence $1^p21^q$, where $p,q$ are
nonnegative integers.

\begin{lemma}\label{lem-fall}
For every $p,q\ge 0$, the matrix $M(\spq p q,2)$ is stack-equivalent to the
matrix $M(\spq {p+q} 0,2)$. Furthermore, if $p\ge 1$, then for every stack polyomino
$\Pi$, there is a bijection $f$ between the $M(\spq p q,2)$-avoiding and
$M(\spq{p+q}0,2)$-avoiding semi-standard fillings of $\Pi$ with these
two properties:
\begin{itemize}
\item $f$ preserves the number of 1-cells in every row.
\item Both $f$ and $f^{-1}$ map $t$-falling fillings to
$t$-falling fillings, for every $t\ge 1$.
\end{itemize}
\end{lemma}

\begin{proof}
Let $M=M(\spq p q,2)$ and $M'=M(\spq {p+q}0,2)$, for some $p,q\ge 0$.  We will
proceed by induction over the number of rows of $\Pi$.  If $\Pi$ has only one
row, then a constant mapping is the required bijection.  Assume now that $\Pi$
has $r\ge 2$ rows, and assume that we are presented with a semi-standard filling
$F$ of~$\Pi$. Let $\Pi^-$ be the diagram obtained from $\Pi$ by erasing the
$r$-th row as well as every column that contains a $1$-cell of $F$ in the
$r$-th row. The filling $F$ induces on $\Pi^-$ a semi-standard filling~$F^-$.

We claim that for every $p,q\ge 0$, a filling $F$ avoids $M$ if and only if
these two conditions are satisfied:
\begin{itemize}
\item[(a)] The filling $F^-$ avoids $M$.
\item[(b)] If the $r$-th row of $F$ contains $m$ $1$-cells in columns
$c_1<c_2<\dotsb<c_m$ and if $m\ge p+q$, then for every $i$ such that $p\le i\le
m-q$, the column $c_i$ is either the rightmost column of the last row of $\Pi$, or it is
directly adjacent to the column $c_{i+1}$ (i.e. $c_i+1=c_{i+1}$).
\end{itemize}
Clearly, the two conditions are necessary. To see that they are sufficient,
note that the first condition guarantees that $F$ does not contain any copy
of $M$ that would be confined to the first $r-1$ rows, whereas the
second condition guarantees that $F$ has no copy of $M$ that would
intersect the $r$-th row.

We now define recursively the required bijection between
$M$-avoiding and $M'$-avoiding fillings. Let $F$ be an
$M$-avoiding filling of $\Pi$, let $F^-$ and $c_1,\dotsc,c_m$ be as
above. By the induction hypothesis, we already have a bijection between
$M$-avoiding and $M'$-avoiding fillings of the shape $\Pi^-$.
This bijection maps $F^-$ to a filling $\tilde F^-$ of $\Pi^-$. Let $\tilde F$
be the filling of $\Pi$ that has the same values as $F$ in the $r$-th row, and
the columns not containing a $1$-cell in the $r$-th row are filled according
to~$\tilde F^-$. Note that $\tilde F$ contains no copy of $M'$ in
its first $r-1$ rows and it contains no copy of $M$ that would
intersect the last row.

If $\tilde F$ has less than $p+q$ 1-cells in the last row, we define $f(F)=\tilde
F$, otherwise we modify $\tilde F$ as follows: for every $i=1,\dotsc,q$, we
consider the columns with indices strictly between $c_{m-q+i}$ and
$c_{m-q+i+1}$ (if $i=q$, we take all columns to the right of $c_m$ that
intersect the last row), we remove these columns from $\tilde F$ and
re-insert them between the columns $c_{p+i-1}$ and $c_{p+i}$ (which used to
be adjacent by condition (b) above). Note that these transformations preserve
the relative left-to-right order of all the columns that do not contain a $1$-cell in their
$r$-th row; in particular, the resulting filling still has no copy of $M'$ in
the first $r-1$ rows. By construction, the filling also satisfies condition
(b) for the values $p'=p+q$ and $q'=0$ used instead of the original $p$ and
$q$; and in particular, it is a $M'$-avoiding filling. This construction
provides a bijection $f$ between $M$-avoiding and $M'$-avoiding fillings.

It is clear that $f$ preserves the number of 1-cells in each row. It remains
to check that if $p\ge 1$, then $f$ preserves the $t$-falling property. Let
us fix $t$, and let $r$ be the number of rows of $\Pi$. If $r<t$ then no
filling of $\Pi$ is $t$-falling. If $r>t$, then $F$ is a $t$-falling filling
if and only if $F^-$ is $t$-falling, so we obtain the required result from
the induction hypothesis and from the fact that the process that transforms
the intermediate filling $\tilde F$ into the final filling $f(F)$ does not
change the relative position of the 1-cells of the first $r-1$ rows.  Finally,
if $r=t$, then $F$ is $t$-falling if and only if $F^-$ is $(t-1)$-falling and
the leftmost 1-cell of the $r$-th row of $F$ is in the leftmost column of
$\Pi$. Both these conditions are preserved by $f$ and $f^{-1}$, provided $p\ge
1$.
\end{proof}

With the help of Lemma~\ref{lem-fall}, we are able to prove several results
about pattern avoidance in partitions. We first state a direct corollary of
Lemmas~\ref{lem-fil}, \ref{lem-poly}, and~\ref{lem-fall}, and
Proposition~\ref{pro-filshift}.

\begin{corollary}\label{cor-spq}
For any partition $\tau$, for any $k\ge 2$, and for any $p,q\ge 0$, the pattern
$12\dotsb k(\tau+k) \spq p q$ is equivalent to $12\dotsb k(\tau+k) \spq{p+q}0$
$12\dotsb k \spq p q (\tau+k)$ is equivalent to $12\dotsb k \spq{p+q}0(\tau+k)$.
\end{corollary}

Next, we present two theorems that make use of the $t$-falling property.
Recall that $\bspq p q =1^p21^q$.

\begin{theorem}\label{thm-spq1}
Let $\tau$ be any partition with $k$ blocks, let $p\ge 1$ and $q\ge 0$. The
pattern $\sigma=\tau (\bspq p q+k)$ is equivalent to $\sigma'=\tau (\bspq{p+q}0+k)$.
\end{theorem}
\begin{proof}
Let $\pi$ be a partition of $[n]$ with $m$ blocks, let $M=M(\pi,m)$. We color
the cells of $M$ red and green, where a cell in row $i$ and column $j$ is
green if and only if the submatrix of $M$ formed by the intersection of the
first $i-1$ rows and $j-1$ columns of $M$ contains $M(\tau,k)$. Let $\Gamma$
be the diagram formed by the green cells of $M$, and let $G$ be the filling
of $\Gamma$ by the values from $M$. Note that $\Gamma$ is an upside-down copy
of a Ferrers shape. It is easy to see that the partition $\pi$ avoids
$\sigma$ if and only if $G$ avoids $M(\bspq p q,2)$, and $\pi$ avoids
$\sigma'$ if and only if $G$ avoids $M(\bspq{p+q}0,2)$.

Let us now assume that $\pi$ is $\sigma$-avoiding. We transform $\pi$
into a $\sigma'$-avoiding partition $\pi'$ by the following
procedure. We first turn the filling $G$ and the diagram $\Gamma$
upside down, which transforms $\Gamma$ into a Ferrers shape
$\overline \Gamma$, and it also transforms the $M(\bspq p
q,2)$-avoiding filling $G$ into an $M(\spq p q,2)$-avoiding
filling $\overline G$ of $\overline \Gamma$. We may then apply the
bijection $f$ of Lemma~\ref{lem-fall} to $\overline G$, ignoring the
zero columns of $\overline G$. We then
turn the resulting filling $\overline G'=f(\overline G)$ of
$\overline \Gamma$ upside down again to obtain an
$M(\bspq{p+q}0,2)$-avoiding filling $G'$ of $\Gamma$. We then fill
the green cells of $M$ with the values of $G'$ while the filling of
the red cells remains the same. We thus obtain a matrix $M'$. The
matrix $M'$ has exactly one 1-cell in each column, so there is a
sequence $\pi'$ over the alphabet $[m]$ such that $M'=M(\pi',m)$.
Note that a cell is green with respect to $M$ if and only if it is
green with respect to~$M'$.

By construction, the sequence $\pi'$ has no subsequence order-isomorphic to
$\sigma'$. We now need to show that $\pi'$ is a restricted-growth sequence.
For this, we will use the preservation of the $t$-falling property. Let $c_i$
be the leftmost 1-cell of the $i$-th row of $M$, let $c'_i$ be the leftmost
1-cell of the $i$-th row of $M'$. Let $s$ be the largest index such that the
cell $c_s$ is red, we set $s=0$ if no such cell exists. Note that all the
cells $c_1,\dotsc,c_s$ are red and all the cells $c_{s+1},\dotsc,c_m$ are
green. We have $c_i=c'_i$ for every $i\le s$. If $s>0$, we also see that all
the green 1-cells are in the columns to the right of $c_s$.

It is easy to see that all the 1-cells above row $s$ are green; in
particular, the filling $\overline G$ is $t$-falling for some $t\ge m-s$. By
Lemma~\ref{lem-fall}, the filling $\overline G'$ is $t$-falling as well. It
follows that the cells $c'_{s+1},\dotsc,c'_m$ form a left-to-right increasing
chain, and since all these cells are to the right of $c'_s$, we see that all
the cells $c'_1,\dotsc,c'_m$ form a left-to-right increasing sequence, hence
$\pi'$ is in canonical sequential form, i.e., $\pi'$ is a partition from
$P(n;\sigma')$.

It is obvious that the above construction can be reversed, which shows that it
is indeed a bijection between $P(n;\sigma)$ and $P(n;\sigma')$.
\end{proof}

The following result is proved by a similar approach, but the argument is
slightly more technical.

\begin{theorem}\label{thm-spq2}
Let $T$ be an arbitrary sequence over the alphabet $[k]$, let $p\ge 1$ and
$q\ge0$. The partition $\sigma=12\dotsb k(\bspq p q +k)T$ is equivalent to
$\sigma'=12\dotsb k (\bspq {p+q}0+k) T$.
\end{theorem}

\begin{proof}
Let $\pi$ be a partition of $[n]$ with $m$ blocks, let $M=M(\pi,m)$. As in
the previous proof, we color the cells of $M$ red and green. A cell in row
$i$ and column $j$ will be green if the submatrix of $M$ formed by rows
$1,\dotsc,i-1$ and columns $j+1,\dotsc,n$ contains $M(T,k)$.

Let $\Gamma$ be
the diagram formed by the green cells and $G$ its filling inherited from $M$.
The partition $\pi$ contains $\sigma$ (or $\sigma'$) if and only if $G$
contains $M(\bspq p q,2)$ (or $M(\bspq{p+q}0,2)$, respectively). The diagram
$\Gamma$ is an upside-down copy of a left-justified stack polyomino. Let
$\Gamma^-$ be the diagram obtained by erasing the zero rows and zero columns of
$\Gamma$, and let $G^-$ be the corresponding filling of $\Gamma^-$.
The upside-down copy of the filling $G^-$ is $r$-falling, where $r$ is the
number of rows of $\Gamma^-$.

We apply the same construction as in the
previous proof and transform $G^-$ into an $M(\bspq{p+q}0,2)$-avoiding filling
$G'$ of $\Gamma^-$ and obtain a matrix $M'=M(\pi',n)$ that avoids
$M(\sigma',k+2)$. It remains to argue that $\pi'$ is in canonical sequential
form.

Let $c_i$ (or $c'_i$) be the leftmost 1-cell in row $i$ of $M$ (or $M'$,
respectively). To prove that $\pi'$ is a partition, we want to show that
$c'_1,\dotsc,c'_m$ form a left-to-right increasing sequence in $M'$. Let us
now fix two row indices $i<j$. We claim that $c'_i$ is left of $c'_j$.

If both $c'_i$ and $c'_j$ are green, then the claim follows from the preservation of
the $r$-falling property. If $c'_j$ is red, then all the green cells below
row $j$ (including $c'_i$) are left of $c_j$.

Finally, assume that $c'_j$ is green and $c'_i$ is red. We have $c'_i=c_i$. In the filling $G$, all the
green 1-cells that are to the left of $c_i$ are also below row $i$; let $x$ be the
number of such 1-cells. Then $x$ corresponds to the number of nonzero columns of
$G$ that are to the left of $c_j$. Since the number of these nonzero columns is
preserved by the mapping $f$, we see that $G'$ also has $x$ 1-cells left of $c_i$.

Since $f$ preserves the number of 1-cells in each row, both $G$ and $G'$ have exactly $x$
1-cells below row $i$. All the 1-cells of $G'$ below row $i$ must appear to
the left of $c_i$, and since there are only $x$ 1-cells of $G'$ to the left
of $c_i$, they must all appear below row $i$, and in particular, all the
green 1-cells above row $i$ (including the cell $c'_j$) appear to the right
of $c_i$.
\end{proof}

\subsection{Patterns equivalent to $12^k13$}
Throughout this subsection, we will assume that $t$ is an arbitrary fixed
integer, and we will deal with the following sets of patterns:
\begin{align*}
\Sigma^+&=\{12^{p+1}12^q32^r; p,q,r\ge 0, p+q+r=t\}\\
\Sigma^-&=\{12^{p+1}32^q12^r; p,q,r\ge 0, p+q+r=t\}\\
\Sigma&=\Sigma^+\cup \Sigma^-
\end{align*}
Our aim is to show that all the patterns in $\Sigma$ are equivalent.
We will use the following definition.

\begin{definition}
Let $\sigma$ be a pattern over the alphabet $\{1,2,3\}$, let $\pi$ be a
partition with $m$ blocks, and let $k\le m$ be an integer. We say that
\emph{$\pi$ contains $\sigma$ at level $k$}, if there are symbols $l,h\in[m]$
such that $l<k<h$, and the partition $\pi$ contains a subsequence $S$ made of the
symbols $\{l,k,h\}$ which is order-isomorphic to $\sigma$.
\end{definition}

For example, the partition $\pi=1231323142221$ contains $\sigma=121223$ at level 3, because
$\pi$ contains the subsequence $131334$, but $\pi$ avoids $\sigma$ at level 2,
because $\pi$ has no subsequence of the form $l2l22h$ with $l<2<h$.

Our plan is to show, for suitable pairs $\sigma, \sigma'\in\Sigma$, that for
every $k$ there is a bijection $f_k$ that maps the partitions avoiding
$\sigma$ at level $k$ onto the partitions avoiding $\sigma'$ at level $k$, while
preserving $\sigma$-avoidance and $\sigma'$-avoidance at all the other levels.
Composing the maps $f_k$ for all possible $k$, we will then obtain a bijection
between $P(n;\sigma)$ and $P(n;\sigma')$.

We need more definitions.
\begin{definition}
Consider a partition $\pi$, and fix a level $k\ge 2$. A symbol
of $\pi$ is called \emph{$k$-low} if it is smaller than $k$ and \emph{$k$-high} if it
is greater than $k$. A \emph{$k$-low cluster} (or \emph{$k$-high cluster}) is a
maximal consecutive sequence of $k$-low symbols (or $k$-high symbols, respectively) in
$\pi$. The \emph{$k$-landscape} of $\pi$ is a word over the alphabet
$\{\Lc,k,\Hc\}$ obtained from $\pi$ by replacing each $k$-low cluster with a
single symbol $\Lc$ and each $k$-high cluster with a single symbol $\Hc$.

A word $w$ over the alphabet $\{\Lc,k,\Hc\}$ is called \emph{$k$-landscape word}
if it satisfies the following conditions:
\begin{itemize}
\item The first symbol of $w$ is $\Lc$, the second symbol of $w$ is $k$.
\item No two symbols $\Lc$ are consecutive in $w$, no two symbols $\Hc$ are
consecutive in $w$.
\end{itemize}
Clearly, the landscape of a partition is a landscape word.

Two $k$-landscape words $w$ and $w'$ are said to be \emph{compatible}, if each of
the three symbols $\{\Lc,k,\Hc\}$ has the same number of occurrences in $w$ as
in $w'$.
\end{definition}

We will often drop the prefix $k$ from these terms, if the value of $k$ is
clear from the context.

To give an example, consider $\pi=1231323142221$: it has five $3$-low
clusters, namely $12$, $1$, $2$, $1$ and $2221$, it has one $3$-high cluster $4$,
and its 3-landscape is $\Lc 3\Lc 3\Lc 3\Lc\Hc \Lc$.

If $w$ and $w'$ are two compatible $k$-landscape words, we have a natural
bijection between partitions with landscape $w$ and partitions with
landscape~$w'$: if $\pi$ has landscape $w$, we map $\pi$ to the partition
$\pi'$ of landscape $w'$ which has the same $k$-low clusters and $k$-high
clusters as $\pi$, and moreover, the $k$-low clusters appear in the same
order in $\pi$ as in $\pi'$, and also the $k$-high clusters appear in the
same order in $\pi$ as in $\pi'$. It is not difficult to check that these
rules define a unique sequence $\pi'$ and this sequence is indeed a partition.
This provides a bijection between partitions of landscape $w$ and partitions
of landscape $w'$ which will be called \emph{the $k$-shuffle from $w$
to $w'$}.

The key property of shuffles is established by the next lemma.

\begin{lemma}\label{lem-shuffle}
Let $w$ and $w'$ be two compatible $k$-landscape words. Let $\pi$ be a
partition with $k$-landscape $w$, let $\sigma$ be a pattern from $\Sigma$, let
$\pi'$ be the partition obtained from $\pi$ by the shuffle from $w$ to $w'$, let $j$
be an integer. The following holds:
\begin{enumerate}
\item If $\sigma$ does not end with the symbol 1 and $j>k$, then $\pi'$
contains $\sigma$ at level $j$ if and only if $\pi$ contains $\sigma$ at level
$j$.
\item If $\sigma$ does not end with the symbol 3 and $j<k$, then $\pi'$
contains $\sigma$ at level $j$ if and only if $\pi$ contains $\sigma$ at level
$j$.
\end{enumerate}
\end{lemma}
\begin{proof}
We begin with the first claim of the lemma. Let us choose
$\sigma\in\Sigma^-$ (the case of $\sigma\in\Sigma^+$ is analogous)
and let us fix $j>k$. Let us write $\sigma=12^{p+1}32^q12^r$, with
$r>0$ by the assumption of the lemma. Assume that $\pi$ contains
$\sigma$ at level $j$. In particular, $\pi$ has a subsequence
$S=lj^{p+1}hj^qlj^r$, with $l<j<h$.

We distinguish two cases: first,
if $k<l$, then all the symbols of $S$ are $k$-high. Since the shuffle
preserves the relative order of high symbols, $\pi'$ contains the
subsequence $S$ as well. Second, if $l\le k$, then the shuffle
preserves the relative order of the symbols $j$ and $h$, which are
all high. Let $x$ and $y$ be the two symbols of $S$ directly
adjacent to the second occurrence of $l$ in $S$ (if $q>0$, both
these symbols are equal to $j$, otherwise one of them is equal to
$h$ and the other to $j$). The two symbols are both high, but they
must appear in different $k$-high clusters. After the shuffle, the
two symbols $x$ and $y$ will again be in different clusters,
separated by a non-high symbol $l'\le k$, and since the first
occurrence of $l'$ in $\pi'$ precedes any occurrence of $j$, the
partition $\pi'$ will contain a subsequence $l'j^{p+1}hj^ql'j^r$,
which is order-isomorphic to $\sigma$.

We see that the shuffle preserves the occurrence of $\sigma$ at
level $j$. Since the inverse of the shuffle from $w$ to $w'$ is the shuffle
from $w'$ to $w$, we see that the inverse of a shuffle preserves the
occurrence of $\sigma$ at level $j$ as well.

The second claim of the lemma is proved by an analogous argument.
\end{proof}

We will use shuffles as basic building blocks for our bijection. The
first example is the following lemma.

\begin{lemma}\label{lem-sigma-}
For every $p,q,r\ge 0$, the pattern $\sigma=12^{p+1}12^q32^r$ is equivalent to the
pattern $\sigma'=12^{p+1}32^q12^r$.
\end{lemma}
\begin{proof}
Let us fix $p,q,r\ge 0$ and define $t=p+q+r$. For a given $k$,  a partition $\pi$
of $[n]$ is called a \emph{$k$-hybrid} if $\pi$ avoids $\sigma'$ at every level
$j< k$ and $\pi$ avoids $\sigma$ at every level $j\ge k$. We will show that
for every $k\in \{2,\dotsc,n-1\}$ there is a bijection $f_k$ between $k$-hybrids and
$(k+1)$-hybrids. Since 2-hybrids are precisely the $\sigma$-avoiding
partitions of $[n]$ and $n$-hybrids are precisely the $\sigma'$-avoiding
partitions of $[n]$, this gives the required result.

Let us fix $k$. Note that a partition $\pi$ contains $\sigma$ at level $k$
if and only if its $k$-landscape $w$ contains a subsequence $k^{p+1}\Lc k^q
\Hc k^r$. Similarly, $\pi$ contains $\sigma$ at level $k$ if and only if $w$
contains a subsequence $k^{p+1}\Hc k^q \Lc k^r$.

Let $\pi$ be a $k$-hybrid with landscape $w$. If $\pi$ has less than $t+1$
occurrences of $k$, then it is also a $(k+1)$-hybrid and we put
$f_k(\pi)=\pi$. Otherwise, we write $w=xyz$, where $x$ is the shortest prefix of $w$
that has $p+1$ symbols $k$ and $z$ is the shortest suffix of $w$ that has $r$
symbols $k$. By assumption, $x$ and $z$ do not overlap (although they may be
adjacent if $q=0$). Let $\overline y$ be the word obtained by reversing the order
of the letters of $y$, define $w'=x\overline y z$. Note that $w'$ is a landscape word
compatible with $w$, and that $w$ avoids $k^{p+1}\Lc k^q\Hc
k^r$ if and only if $w'$ avoids $k^{p+1}\Hc k^q\Lc k^r$. We apply to $\pi$ the
shuffle from $w$ to $w'$ which transforms it into a partition $\pi'=f_k(\pi)$.

Using Lemma~\ref{lem-shuffle}, it is easy to check that $\pi'$ is a
$(k+1)$-hybrid, which shows that $f_k$ is the required bijection.
\end{proof}

Another result in the same spirit is the following lemma.

\begin{lemma}\label{lem-sigma+}
For every $p,q,r\ge 0$, the pattern $\sigma=12^{p+2}1 2^q 3 2^r$ is equivalent to
the pattern $\sigma'=12^{p+1}1 2^q 3 2^{r+1}$.
\end{lemma}
\begin{proof}
We follow the same argument as in Lemma~\ref{lem-sigma-}. As before, a
\emph{$k$-hybrid} is a partition that avoids $\sigma'$ at every level $j<k$ and
that avoids $\sigma$ at every level $j\ge k$. We will present a bijection $f_k$
between $k$-hybrids and $(k+1)$-hybrids.

Note that $\pi$ avoids $\sigma$ at level $k$
if and only if its landscape $w$ avoids $k^{p+2}\Lc k^q \Hc k^r$.

Fix a $k$-hybrid $\pi$ with a landscape $w$. If $\pi$
has less than $p+2+q+r$ occurrences of $k$, then it is also a $(k+1)$-hybrid and
we define $f_k(\pi)=\pi$; otherwise, we write $w=xSyz$ where $x$ is the
shortest prefix of $w$ that has $p+1$ occurrences of $k$, $z$ is the shortest
suffix with $r$ occurrences of $k$, $S$ is the subword that starts just after
the $(p+1)$-th occurrence of $k$ and ends immediately after the $(p+2)$-th
occurrence of $k$. We define $w'=xy\overline Sz$, where $\overline S$ is the reversal
of $S$.

Note that in the definition of $w'$, we need to take $w'=xy\overline Sz$ instead
of the seemingly more natural definition $w=xySz$. This is because in general, the string
$xySz$ need not be a landscape word, since it may
contain to consecutive occurrences of either $\Lc$ or $\Hc$. Our definition
guarantees that $w'$ is a correct landscape word, and that $w'$ avoids $k^{p+1}\Lc k^q
\Hc k^{r+1}$ if and only if $w$ avoids $k^{p+2}\Lc k^q \Hc k^r$ (which is if and
only if $y$ avoids $\Lc k^q\Hc$).

The rest of the argument is the same as in the previous lemma.
\end{proof}

We are now ready to state and prove the main result of this subsection.

\begin{theorem}
All the patterns in the set $\Sigma$ are equivalent.
\end{theorem}
\begin{proof}
By Corollary~\ref{cor-spq}, we already know that for any $p,q\ge 0$, the pattern
$12^{p+1}12^{q}3$ is equivalent to the pattern $12^{p+q+1}13$. This, together
with the two previous lemmas gives the required result by transitivity.
\end{proof}

\subsection{More `landscape' patterns}

We will show that with a little bit of additional effort, the previous
argument involving landscapes can be adapted to prove, for every
$p,q\ge 0$, the following equivalences:
\begin{itemize}
\item $1232^{p}142^{q}\sim 12312^p 42^q$
\item $1232^p 412^q \sim 1232^p42^q 1$
\item $123^{p+1}413^q\sim 12343^p13^q$
\item $123^{p+1}143^q \sim 123^{p+1}13^q 4$
\end{itemize}

Throughout this subsection, we will say that $\tau$ is a \emph{1-2-4 pattern}
if $\tau$ has the form $123S$ where $S$ is a sequence that has exactly one
occurrence of the symbol 1, exactly one occurrence of the symbol 4, and all
its remaining symbols are equal to 2, and furthermore, the symbol 4 is neither
the first nor the last symbol of $S$. Similarly, a \emph{1-3-4 pattern} is a
pattern of the form $123S$ where $S$ has one occurrence of 1 and of 4, and all
its other symbols are equal to 3, and furthermore, the symbol 1 is not the last symbol of $S$.

We have decided to exclude the patterns of the form $1232^p12^q4$,
$12342^p12^q$ and $123 3^p43^q1$ from the set of 1-2-4 and 1-3-4
patterns defined above, because some of the arguments we will need
in the following discussion (namely in Lemma~\ref{lem-shuffle2})
would become more complicated it these special types of patterns
were allowed. We need not be too concerned about this constraint,
because we have already dealt with the patterns of the three
excluded types in Corollary~\ref{cor-spq} and
Theorem~\ref{thm-spq2}.

For our arguments, we need to extend some of the terminology of the previous
subsection to cover the new family of patterns.
Let $\tau$ be a 1-2-4 pattern, let $k$ be a natural number, and let $\pi$ be
a partition. We say that $\pi$ \emph{contains $\tau$ at level $k$}, if $\pi$
has a subsequence $T$ order-isomorphic to $\tau$ such that the occurrences of
the symbol 2 in $\tau$ correspond to the occurrences of the symbol $k$ in $T$.
Similarly, if $\tau$ is a 1-3-4 pattern, we say that a partition $\pi$
\emph{contains $\tau$ at level $k$} if $\pi$ has a subsequence $T$
order-isomorphic to $\tau$ with the symbol $k$ in $T$ corresponding to the
symbol $3$ in $\tau$.

Our aim is to prove an analogue of Lemma~\ref{lem-shuffle} for 1-2-4 and 1-3-4
patterns. Unfortunately, general $k$-shuffles may behave badly with respect to
the avoidance of these patterns. However, we will define special types of
$k$-shuffles that have the properties we need. We first introduce some new
definitions.

Let $w$ be a $k$-landscape word. We say that two occurrences of the symbol
$\Hc$ in $w$ are \emph{separated} if there is at least one occurrence of $\Lc$
between them. Similarly, two symbols $\Lc$ are separated if there is at least
one $\Hc$ between them. As an example, consider the $k$-landscape word $w=\Lc k
\Lc k \Hc kk \Hc \Lc k \Hc$. In $w$, neither the first two occurrences of
$\Lc$ nor the first two occurrences of $\Hc$ are separated; however, the
second and third occurrence of $\Hc$, as well as the second and third
occurrence of $\Lc$ are separated. We also say that two $k$-high clusters of a
partition are separated if there is at least one low cluster between them and
similarly, two low clusters are separated if there is a high cluster between
them.

Let $w$ and $w'$ be two $k$-landscape words. We say that $w$ and $w'$ are
\emph{$\Hc$-compatible} if they are compatible, and
moreover, they have the property that for any $i,j$, the $i$-th and $j$-th
occurrence of $\Hc$ in $w$ are separated if and only if the $i$-th
and $j$-th occurrence of $\Hc$ in $w'$ are separated. An $\Lc$-compatible pair of
words is defined analogously.

For example, the two compatible words $w=\Lc k \Hc kk \Hc \Lc$ and $w'=\Lc k \Hc
k \Lc \Hc k$ are $\Lc$-compatible (since the two occurrences of $\Lc$ are
separated in both words) but they are not $\Hc$-compatible (the two symbols
$\Hc$ are not separated in $w$ but they are separated in $w'$).

We are now ready to prove the following key lemma.

\begin{lemma}\label{lem-shuffle2}
Let $k$ be an integer. The following holds:
\begin{itemize}
\item[(1)]
Let $w$ and $w'$ be two $\Lc$-compatible $k$-landscape words, and
let $\tau$ be a 1-2-4 pattern. Let $\pi$ be an arbitrary partition,
and let $\pi'$ be the partition obtained from $\pi$ by the
$k$-shuffle from $w$ to~$w'$. For every $j<k$, $\pi$ contains $\tau$
at level~$j$ if and only if $\pi'$ contains $\tau$ at level~$j$.
Moreover, if the last symbol of $\tau$ is equal to 2, then the
previous equivalence also holds for every $j>k$.
\item[(2)]
Let $w$ and $w'$ be two $\Hc$-compatible $k$-landscape words, and
let $\tau$ be a 1-3-4 pattern. Let $\pi$ be an arbitrary partition,
and let $\pi'$ be the partition obtained from $\pi$ by the
$k$-shuffle from $w$ to~$w'$. For every $j>k$, $\pi$ contains $\tau$
at level~$j$ if and only if $\pi'$ contains $\tau$ at level~$j$.
Moreover, if the last symbol of $\tau$ is equal to 3, then the
previous equivalence also holds for every $j<k$.
\end{itemize}
\end{lemma}
\begin{proof}
We first prove (1). Assume that $\pi$ contains a 1-2-4 pattern
$\tau$ at level $j$. If $j>k$, it is easy to see that the occurrence
of $\tau$ is preserved by the shuffle as long as $\tau$ does not end
with a 1: we may use the same argument as in the proof of the first
part of Lemma~\ref{lem-shuffle}. Assume now that $j<k$. Let us write
$\tau=1232^p42^q12^r$ (the case when $\tau$ has the form
$1232^p12^q42^r$ is analogous). By assumption, $\pi$ contains a
subsequence $T$ order-isomorphic to $\tau$, with the symbol 2 of
$\tau$ corresponding to the symbol $j$ in $T$. Let us label the
$1+p+q+r$ occurrences of $j$ in $T$ by $j_0,j_1,\dotsb, j_{p+q+r}$,
in their natural left-to-right order. Let $a<b<c$ denote the symbols
of $T$ that correspond respectively to the symbols $1,3$ and $4$ in
$\tau$; we label the two occurrences of $a$ in $T$ by $a_0$ and
$a_1$. With this notation, we may write $T$ as follows:
\[
T=a_0j_0bj_1\dotsb j_p c j_{p+1}\dotsb j_{p+q} a_1 j_{p+q+1}\dotsb j_{p+q+r}.
\]
Now, we distinguish several cases, based on the relative order of $b,c$ and $k$:
If $c<k$, then all the symbols of $T$ are $k$-low and their relative position is
preserved by the shuffle, which means that $T$ is also a subsequence of $\pi'$.

If $c=k$, then the symbols $a<j<b$ are $k$-low. Let $x$ and $y$ be the two
symbols adjacent to $c$ in $T$ (typically $x=j_p$ and $y=j_{p+1}$, unless $q$
is zero, in which case $y=a_1$; recall that $c$ cannot directly follow $b$
and it cannot be the last element of $T$ by the definition of 1-2-4 pattern).
The elements $x$ and $y$ are low and they appear in two distinct low
clusters. After the shuffle, the occurrences of $a, b$ and $j$ in $T$ have
the same relative order, and the elements $x$ and $y$ still belong to
different clusters, which means that $\pi'$ contains a symbol greater than
$b$ between $x$ and $y$. This shows that $\pi'$ has a subsequence
order-isomorphic to $\tau$.

If $c>k$ and $b<k$, the argument from the previous paragraph applies as well.

It remains to consider the most complicated case: $c>k$ and $b\ge k$. This is
when we first use the $\Lc$-compatibility assumption. Since $b$ is not
$k$-low, we know that $j_1$ does not belong to the leftmost low cluster. Let
$x$ and $y$ be the two symbols adjacent to $c$ in $T$; by the definition
of 1-2-4 patterns, $x$ and $y$ are both $k$-low.

We know that $x$ and $y$ belong to distinct low clusters, and that their
clusters are separated, since $c$ is high. The shuffle preserves these
properties; in particular, in $\pi'$, the symbol $j_1$ does not belong to the
leftmost low cluster, which means that there is at least one non-low symbol
appearing in $\pi'$ before $j_1$. Since $\pi'$ is a partition in its
canonical sequential form, this implies that all the symbols $1,2,\dotsb, k$
appear in $\pi'$ in this order before $j_1$. Let $a', j'$ and $k'$
denote respectively the leftmost occurrences of $a, j$ and $k$ in $\pi'$. We
also know, from the $\Lc$-compatibility of $w$ and $w'$, that in $\pi'$ the
two symbols $x$ and $y$ appear in distinct and separated low clusters. In
particular, $\pi'$ contains a $k$-high symbol $c'$ between $x$ and $y$.
Putting it all together, we see that $\pi'$ contains the subsequence
\[
T'=a'j'k' j_1\dotsb j_p c' j_{p+1}\dotsb j_{p+q} a_1 j_{p+q+1}\dotsb j_{p+q+r},
\]
which is order isomorphic to $\tau$.

In all the cases, we see that if $\pi$ contains a 1-2-4 pattern
$\tau$ at level $j$, then $\pi'$ contains the same pattern at the
same level as well. The same proof also applies to the reverse
shuffle from $\pi'$ to $\pi$. This completes the proof of (1).

Claim (2) is proved by a similar argument. Let $\tau$ be a 1-3-4
pattern of the form $123^{p+1}13^q43^r$ (the case when
$\tau=123^{p+1}43^q13^r$ is analogous and easier). Assume that $\pi$
contains $\tau$ at level $j$, witnessed by a sequence $T$ of the
form
\[
T= a_0b j_0 j_1 \dotsb j_p a_1 j_{p+1}\dotsb j_{p+q} c j_{p+q+1}\dotsb
j_{p+q+r},
\]
with $a<b<j<c$.

If $j<k$, we apply the same argument as in the proof of the second claim of
Lemma~\ref{lem-shuffle}, to prove that if $\tau$ does not end with 4, the
occurrence of $\tau$ is preserved by the shuffle.

Next, we assume that $j>k$ and we distinguish several cases based on the
relative order of $a,b$ and $k$.

If $a>k$, then all the symbols of $T$ are $k$-high and their order is preserved
by the shuffle.

If $a=k$, or if $a<k$ and $b>k$, we let $x$ and $y$ denote the two symbols
adjacent to $a_1$ in $T$, and we observe that $\pi'$ has a non-high element $a'$ between $x$ and
$y$. The first occurrence of $a'$ in $\pi'$ must appear to the left of any $k$-high
symbol, hence $\pi'$ contains a subsequence $a'bj^{p+1}a'j^q c j^r$
order-isomorphic to~$\tau$.

If $a<k$ and $b\le k$, we define $x$ and $y$ as in the previous paragraph. This time, $x$ and $y$
belong to two separated high clusters, so $\pi'$ has a $k$-low element $a'$ between $x$
and $y$, and in particular,  $\pi'$ contains the subsequence $a'kj^{p+1}a'j^q c
j^r$.
\end{proof}

With the help of Lemma~\ref{lem-shuffle2}, we may prove all the equivalence
relations announced at the beginning of this section. We split the proofs into
four lemmas and then summarize the results in a theorem.

\begin{lemma}\label{lem-124}
Let $p,q\ge 1$. The pattern $\tau=1232^p412^q$ is equivalent to
$\tau'=1232^p42^q1$.
\end{lemma}
\begin{proof}
For an integer $k$ we say that a partition $\pi$ is a \emph{$k$-hybrid} if $\pi$
avoids $\tau$ at level $j$ for every $j\ge k$ and it avoids $\tau'$ at level
$j$ for every $j< k$. To prove the claim, it is enough to establish a bijection $f_k$
between $k$-hybrids and $(k+1)$-hybrids.

We say that a $k$-high cluster of $\pi$ is \emph{extra-high} if it contains a
symbol greater than $k+1$. We claim that $\pi$ contains $\tau$ at level $k$ if and
only if by scanning $\pi$ in the left-to-right direction we may find (not necessarily
consecutively) the leftmost high cluster, followed by $p$ occurrences of the symbol $k$,
followed by an extra-high cluster, followed by a low cluster, followed by $q$
occurrences of $k$. To see this, it suffices to notice that the leftmost high
cluster contains the symbol $k+1$, and to the left of this cluster we may always
find all the symbols $12\dotsb k$ in the increasing order.

By a similar argument, we see that $\pi$ contains $\tau'$ at level $k$ if and
only if it contains, left-to-right, the leftmost high-cluster, $p$ occurrences of
$k$, an extra-high cluster, $q$ occurrences of $k$ and a low cluster.

Now assume that $\pi$ is a $k$-hybrid partition. Let us try to find the leftmost
extra-high cluster $\Hc'$ of $\pi$ with the property that between $\Hc'$ and the
leftmost high cluster of $\pi$ there are at least $p$ occurrences of $k$. If no
such cluster exists, or if $\pi$ has less than $q$ symbols equal to $k$ to the
right of $\Hc'$, then $\pi$ avoids both $\tau$ and $\tau'$ at level $k$, and we define
$f_k(\pi)=\pi$.

Otherwise, let $w$ be the $k$-landscape of $\pi$. We will
decompose $w$ into the following concatenation:
\[
w=x\Hc'y k_q S_1 k_{q-1} S_2 \dotsb k_1 S_q,
\]
where $\Hc'$ represents the extra-high cluster defined above, $x$ is the prefix
of $w$ ending just before $\Hc'$, $k_i$ represents the $i$-th symbol $k$ in
$\pi$, counted from the right, $y$ is the subword of $w$ between $\Hc'$ and
$k_q$, and $S_i$ is the subword of $w$ between $k_{q-i+1}$ and $k_{q-i}$ with
$S_q$ being equal to the suffix of $w$ to the right of $k_1$. By construction,
none of the $S_i$'s contains the symbol $k$, so each of them is an alternating
sequence over the alphabet $\{\Lc,\Hc\}$, possibly empty. Since $\pi$ avoids
$\tau$ at level $k$, the subword $y$ does not contain the symbol $\Lc$.

As the next step, we decompose $S_1$ into two parts $S_1= H^* S_1^-$ as
follows: if the first letter of $S_1$ is $\Hc$, then we put $H^*=\Hc$ and
$S_1^-$ is equal to $S_1$ with the first letter removed; on the other hand,
if $S_1$ does not start with $\Hc$, then $H^*$ is the empty string and
$S_1^-=S_1$.

Now, let us define the word $w'$ as follows:
\[
w'=x\Hc' S_1^- k_1 S_2 k_2 S_3 k_3\dotsb k_{q-1} S_q k_q H^* y.
\]
It is not difficult to check that $w'$ is a landscape word (note that neither
$y$ nor $S_1^-$ can start with the symbol $\Hc$), and that $w'$ is
$\Lc$-compatible with $w$ (recall that $y$ contains no $\Lc$).

Let $\pi'$ be the partition obtained from $\pi$ by the shuffle from $w$ to $w'$.

Note that the prefix of $\pi$ up to the cluster $\Hc'$ (inclusive) is not
affected by the shuffle, because the words $w$ and $w'$ share the same prefix
up to the symbol $\Hc'$. In particular, the shuffle preserves the property
that $\Hc'$ is the leftmost extra-high cluster with at least $p$ symbols $k$
between $\Hc'$ and the leftmost high cluster of $\pi'$. It is routine to
check that $\pi'$ avoids $\tau'$ at level $k$. By Lemma~\ref{lem-shuffle2},
$\pi'$ is a $(k+1)$-hybrid partition. With these observations, it is easy to
see that for any given $(k+1)$-hybrid partition $\pi'$, we may uniquely
invert the procedure above and obtain a $k$-hybrid partition $\pi$.

Defining $f_k(\pi)=\pi'$, we obtain the required bijection between
$k$-hybrids and $(k+1)$-hybrids.
\end{proof}

The proofs of the following three lemmas follow the same basic argument as the
proof of Lemma~\ref{lem-124} above. The only difference is in the decompositions
of the corresponding landscape words $w$ and $w'$. We omit repeating the common parts of the
arguments and concentrate on pointing out the differences.

\begin{lemma}
Let $p,q\ge 1$. The pattern $\tau=1232^p142^q$ is equivalent to
$\tau'=12312^p42^q$.
\end{lemma}
\begin{proof}
A partition $\pi$ contains $\tau$ at level $k$ if and only if it
contains, in left-to-right order, the leftmost high cluster, $p$
copies of $k$, a low cluster, an extra-high cluster, and $q$ copies
of $k$. Similar characterization applies to $\tau'$.

Let $\Hc_1$ denote the leftmost high cluster of $\pi$, let $\Hc'$ denote the
rightmost extra-high cluster of $\pi$ that has the property that there are at
least $q$ occurrences of $k$ to the right of $\Hc'$. If $\Hc'$ does not exist,
or if there are less than $p$ symbols $k$ between $\Hc_1$ and $\Hc'$, then
$\pi$ contains neither $\tau$ nor $\tau'$ at level $k$ and we put
$f_k(\pi)=\pi$. Otherwise, let $w$ be
the landscape of $\pi$, and let us write
\[
w=x\Hc_1 S_1 k_1 S_2 k_2 \dotsc S_p k_p y \Hc' z
\]
where none of the $S_i$ contains $k$, and $y$ avoids $\Lc$. Define
$S_p^-$ and $H^*$ by writing $S_p=S_p^- H^*$ where $S_p^-$ does not
end with the letter $\Hc$ and $H^*$ is equal either to $\Hc$ or to the empty
string, depending on whether $S_p$ ends with $\Hc$ or not.

Now we write
\[
w=x\Hc_1 \overline y k_1 H^* S_1 k_2 S_2\dotsb k_p S_p^- \Hc' z,
\]
where $\overline y$ is the reversal of $y$. The rest of the proof is analogous to
Lemma~\ref{lem-124}.
\end{proof}

We now apply the same arguments to 1-3-4 patterns.

\begin{lemma}\label{lem-134}
For any $p\ge 0$ and $q\ge 1$, the pattern $\tau=123^{p+1}13^q4$ is equivalent to
the pattern $\tau'=123^{p+1}143^q$.
\end{lemma}
\begin{proof}
As usual, a $k$-hybrid is a partition that avoids $\tau$ at every
level $j\ge k$ and that avoids $\tau'$ at every level below $k$.

Let us say that a $k$-cluster of a partition $\pi$ is \emph{extra-low} if it
contains a symbol smaller than $k-1$. A partition contains $\tau$ at level $k$ if and only if it has $p+1$
occurrences of $k$ followed by an extra-low cluster, followed $q$ symbols $k$,
followed by a high cluster; similarly, a partition contains $\tau'$ at level
$k$ if and only if it has $p+1$ copies of $k$, followed by an extra-low
cluster, followed by a high cluster, followed by $q$ copies of $k$.

Assume $\pi$ is a $k$-hybrid partition. Let $\Lc'$ denote the leftmost
extra-low cluster of $\pi$ that has at least $p+1$ copies of $k$ to its left.
If $\Lc'$ does not exist, or if it has less than $q$ copies of $k$ to its
right, we put $f_k(\pi)=\pi$; otherwise, we decompose the landscape word $w$
of $\pi$ as follows:
\[
w=x\Lc' S_1 k_1 S_2 k_2\dotsb S_{q-1} k_{q-1} S_q k_q y,
\]
where the $S_i$ do not contain $k$, and by assumption, $y$ avoids $\Hc$.
Next, we write $y=L^* y^-$ where $L^*$ is an empty string or a
single symbol $\Lc$, and $y^-$ does not start with $\Lc$. We
define $w'$:
\[
w'=x\Lc' y^- k_1 L^* S_1 k_2 \dotsb S_{q-1} k_q S_q.
\]
The words $w$ and $w'$ are $\Hc$-compatible. We define the bijection between
$k$-hybrids and $(k+1)$-hybrids in the usual way.
\end{proof}

\begin{lemma}
For every $p\ge 0$ and $q\ge 1$, the pattern $\tau = 123^{p+1}413^q$ is
equivalent to the pattern $\tau'=12343^{p}13^q$.
\end{lemma}
\begin{proof}
As before, take $\pi$ to be a $k$-hybrid partition. Let $\Lc'$ be the
rightmost extra-low cluster that has at least $q$ copies of $k$ to its right.
If $\Lc'$ has at least $p+1$ copies of $k$ to its left, we perform the
following decomposition of the landscape $w$ of $\pi$:
\[
w=\Lc k_1 S_1 k_2 S_2\dotsb k_p S_p k_{p+1} y \Lc' z.
\]
Next, we write $S_p=S_p^- L^*$ with the usual meaning and define
\[
w'=\Lc k_1 L^* \overline y k_2 S_1 k_3 S_2\dotsb S_{p-1} k_{p+1} S_p^-
\Lc' z.
\]
The rest is the same as before.
\end{proof}

We summarize our results:

\begin{theorem}
For every $p,q\ge 0$, we have the following equivalences:
\begin{enumerate}
\item $1232^{p}142^{q}\sim 12312^p 42^q$
\item $1232^p 412^q \sim 1232^p42^q 1$
\item $123^{p+1}413^q\sim 12343^p13^q$
\item $123^{p+1}143^q \sim 123^{p+1}13^q 4$
\end{enumerate}
\end{theorem}
\begin{proof}
If $p$ and $q$ are both positive, the results follow directly from the four
preceding lemmas.

If $p=0$, the first and the third claim are trivial, the second one is a
special case of Corollary~\ref{cor-spq}, and the fourth is covered by
Lemma~\ref{lem-134}.

If $q=0$, the first claim is a special case of Corollary~\ref{cor-spq}, the
second and the fourth are trivial, and the third follows from
Theorem~\ref{thm-spq2}.
\end{proof}

\section{The avoidance of four-letter patterns}\label{sec-1123}
In this section, we will complete the classification of the
equivalence classes of the patterns of length four, by proving the
equivalence $1212\sim 1123$. Unlike in the previous arguments, we do
not present a direct bijection between pattern-avoiding classes, but
rather we prove that $p(n;1123)$ is equal to the $n$-th Catalan
number. Since it is well known that noncrossing partitions are
enumerated by the $n$-th Catalan number as well, this will yield the
desired equivalence. All the other equivalent pairs of patterns of
length four are covered by the general theorems proved in the
previous parts of the paper (see Table~\ref{tab2}).
\begin{table}[h]
\begin{center}
\begin{tabular}{ l|l }
  $\tau$ & $p(n;\tau)$ \\ \hline\hline
  $1111$ & \cite[Sequence A001680]{oeis} (see Equation~\eqref{eq-1m})\\
  $1112$, $1121$, $1211$, $1222$ & \cite[Sequence A005425]{oeis} (see Corollary~\ref{co1222}) \\
  $1122$         & $1,1,2,5,14,42,133,441,\dotsc$ \\
  $1123$, $1212$, $1221$         & $\frac{1}{n+1}\binom{2n}{n}$ \\
  $1213$, $1223$, $1231$, $1232$, $1233$, $1234$ &\cite[Sequence A007051]{oeis}
  (see Equation~\eqref{eq-123k})\\[5pt]
\end{tabular}
\caption{Number of partitions in $P(n;\tau)$, where $\tau\in P(4)$.}\label{tab2}
\end{center}
\end{table}
\subsection{Enumeration of $1123$-avoiding partitions}
As we said before, our aim is to show that the 1123-avoiding partitions of
$[n]$ are enumerated by the $n$-th Catalan number, i.e.,
$p(n;1123)=\frac{1}{n+1}\binom{2n}{n}$.

We achieve this by proving that 1123-avoiding partitions of $n$ are in bijection
with Dyck paths of semilength $n$. A \emph{Dyck path of semilength $n$} is a nonnegative
path on the two-dimensional integer lattice from $(0,0)$ to $(2n,0)$ composed of
\emph{up-steps} connecting $(x,y)$ to $(x+1,y+1)$ and \emph{down-steps} connecting $(x,y)$ to
$(x+1,y-1)$. It is well known that these paths are enumerated by Catalan
numbers.

We will need the following refinement: let $T(n,k)$ denote the set of
1123-avoiding partitions $\pi$ of the set $[n]$ with the property that
$\pi_n=k$. Let $t(n,k)$ be the cardinality of $T(n,k)$. We will prove that
$t(n,k)$ is equal to the number of Dyck paths of semilength $n$ whose last
up-step is followed by exactly $k$ down-steps.

Let $T'(n,k)$ denote the set of Dyck paths of semilength $n$ such that their
last up-step is followed by exactly $k$ down-steps, let $t'(n,k)$ be the
cardinality of  $T'(n,k)$. We remark that standard bijections between Dyck
paths and pattern-avoiding permutations show that $t'(n,k)$ is also equal to
the number of 123-avoiding permutations $(\pi_1,\pi_2,\dotsc,\pi_n)$ such that
$\pi_n=k$. Our aim is to prove the following result:
\begin{theorem}\label{thm-main}
For every $n,k$, $t(n,k)$ is equal to $t'(n,k).$
\end{theorem}

Before starting the proof of the theorem, we introduce more
definitions.

\begin{definition}
A \emph{123-avoiding sequence} is a sequence $s_1,s_2,\dotsc,s_\ell$ of
positive integers, such that there are no three indices $i<j<k$ that would
satisfy $s_i<s_j<s_k$. We define the \emph{rank} of a sequence to be equal to
$\ell+m-1$, where $\ell$ is the length of the sequence and
$m=\max\{s_i,i=1,\dotsc,\ell\}$ is the largest element of the sequence.
\end{definition}

For example, there are five 123-avoiding sequences of rank 3: those are the
sequences (1,1,1), (1,2), (2,1), (2,2), and (3). There are fourteen
123-avoiding sequences of rank 4: (1,1,1,1), (1,1,2), (1,2,1), (1,2,2),
(2,1,1), (2,1,2), (2,2,1), (2,2,2), (1,3), (2,3), (3,1), (3,2), (3,3), and~(4).

The proof of Theorem~\ref{thm-main} is divided into the following two claims:

\begin{claim}\label{cla-first}
A 1123-avoiding partition $\pi$ of $[n]$ with $m$ blocks has the following form:
\begin{equation}\label{eq-tail}
\pi = 123\dotsb (m-2)(m-1) S
\end{equation}
where $S$ is a 123-avoiding sequence of rank $n$, with maximum element $m$.
Conversely, If $S$ is any 123-avoiding sequence of rank $n$ with maximum element
$m$ then $\pi$ defined by the formula~\eqref{eq-tail} is a canonical sequence of a
1123-avoiding partition of $[n]$.

In particular, the number of 123-avoiding sequences of rank $n$ with last
element $k$ is equal to $t(n,k)$.  \end{claim}

\begin{claim}\label{cla-second}
The numbers $t(n,k)$ satisfy the following recurrences:
\begin{align}
t(1,1)&=1   \label{r-1} \\
t(n,k)&=0 \quad\text{if}\quad k<1 \quad\text{or}\quad k>n \label{r-2}\\
t(n,k)&=\sum_{j=k-1}^{n-1} t(n-1,j) \quad\text{for}\quad n\ge 2, n\ge k\ge 1  \label{r-3}
\end{align}

\end{claim}

It is easy to see the recurrences
\eqref{r-1},\eqref{r-2} and \eqref{r-3} would all hold if $t$ were replaced by
$t'$: given a Dyck path from
$T'(n,k)$, we erase its last up-step and the following down-step, to obtain a
Dyck path from $\cup_{j=k-1}^{n-1} T'(n-1,j)$; in particular,
Claim~\ref{cla-second} implies that
$t(n,k)=t'(n,k)$. Thus, the two claims also show that the number of 123-avoiding sequences of rank
$n$ is exactly the $n$-th Catalan number.

\begin{proof}[Proof of Claim~\ref{cla-first}]
Let $\pi$ be a $1123$-avoiding partition of $[n]$ with $m$ blocks. Observe that
for every $i\in[m-1]$, the symbol $\pi_i$ is equal to $i$, otherwise $\pi$ would
contain the forbidden pattern. It follows that $\pi$ can be decomposed according
to the equality $\pi = 123\dotsb (m-2)(m-1) S$, where
the sequence $S$ has length $l=n-m+1$ and maximum element equal
to $m$, hence $S$ has rank~$n$. Also, the last element of $S$ is equal to
$k$.

It remains to check that $S$ is 123-avoiding: if $S$ contained a subsequence $xyz$
for $x<y<z$ then the original partition would contain a subsequence $xxyz$,
which is forbidden. It follows that $S$ obtained from a 1123-avoiding partition
$\pi$ has all the required properties.

The ``converse'' part of the claim is equally easy to verify, and we omit it.
\end{proof}

Claim~\ref{cla-first} motivates the following definition;
\begin{definition} Let $\pi$ be a 1123-avoiding partition of $[n]$ with $m$
blocks. The 123-avoiding sequence $S$ obtained from the decomposition according
to the formula~\eqref{eq-tail} will be called \emph{the tail of $\pi$}. Let
$T_0(n,k)$ be the set of all the tails of the partitions from $T(n,k)$;
equivalently, $T_0(n,k)$ is the set of 123-avoiding sequences of rank $n$ with
the last element equal to $k$.
\end{definition}

We are now ready to prove Claim~\ref{cla-second}.

\begin{proof}[Proof of Claim~\ref{cla-second}]
Only the recurrence \eqref{r-3} is nontrivial. To prove the recurrence, we need
a bijection from $T(n,k)$ to $\cup_{j=k-1}^{n-1} T(n-1,j)$. It is more
convenient to work with the tails of the partitions, i.e., to describe a
bijection between $T_0(n,k)$ and $\cup_{j=k-1}^{n-1} T_0(n-1,j)$. We will
construct the required bijection as the union of two injective maps
$f_1$ and
$f_2$ with the property that the domains of $f_1$ and $f_2$
form a disjoint partition of $T_0(n,k)$ and their ranges form a disjoint
partition of $\cup_{j=k-1}^{n-1} T_0(n-1,j)$.

Let $S\in T_0(n,k)$ be a 123-avoiding sequence of length $\ell$. The sequence $S$
can be uniquely decomposed into a concatenation of the form $S=S_01^b k$, where
$S_0$ is a (possibly empty) prefix of $S$ whose last element is different from~1.

We distinguish the following two cases:

\textbf{Case 1.} If $S_0$ is nonempty and the last element of $S_0$ is greater
than or equal to $k$, we define $S'=f_1(S)=S_0(k-1)^b$. It is easy to see that
$f_1$ is injective. By construction, the length
of $S'$ is equal to $\ell -1$ and the maximum of $S'$ has the same value
as the maximum of $S$, so $S'$ has rank $n-1$. It is also easy to check that
$S'$ is 123-avoiding, and that $S'\in \cup_{j=k-1}^{n-1} T_0(n-1,j)$.

\textbf{Case 2.} We now deal with the case when $S_0$ is empty, or the last element of $S_0$ is smaller than
$k$. We first observe that all the elements of $S_0$ are greater
than 1. Indeed, the last element of $S_0$ is never equal to 1 by definition, and
if $S_0$ contained an element 1 before the last one, then $S$ would contain a
subsequence $1jk$, where $j$ is the last element of $S_0$ and $k$ the last
element of $S$; however, this is impossible, because $S$ avoids 123. Now, we
define $S'=f_2(S)=(S_0-1) (k-1)^{b+1}$, where $S_0-1$ denotes the sequence $S_0$
with all its elements decreased by 1. Note that the last symbol of $S_0-1$ is
smaller than $k-1$, which implies that $f_2$ is an injective map.
Clearly, the length of $S'$ is $\ell$ and the maximum of $S'$ is one less than
the maximum of $S$, so $S'$ has rank $n-1$.  Also, $S'$ is easily seen to be
123-avoiding, and $S'\in T_0(n-1,k-1)$.

To finish the proof, we need to check that every sequence
$S'\in\cup_{j=k-1}^{n-1} T_0(n-1,j)$ is in the range of exactly one of the
two injections $f_1$ and $f_2$. To see this, express $S'$ as the concatenation
$S'=S'_0 (k-1)^c$, where $c\ge 0$ and $S'_0$ is a (possibly empty) prefix of
$S'$ whose last element is different from $k-1$. The sequence $S'$ is in the
range of $f_1$ if and only if $S_0'$ is nonempty and the last element of $S_0'$
is at least $k$; in such case, we have $S'=f_1(S_0' 1^c k)$. On the other hand,
if $S'_0$ is empty or if its last element is smaller than $k-1$, then necessarily
$c\ge 1$ and $S'=f_2( (S'_0+1)1^{c-1} k)$. This completes the proof.
\end{proof}

We may use Theorem~\ref{thm-main} to derive the closed-form expression for
$t(n,k)$. Since the number of Dyck paths that end with an up-step followed by
$k$ down-steps is equal to the number of non-negative lattice paths from
$(0,0)$ to $(2n-k-1,k-1)$, we may apply standard arguments for the enumeration
of non-negative lattice paths to obtain the formula
\[
t(n,k)=\frac{k}{n}\binom{2n-k-1}{n-1}.
\]

Let us remark that it is possible to obtain a closed form formula for $t(n,k)$
without Theorem~\ref{thm-main}, using directly the recurrences from
Claim~\ref{cla-second}. This is achieved by the kernel method techniques as
described, e.g., in~\cite{Man1}. We omit the details here.

\section{The patterns of size five}\label{sec-12112}

For a full characterization of the equivalence of patterns up to size seven,
we need to consider one more isolated case, namely the pattern $12112$. Our
aim is to show that this pattern is equivalent to the three patterns $12221$,
$12212$, and $12122$. Note that the latter three patterns are all equivalent
by Corollary~\ref{cor-spq}; it is thus sufficient to show that $12112\sim
12212$. The proof we are about to present is rather long and occasionally
technical, however, it is well worth the effort, since the equivalence of
$12112$ and $12212$ has several consequences related to fillings of Ferrers
shapes and pattern-avoiding graphs. We will discuss these consequences in
greater detail at the end of this section.

We remark that contrary to the case of the three equivalent patterns
$12221\sim 12212 \sim 12122$, whose equivalence was obtained as a consequence
of the Ferrers-equivalence of the corresponding matrices $M(2,221)\eqf
M(2,212)\eqf M(2,122)$, the proof involving the pattern $12112$ does not use
the notion of Ferrers equivalence. In fact, the matrix $M(2,112)$ is not
Ferrers-equivalent to the three matrices above.

\subsection{Introduction}
We will first introduce the basic terminology and notation that we will use throughout the
proof.

Let $S=s_1s_2\dotsc s_n$ be a sequence of length $n$ over the alphabet $[m]$,
such that every element of $[m]$ appears in $S$ at least once. For $i\in[m]$ let
$f_i$ and $l_i$ denote the index of the first and the last symbol of $S$ that
is equal to $i$; formally, $f_i=\min\{j;\ s_j=i\}$ and $l_i=\max\{j;\
s_j=i\}$.

\begin{definition}\label{def-semican}
For $k\in[m]$, we say that the sequence $S$ is a \emph{$k$-semicanonical
sequence} (or $k$-sequence for short), if $S$ has the following properties:
\begin{itemize}
\item For every $i,i'$ such that $1\le i<k$ and $i<i'$, we have $f_i<f_{i'}$.
\item For every $i,i'$ such that $k\le i<i'\le m$, we have $l_i<l_{i'}$.
\end{itemize}
\end{definition}

Note that $m$-semicanonical sequences are precisely the canonical sequences of
partitions of $[n]$ with $m$ blocks (i.e., the sequences satisfying
$f_i<f_{i+1}$ for $i\in[m-1]$), while the $1$-canonical sequences are precisely the sequences
satisfying $l_i<l_{i+1}$ for $i\in[m-1]$. Also, for every fixed $k\in[m]$ and a fixed partition
$\pi$ with $m$ blocks, there is exactly one $k$-sequence $S$ with
the property $s_i=s_j \iff \pi_i=\pi_j$; this sequence $S$ can be obtained from
$\pi$ by preserving the numbering of the first $k-1$ blocks of $\pi$, and by
numbering the remaining $m-k+1$ blocks in the increasing order of their largest
element.

In particular, assuming
$n$ and $m$ are fixed, the number of $k$-sequences is independent of $k$, and
each partition of $[n]$ with $m$ blocks is
represented by a unique $k$-sequence.  To prove the equivalence $12112\sim
12212$, we will exploit a remarkable property of the pattern $12112$,
described by the following key lemma.

\begin{lemma}\label{lem-semican}
For every fixed $n$ and $m$, the number of $12112$-avoiding $k$-sequences is
independent of~$k$. Thus, for every $k\in[m]$, the number of $12112$-avoiding
$k$-sequences of length $n$ with $m$ symbols is equal to the
number of $12112$-avoiding partitions of $n$ with $m$ blocks.  \end{lemma}

We stress that a $12112$-avoiding $k$-sequence can actually
represent a partition that contains $12112$ in its canonical representation.

Before we present the proof of Lemma~\ref{lem-semican}, let us explain how the
lemma implies the equivalence $12112\sim 12212$.

\begin{theorem}\label{thm-12112}
The pattern $12112$ is equivalent to $12212$. In fact, for every $m$ and $n$, there is a bijection
between $12112$-avoiding partitions of $[n]$ with $m$ blocks and
$12212$-avoiding partitions of $[n]$ with $m$ blocks.
\end{theorem}
\begin{proof} Fix $m$ and $n$. We know that the $12112$-avoiding partitions of $[n]$
with $m$ blocks correspond precisely to $m$-semicanonical sequences over $[m]$
of length $n$, and by Lemma~\ref{lem-semican}, these sequences are in
bijection with $1$-semicanonical $12112$-avoiding sequences of the same
length and alphabet. It remains to provide a bijection between the
$12112$-avoiding 1-sequences and the $12212$-avoiding partitions, which is
done as follows: take a $1$-semicanonical $12112$-avoiding sequence $S$ with
$m$ symbols and length $n$, reverse the order of letters in $S$, and then
replace each symbol $i$ of the reverted sequence by the symbol $m-i+1$. It is
easy to check that this transform is an involution which maps
$12112$-avoiding 1-sequences onto $12212$-avoiding $m$-sequences, which are
precisely the $12212$-avoiding partitions of $[n]$ with $m$ blocks.
\end{proof}

It now remains to prove Lemma~\ref{lem-semican}. For the rest of the proof,
let as assume that $m$ and $n$ are fixed, and that each sequence we
consider has length $n$ and $m$ distinct symbols, unless otherwise noted.

In the following arguments, it is often convenient to represent a sequence
$S=s_1\dotsb s_n$ by the matrix $M(S,m)$ (recall that $M(S,m)$ is the 0-1
matrix with a 1-cell in row $i$ and column $j$ if and only if $s_j=i$). A
matrix representing a $k$-sequence will be called \emph{$k$-semicanonical
matrix} (or just $k$-matrix), and a matrix representing a $12112$-avoiding
sequence will be simply called \emph{$12112$-avoiding
matrix}. In accordance with earlier terminology, we will use the term
\emph{sparse matrix} for a 0-1 matrix with at most one 1-cell in each
column, and we will use the term \emph{semi-standard matrix} for a
0-1 matrix with exactly one 1-cell in each column. For a 0-1 matrix $M$,
we let $f_i(M)$ and $l_i(M)$ denote the column-index of the first and the last
1-cell in the $i$-th row of $M$. We will write $f_i$ and $l_i$ instead of
$f_i(M)$ and $l_i(M)$ if there is no risk of confusion.

Before we formulate the proof of our key lemma, let us present a brief sketch
of the main idea. Assume we want to build a bijection that transforms a
$(k+1)$-matrix $M$ into a $k$-matrix (ignoring $12112$-avoidance for a
while). Such bijection is easy to obtain: assume that the last 1-cell in row
$k$ is in column $c$, let us call the row $k$ \emph{the key row of $M$}. If
the last 1-cell in row $k+1$ appears to the right of column $c$, then $M$ is
already a $k$-matrix and we do not need to do anything; on the other hand, if
row $k+1$ has no 1-cell to the right of $c$, we swap the key row $k$ with
the row $k+1$, to obtain a new matrix $M'$ whose key row is now the row
$k+1$. We now repeat the same procedure: we compare the position of the last
1-cell in the key row $k+1$ (which is in the column $c$, as we know) with the
last 1-cell in the row $k+2$, and if necessary, we swap the key row with the
row directly above it, until we reach the situation when the key row is
either the topmost row of the matrix, or the row above the key row has a
1-cell to the right of column $c$. This procedure transforms the original
$k+1$ matrix into a $k$-matrix. Also, the procedure is invertible (note that
the first 1-cell of the key row is always to the left of any other 1-cell in
the rows $k, k+1, \dotsc,m$).

Unfortunately, this simplistic approach does not preserve $12112$-avoidance.
However, we will present a more sophisticated algorithm, which follows the same
basic structure as the procedure above, but instead of simply swapping the key
row with the row above it, it performs a more complicated step. The description
of this step is the fundamental ingredient of our proof.

To formalize our argument, we need to introduce more definitions. Let $M$ be a
0-1 matrix with exactly one 1-cell in each column and at least one 1-cell in
each row, and let us write $f_i=f_i(M)$ and $l_i=l_i(M)$. Let $k,q,p$ be
three row-indices of $M$, with $k\le p\le q$. We will say that $M$ is a
\emph{$(k,p,q)$-matrix}, if $M$ has the following form:
\begin{itemize}
\item The matrix obtained from $M$ by erasing row $p$ is a $k$-semicanonical matrix with $m-1$ rows.
\item For each $i<k$, we have $f_i<f_p$. For every $j\ge k$, $j\neq p$, we have $f_p<f_j$.
\item The number $q$ is determined by the relation $q=\max\{j;\ l_j\le l_p\}$. By the first condition, this implies that $l_j\le l_p$ for every $j\in\{k,k+1,\dotsc,q\}$.
\end{itemize}
In a $(k,p,q)$-matrix, row $p$ will be called \emph{the key row}.

Intuitively, a $(k,p,q)$-matrix is an intermediate stage of the above-described
procedure which transforms a $(k+1)$-matrix into a $k$-matrix by moving the key
row towards the top; the number $p$ is the index of the key row in a given step
of the procedure, while the number $q$ is the topmost row that needs to be
swapped with the key row to produce the required $k$-matrix.

In particular, a matrix $M$ is $(k+1)$-semicanonical if and only if it is a $(k,
k,q)$-matrix, and $M$ is $k$-semicanonical if and only if it is a $(k,q,q)
$-matrix.

As an example, consider the sequence $S=1331232431$ with $n=10$ and
$m=4$. This sequence corresponds to the following matrix $M=M(S,4)$.
\begin{align*}
M=\left(\begin{array}{l}
0000000100\\
0110010010\\
0000101000\\
1001000001\\
\end{array}\right) &&
M'=\left(\begin{array}{l}
0110010010\\
0000000100\\
0000101000\\
1001000001\\
\end{array}\right)
\end{align*}
The matrix $M$ is a $(2,3,4)$-matrix. If we exchange the third row (which acts
as the key row) with the fourth row, we obtain a $(2,4,4)$-matrix $M'$
representing the $2$-sequence $S'=1441242341$. The matrix $M'$ can also be
regarded as a $(1,1,4)$-matrix, with the key row at the bottom.

In general, a matrix $M$ is $(k+1)$-semicanonical if and only if it is a $(k,
k,q)$-matrix, and $M$ is $k$-semicanonical if and only if it is a $(k,q,q)
$-matrix. To prove Lemma~\ref{lem-semican}, we will prove the following
lemma.

\begin{lemma}\label{lem-kpq}
For arbitrary $k\le p<q$, there is a bijection $\phi$ between $12112$-avoiding $(k,p,q)
$-matrices and $12112$-avoiding $(k,p+1,q)$-matrices.
\end{lemma}

Clearly, Lemma~\ref{lem-kpq} implies Lemma~\ref{lem-semican}. Before we
construct the bijection $\phi$ and prove its correctness, we need to prove several
basic properties of the $12112$-avoiding $(k,p,q)$-matrices.

\subsection{Tools of the proof}

We will use the following terminology: if $x\in[m]$ is a row of a matrix $M$, then an
\emph{$x$-column} is a column of $M$ that has a 1-cell in row $x$. Similarly, if
$X\subseteq[m]$ is a set of rows of $M$, we will say
that a column $j$ is an \emph{$X$-column} if $j$ has a 1-cell in a row belonging
to $X$.

If $x,y$ is a pair of rows of $M$ with $x<y$, we will say that \emph{$M$ contains
$12112$ in $(x,y)$} if the submatrix of $M$ induced by the pair of rows $x,y$
contains $12112$. If $X$ and $Y$ are two sets of rows, we will say that $M$
\emph{contains $12112$ in $(X,Y)$} if there is an $x\in X$ and $y\in Y$ such
that $x<y$ and $M$ contains $12112$ in $(x,y)$.

Throughout this section, we will assume that $k,p,q$ are fixed, and that
$k\le p<q$.

We now state a pair of simple but useful observations. Their proofs are
straightforward, and we omit them.

\begin{observation}\label{obs-f}
Let $M$ be a sparse 0-1 matrix, and let $x<y$ be two rows of $M$, such that
$f_x<f_y$. The matrix $M$ avoids $12112$ in $(x,y)$ if and only if $M$ has
at most one $x$-column $s$ satisfying $f_y<s<l_y$. If such a unique column $s$ exists, we
will say that $s$ \emph{separates} row $y$. The $y$-columns
that are to the left of the separating column $s$ will be called
\emph{front $y$-columns (with respect to row~$x$)} and their
1-cells will be called \emph{front 1-cells}, and similarly, the $y$-columns to
the right of $s$ will be called \emph{rear $y$-columns} and their 1-cells are
\emph{rear 1-cells}. If there is no such separating column, then we will assume
that all the $y$-columns and their 1-cells are front.  \end{observation}

\begin{observation}\label{obs-l}
Let $M$ be a sparse 0-1 matrix, and let $x<y$ be a pair of rows such that $l_x<l_y$. Let $t$
be the number of 1-cells in row $x$, and let $c_i$ be the
$i$-th $x$-column, i.e., $f_x=c_1<c_2<\dotsb<c_t=l_x$. The
matrix $M$ avoids $12112$ in $(x,y)$, if and only if every $y$-column
appears either to the left of column $c_1$, or between the columns $c_{t-1}$
and $c_{t}$, or to the right of column $c_t$. These three types of $y$-columns
(and their 1-cells) will be called \emph{left, middle, and right $y$-columns (or 1-cells) with respect to row $x$}.
\end{observation}

Next, we prove a lemma that will greatly simplify our task of constructing the
bijection $\phi$ and proving its correctness.

\begin{lemma}\label{lem-low}
Let $M$ be a $12112$-avoiding $(k,p,q)$-matrix, Let $j$ be a row of $M$ with
$k\le j\le p$. Let $M'$ be a sparse
0-1 matrix of the same size as $M$, with the property that for every
$i\not\in\{j,j+1,\dotsc,q\}$, the $i$-th row of $M$ is equal to the $i$-th row
of $M'$. If $M'$ has a copy of the pattern $12112$ in a pair of rows $x<y$, then $j\le
x\le q$.
\end{lemma}
\begin{proof}
Let $M$ and $M'$ be as above. We will call the
rows $\{j,j+1,\dotsc,q\}$ \emph{mutable}, and the remaining rows will be
called \emph{constant}.

Assume that $M'$ has a copy of the forbidden pattern in the rows $x<y$.
Clearly, at least one of the two rows $x,y$ must be mutable, and in particular,
we must have $x\le q$. The lemma claims that $x$ must be mutable. For
contradiction, assume that $x<j$. We now distinguish two cases.

\textbf{The case $x<k$:} Necessarily, $y$ is one of the mutable rows. From
the definition of the $(k,p,q)$-matrix, we obtain that all the columns of $M$
to the left of $f_p(M)$ and to the right of $l_p(M)$ contain a 1-cell in one
of the constant rows. Since $M'$ is sparse, we conclude that in $M'$, all the
1-cells in the mutable rows can only appear in the columns $i$ such that $f_p(M)\le
i\le l_p(M)$.

Now, we apply Observation~\ref{obs-f} to the
rows $x$ and $p$ in the matrix $M$, to conclude that $M$ (and hence also $M'$) has at most one
$x$-column $s$ such that $f_p(M)\le s\le l_p(M)$, and therefore $M'$ also has at most
one $x$-column between $f_y(M')$ and $l_y(M')$. By Observation~\ref{obs-f}, this shows that $x$ cannot
form the forbidden pattern with any of the mutable rows $y$ of $M'$.

\textbf{The case $k\le x<j$}. As before, we have $y\in\{j,\dotsc,q\}$. Let
$c_1<c_2<\dotsb<c_t$ be the $x$-columns of $M$ (and hence of $M'$ as well, since
$x$ is constant). For any
mutable row $i$, we have $l_x(M)<l_i(M)$ by the definition of $(k,p,q)$-matrix.
Observation~\ref{obs-l}, applied to
the pair of rows $x,i$ in $M$, tells us that all the $i$-columns of $M$
appear either to the left of $c_1$ or to the right of $c_{t-1}$. In
particular, all the 1-cells between the columns $c_1$ and $c_{t-1}$ belong to
the constant rows. This implies that $M'$ can have no occurrence of
$12112$ in the two rows $x<y$.
\end{proof}

We will now describe a simple operation on $12112$-avoiding pairs of rows. This
operation, which we will call \emph{pseudoswap} will play an important part in
the construction of the bijection~$\phi$.

Assume that $M$ is a sparse matrix
with a pair of adjacent rows $x<y$ (where $y=x+1$) that avoids $12112$ in
$(x,y)$. Assume furthermore that $f_x<f_y\le l_y<l_x$.  The pseudoswap of the
two rows is performed as follows.

\begin{figure}
\includegraphics[scale=0.8]{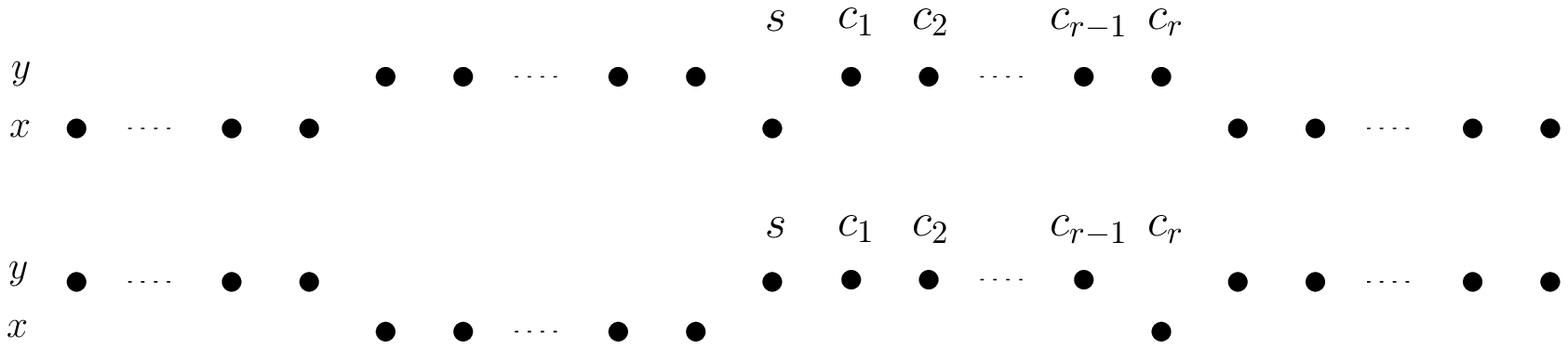}
\caption{The illustration of the hard case of the pseudoswap operation.}
\label{fig-pswap}
\end{figure}

\begin{description}
\item[Easy case]
If the row $y$ is not separated by an $x$-column (in the sense of
Observation~\ref{obs-f}), or if $M$ has at most one rear $y$-column with
respect to row $x$, the pseudoswap is performed by changing the order of the
two rows; in other words, each 1-cell in row $x$ moves into row $y$ in the
same column, and vice versa.

\item[Hard case]
Assume $M$ has an $x$-column $s$ separating $y$, and that it has $r>1$ rear
$y$-columns $c_1<c_2<\dotsb<c_r$ (see Figure~\ref{fig-pswap}). In this case,
the pseudoswap preserves the position of all the 1-cells in columns
$c_1,\dotsc,c_{r-1}$ (i.e., the 1-cells in these columns remain in row $y$),
and all the other 1-cells of in rows $x,y$ are moved from $x$ to $y$ and vice
versa. Note that after the pseudoswap is performed, the columns
$s<c_1<c_2<\dotsb<c_{r-1}$ all contain a 1-cell in row $y$, and these $r$
1-cells are precisely the middle 1-cells of $y$ with respect to $x$.  (in the
sense of Observation~\ref{obs-l}).  \end{description}

Let $M'$ be the matrix obtained from $M$ by the pseudoswap. It can be routinely
checked that $M'$ avoids $12112$ in $(x,y)$. Let us write $f'_i$ for $f_i(M')$
and $l'_i$ for $l_i(M')$. Clearly, $f'_x=f_y$ and $f'_y=f_x$, and also
$l'_x=l_y$ and $l'_y=l_x$. Also, if $M$ has $r\ge 0$ rear cells in row $y$, then
$M'$ has $r$ middle cells in row $y$.

It is not difficult to see that the pseudoswap can be inverted in the
following way: let $M'$ be a sparse matrix avoiding $12112$ in two adjacent rows
$x<y$, such that $f'_y<f'_x\le l'_x<l'_y$. We again
distinguish two cases: if $M'$ has less than two middle $y$-columns, we invert
the easy case of the pseudoswap by exchanging the two rows; on the other
hand, if $M'$ has $r>1$ middle $y$-columns $m_1<\dotsb< m_r$, we
invert the hard case by preserving the position of the 1-cells in columns
$m_2,m_3,\dotsc,m_r$ and inverting all the other $\{x,y\}$-columns.

We will be mostly interested in the situation when the pseudoswap is applied
to the pair of rows $(p,p+1)$ in a $(k,p,q)$-matrix with $p<q$. It is not hard to see that
this operation yields a $(k,p+1,q)$-matrix. Let us now look in more detail at
the situation related to the hard case of the pseudoswap. Recall that if $X$ and $Y$ are two sets
of rows of $M$, we say that \emph{$M$ avoids $12112$ in $(X,Y)$}, if there
is no $x\in X$ and $y\in Y$ such that $x<y$ and the two rows $x,y$ contain a
copy of $12112$.

\begin{figure}
\includegraphics[scale=0.85]{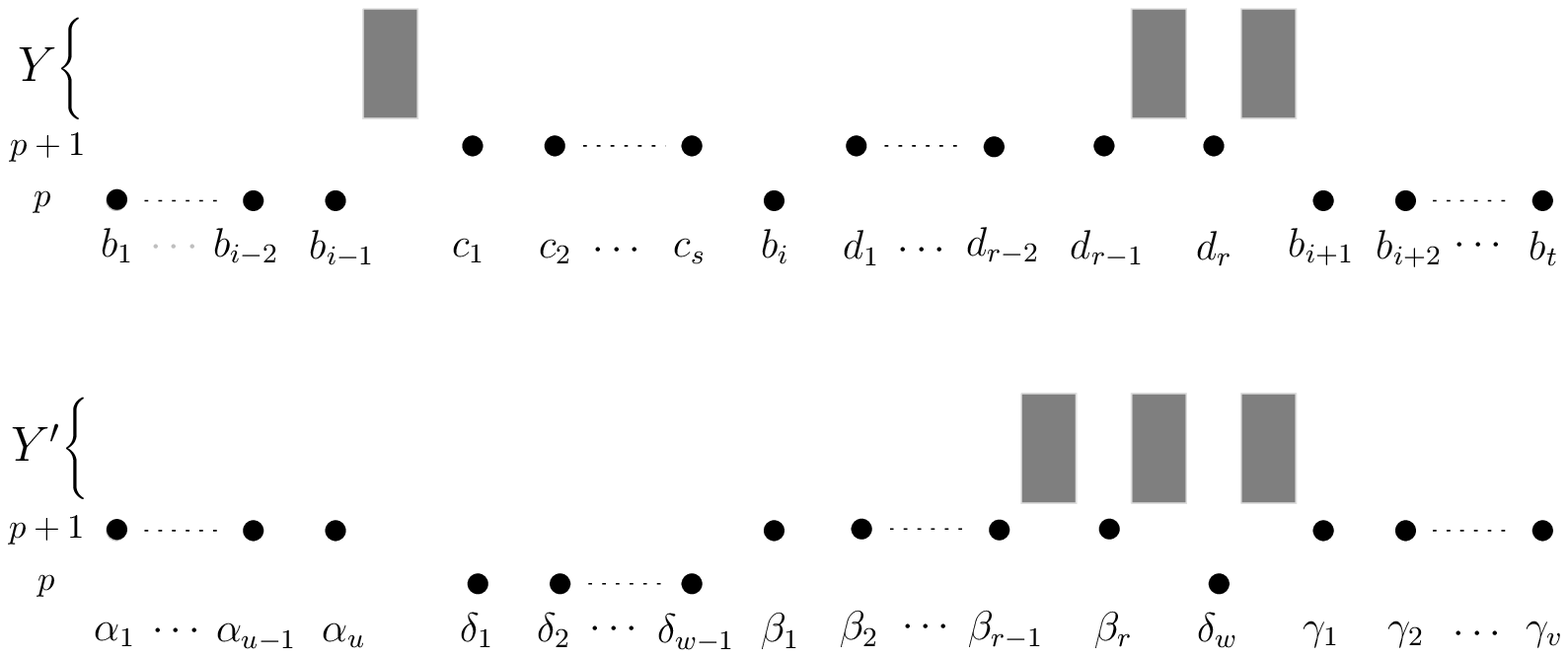}
\caption{Illustration of Lemma~\ref{lem-high}, part $(a)$ is above, part
$(b)$ below. The black dots correspond to 1-cells in rows $p$ and $p+1$, and
the shaded rectangles correspond to possible positions of the 1-cells in the
rows of $Y$ or $Y'$.}
\label{fig-l4}
\end{figure}

The following technical lemma is illustrated on Figure~\ref{fig-l4}.

\begin{lemma}\label{lem-high}\
\begin{itemize}
\item[(a)] Let $M$ be a $(k,p,q)$-matrix that has no copy of $12112$ in the two
rows $p,p+1$. Let $f_p(M)=b_1<b_2<\dotsb<b_t=l_p(M)$ be the $p$-columns of~$M$. Assume that the row
$p+1$ is separated by the column $b_i$, and that it has $r\ge 2$ rear 1-cells.
Let $c_1<c_2<\dotsb<c_s$ be the front $(p+1)$-columns and
let $d_1<d_2<\dotsb<d_r$ be the rear $(p+1)$-columns. By Observation~\ref{obs-f},
we have the inequalities
\[
b_1<\dotsb<b_{i-1}<c_1<\dotsb<c_s<b_i<d_1<\dotsb<d_r<b_{i+1}<\dotsb<b_t.
\]

Let $X=\{p,p+1\}$ and let
$Y$ be the set of all the rows above $p+1$ that contain at least one
1-cell to the left of the column $d_{r-1}$; formally,
\[
Y=\{y>p+1;\ f_y(M)<d_{r-1}\}.
\]
The matrix $M$ avoids $12112$ in $(X,Y)$ if and only if each $Y$-column $y$
satisfies one of the following three inequalities:
\begin{enumerate}
\item $b_{i-1}<y<c_1=f_{p+1}$
\item $d_{r-1}<y<d_r$
\item $d_r<y<b_{i+1}$
\end{enumerate}
The rows in $Y$ are precisely the rows above $p+1$ that are
separated by the $p$-column~$b_i$.
\item[(b)] Let $M'$ be a $(k,p+1,q)$-matrix that avoids $12112$ in
$(p,p+1)$. Let
$\alpha_1<\dotsb<\alpha_u<\beta_1<\dotsb<\beta_r<\gamma_1<\dotsb<\gamma_v$ be the
$(p+1)$-columns of $M'$, where the $\alpha_i$, $\beta_i$ and $\gamma_i$ denote respectively the left,
middle and right $(p+1)$-columns with respect to row $p$. Assume that there
are at least two middle 1-cells. Let $\delta_1<\dotsb<\delta_w$ be the $p$-columns of
$M'$. By Observation~\ref{obs-l}, we have the inequalities
\[
\alpha_1<\dotsb<\alpha_u<\delta_1<\dotsb<\delta_{w-1}<\beta_1<\dotsb<\beta_r<\delta_w<\gamma_1<\dotsb<\gamma_v.
\]

Let $X=\{p,p+1\}$ and let $Y'$ be the set of all the rows above $p+1$ that
contain at least one 1-cell to the left of column $\beta_r$. The matrix $M'$ avoids
$12112$ in $(X,Y')$ if and only if each $Y'$-column $y$ satisfies one of the
following three inequalities:
\begin{enumerate}
\item $\beta_{r-1}<y<\beta_r$
\item $\beta_r<y<\delta_w$
\item $\delta_w<y<\gamma_1$
\end{enumerate}
The rows in $Y'$ are precisely the rows above $p+1$ that are
separated by the $(p+1)$-column $\beta_r$.
\end{itemize}
\end{lemma}

\begin{proof}
Let us consider part $(a)$. Fix a row $y\in Y$. By the definition of
$(k,p,q)$-matrix, we have $d_r< l_y$. By Observation~\ref{obs-l}, we see
that $M$ avoids $12112$ in $(p+1,y)$ if and only if every $y$-column $j$
satisfies either $j<c_1=f_{p+1}$, $d_{r-1}<j<d_r$, or $j>d_r$. The first $y$-column satisfies
$f_y<d_{r-1}$ by the definition of $Y$, and hence $f_y<c_1=f_{p+1}<b_i$. Since
$l_y>l_{p+1}=d_r>b_i$, we see that if the pair of rows $(p+1,y)$ avoids $12112$,
then $y$ is separated by $b_i$ and the two rows $(p,y)$ avoid $12112$ if and
only if $b_{i-1}<f_y<l_y<b_{i+1}$. This proves part $(a)$ of the lemma.

The proof of part $(b)$ is analogous and we omit it.
\end{proof}

\subsection{The bijection}
We are now ready to present the bijection $\phi$. Let $M$ be a
$12112$-avoiding $(k,p,q)$-matrix with $p<q$, and let us write $f_i$ and $l_i$
for $f_i(M)$ and $l_i(M)$. By the
definition of $(k,p,q)$-matrix and by the assumption $p<q$, we know that
$f_p<f_{p+1} \le l_{p+1}<l_p$, so we may perform the pseudoswap of the rows $p$
and $p+1$ in~$M$. Let $M'$ be the $m\times n$
matrix obtained from $M$ by this pseudoswap. Let $f'_i=f_i(M')$ and $l'_i=l_i(M')$.
Note that $f'_i=f_i$ and $l'_i=l_i$ for every $i\not\in\{p,p+1\}$.

Clearly, $M'$ is a $(k,p+1,q)$-matrix. We now distinguish two cases,
depending on whether the pseudoswap we performed was easy or hard.

\textbf{Easy case:} If the row $p+1$ of $M$ has at most one rear 1-cell with
respect to row $p$, then $M'$ is $12112$-avoiding, and we may define
$\phi(M)=M'$. Indeed, from the definition of the pseudoswap we know that
$M'$ cannot contain a copy of $12112$ in the rows $(p,p+1)$, and
since we are performing the easy case of the pseudoswap, we cannot create any
new copy of the forbidden pattern that would intersect the remaining $m-2$ rows.

\textbf{Hard case:} Assume that the row $p+1$ of $M$
has $r>1$ rear 1-cells. Let $b_1<\dotsb<b_t$, $c_1<\dotsb<c_s$,
$d_1<d_2<\dotsb<d_r$, and $Y$ have the same meaning as in part $(a)$ of
Lemma~\ref{lem-high}. Let $Y_1$, $Y_2$ and $Y_3$ denote, respectively, the
$Y$-columns that lie between $b_{i-1}$ and $c_1$, between $d_{r-1}$ and $d_r$,
and between $d_r$ and $b_{i+1}$.

The bijection $\phi$ is now constructed in two steps. In
the first step, we perform the pseudoswap of the rows $p$ and $p+1$. Let $M'$ be
the result of this first step. Let us now apply the notation of part $(b)$ of
Lemma~\ref{lem-high} to the matrix $M'$; see Fig.~\ref{fig-l4}. Note that $d_{r-1}=\beta_r$, and hence
$Y=Y'$. Part $(b)$ of Lemma~\ref{lem-high} requires that all the $Y'$-columns of
a $12112$-avoiding $(k,p+1,q)$-matrix fall into one of the three groups:
\begin{itemize}
\item columns between $\delta_w<y<\gamma_1$. In $M'$, we have $\delta_w=d_r$ and
$\gamma_1=b_{i+1}$, so these columns are precisely the columns $Y_3$.
\item columns between $\beta_r<y<\delta_w$. In $M'$, these are precisely the
columns $Y_2$.
\item columns between $\beta_{r-1}<y<\beta_r$. In $M'$, there are no $Y$-columns
in this range.
\end{itemize}
On the other hand, if $Y_1$ is nonempty, then these columns violate the
inequalities of part $(b)$ in Lemma~\ref{lem-high}, showing that $M'$ is not
$12112$-avoiding. We will now apply the second
step of the bijection~$\phi$, which will exchange the relative order of the
columns $Y_1$ and some of the $\{p,p+1\}$-columns of $M'$. Intuitively, we will
move the columns $Y_1$ into the gap between $\beta_{r-1}$ and $\beta_r$, to
transform $M'$ into a matrix that satisfies the inequalities of
Lemma~\ref{lem-high}. This
operation only affects the columns $Y_1$ and the $\{p,p+1\}$-columns; all the
other cells of the matrix $M'$ remain unchanged.

Formally, the second step of the
bijection is performed as follows. Consider the submatrix of $M'$ induced by the columns
$Y_1$ and the columns
$Z=\{\delta_1<\dotsb<\delta_{w-1}<\beta_1<\dotsb<\beta_{r-1}\}$.
Note that the columns $Y_1$ are to the left of any column of $Z$. Now
we rearrange the columns inside this submatrix, such that the columns of the set
$Z$ precede the columns of $Y_1$. The relative order of the columns in $Z$, as
well as the relative order of the columns in $Y_1$, is preserved. This operation transforms the matrix $M'$ into
a matrix $M''$ that satisfies the conditions of part $(b)$ of
Lemma~\ref{lem-high}. We now define $\phi(M)=M''$.

Since $M''$ is clearly a $(k,p+1,q)$-matrix, it remains to check that $M''$
avoids $12112$. Let $x<y$ be a pair of rows of $M''$. We want to check that
$M''$ avoids $12112$ in these two rows. Let us consider the following cases
separately.

\textbf{The case $x<p$:} Since the rows below row $p$ are unaffected by
$\phi$, by Lemma~\ref{lem-low}, we know that $M''$ avoids $12112$ in
the rows $(x,y)$.

\textbf{The case $x=p, y=p+1$:} The properties of pseudoswap guarantee
that
 $M''$ avoids $12112$ in these two rows.

\textbf{The case $x\in X=\{p,p+1\}$ and $y\in Y'$:} By construction,
$M''$ satisfies the inequalities of part $(b)$ of
Lemma~\ref{lem-high}, and thus it avoids $12112$ in $(X,Y)$.

\textbf{The case $x\in X=\{p,p+1\},$ $y\not\in Y'$ and $y>p+1$:} By the
definition of $Y'$, we have $f_y(M'')=f_y(M)>d_{r-1}=\beta_r$. In
any column to the right of $\beta_r$ the mapping $\phi$ acts by
exchanging the rows $p$ and $p+1$. It is easy to check that this
action cannot create a copy of $12112$ in $(x,y)$ (note that in any
of the three matrices $M$, $M'$ and $M''$, both
the rows $p$ and $p+1$ have a 1-cell to the left of $\beta_r$ ).

\textbf{The case $y>x>p+1$:} The submatrix of $M''$ induced by the rows
above $p+1$ only differs from the corresponding submatrix of $M$ by
the position of the zero columns. Thus, it cannot contain any copy of
$12112$.

This shows that $\phi(M)$ is indeed a $12112$-avoiding $(k,p+1,q)$-matrix.

It is routine to check that the mapping $\phi$ can be inverted, which shows that
$\phi$ is indeed the required bijection.

\subsection{Consequences}

Theorem~\ref{thm-12112} has several consequences for pattern-avoiding
fillings of Ferrers shapes and pattern-avoiding ordered graphs.

By Lemma~\ref{lem-fil2}, there is a bijection
between $12112$-avoiding partitions of $[n]$ with $m$ blocks and
semi-standard fillings of Ferrers shapes with $n-m$ columns and at most $m$
rows that avoid $M(2,112)$; similarly, there is an analogous bijection
between $12212$-avoiding partitions and $M(2,212)$-avoiding fillings of
Ferrers shapes. Thus, we obtain the following direct consequence of Theorem~\ref{thm-12112}.

\begin{corollary}\label{cor-112}
For every $r$ and $c$, there is a bijection between the $M(2,112)$-avoiding
semi-standard fillings of all the Ferrers shapes with $r$ rows and $c$
columns and the $M(2,212)$-avoiding semi-standard fillings of all the Ferrers
shapes with $r$ rows and $c$ columns.
\end{corollary}

It would be tempting to assume that for a given Ferrers shape $F$, the
$M(2,112)$-avoiding semi-standard fillings of $F$ are in bijection with the
$M(2,212)$-avoiding semi-standard fillings of $F$, i.e., that the two
matrices $M(2,112)$ and $M(2,212)$ are Ferrers-equivalent. However, as we
already mentioned in the introduction of Section~\ref{sec-12112}, this is not the case. For instance,
the Ferrers shape $F$ with five columns of height 4 and one column of height
2 has $866$ $M(2,112)$-avoiding fillings but only $865$ $M(2,212)$-avoiding
fillings. Thus, the bijection of Corollary~\ref{cor-112} in general cannot
preserve the shape of the underlying diagram.

Let us now describe a well-known and useful correspondence between 0-1 fillings
of Ferrers shapes, and graphs with linearly ordered vertex sets; the
correspondence has been used , e.g., in \cite{adm} or \cite{kra}.

Every 0-1 filling
$F$ of a Ferrers shape with $p$ columns and $q$ rows can be represented by a
graph with $p+q$ linearly ordered vertices, defined in the following way: the
graph has two kinds of vertices, called \emph{right vertices $r_1,\dotsc,r_p$}
and \emph{left vertices $l_1,\dotsc,l_q$}. The $i$-th column of $F$ is
associated with the $i$-th right vertex $r_i$, and the $j$-th row of $F$ is
associated with the $j$-th left vertex~$l_j$. All the vertices are linearly
ordered by a left-to-right relation $<$ with the properties $r_1<\dotsb<r_p$,
$l_1<l_2<\dotsb < l_q$, and furthermore, $l_j<r_i$ if and only if row $j$ intersects
column $i$ inside $F$. The edge-set of the graph is determined by the 1-cells
of $F$ in the natural way: a 1-cell in row $j$ and column $i$ corresponds to
the edge between $l_j$ and $r_i$. Note that if $l_j$ and $r_i$ are connected
by an edge, then $l_j<r_i$.

\begin{figure}
\includegraphics[scale=0.75]{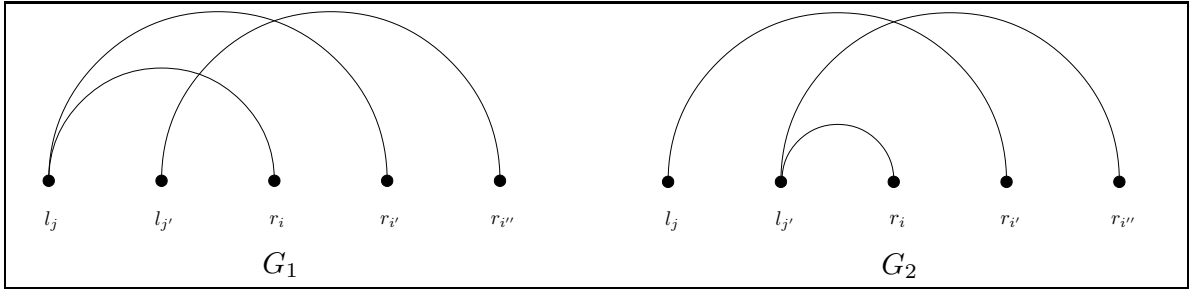}
\caption{The ordered graphs $G_1$ and $G_2$ corresponding to the filling
patterns $M(2,112)$ and $M(2,212)$.}\label{fig-g}
\end{figure}

In this representation, the semi-standard fillings of Ferrers shapes correspond precisely to the ordered
graphs with the property that every right vertex is connected to precisely
one left vertex, whereas the degrees of the left vertices can be arbitrary.
In accordance with our terminology for fillings, we will call such graphs
\emph{semi-standard}. Pattern-avoidance of semi-standard graphs has been
studied by A.~de Mier~\cite{adm2}, who considered the avoidance of crossings and
nestings with prescribed size within this class of graphs. The $M(2,112)$
avoiding fillings of $F$ correspond precisely to ordered graphs which avoid a
subgraph $G_1$ with five vertices $l_j<l_{j'}<r_i<r_{i'}<r_{i''}$ and three
edges $l_j r_i, l_j r_{i'},$ and $l_{j'} r_{i''}$. Similarly, the fillings
avoiding $M(2,212)$ correspond to graphs avoiding the subgraph $G_2$ with
vertices $l_j<l_{j'}<r_i<r_{i'}<r_{i''}$ and edges $l_j r_{i'}, l_{j'}
r_{i},$ and $l_{j'} r_{i''}$ (see Figure~\ref{fig-g}).

Theorem~\ref{thm-12112} then immediately yields
the following result.

\begin{corollary}
There is a bijection between semi-standard $G_1$-avoiding graphs and
semi-standard $G_2$-avoiding graphs that preserves the number of left vertices
and right vertices.
\end{corollary}

Whether this result can be extended to more general classes of graphs or more
general pairs of patterns is at this point an open problem.

\section{Concluding remarks}

In Appendix A, B and C, we present the results of the computer enumeration of
partitions avoiding fixed patterns of size five, six and seven,
respectively. Each row of the tables corresponds to one equivalence class.

In Table~\ref{tab-classes}, we present the total number of
equivalence classes of patterns of length $1,2,\dotsc,7$.

\begin{table}[htb]
\begin{center}
\begin{tabular}{|c|c|c|c|c|c|c|c|}
\hline $n$&1&2&3&4&5&6&7\\ \hline\hline number of classes of
patterns of size $n$ & 1&1&2&5&21&114&617 \\ \hline
\end{tabular}\\[4pt]
\caption{Number of equivalence classes of patterns of size
$1,2,\dotsc,7$}\label{tab-classes}
\end{center}
\end{table}

To provide an accurate asymptotic estimate of the number of equivalence
classes of patterns of a given size remains out of reach of our methods.

Let us remark that our computer enumeration has revealed several pairs of
non-equivalent patterns $\tau\not\sim\tau'$ whose growth functions
$p(n;\tau)$ and $p(n;\tau')$ coincide for several small values of~$n$. For
instance, the growth functions of the two patterns $\tau=1234415$ and
$\tau'=1234152$ are equal for $n<15$; in other words, the value of $n=15$ is
the smallest witness of the non-equivalence of the two patterns. It is an
interesting open problem to find, for a given~$l$, a common upper bound for all the
smallest witnesses demonstrating the non-equivalence of the non-equivalent
pairs of patterns of length $l$. Note that for any $k$ and for $\tau, \tau'$
chosen as above, the pair of non-equivalent patterns $12\dotsb k (\tau+k)$ and
$12\dotsb k(\tau'+k)$ of length $k+7$ requires a witness of length
$k+15$ (this follows from Theorem~\ref{thm-lift}).

\section*{Appendix A: Table of patterns of length five}

\footnotesize
\begin{longtable}{l|l}
\caption{Number of partitions in $P(n;\tau)$, where $\tau\in
P(5)$.}\label{tab5}\\
$\tau$ & $\{p(n;\tau)\}_{n=6}^{11}$ \\
\hline
\endhead
\hline
\multicolumn{2}{r@{}}%
continue\\
\endfoot
\endlastfoot
${12314}, {12324}, {12334}, {12341}, {12342}, {12343}, {12344}, {12345}$&$187,715,2795,11051,43947,175275$  \\
  \hline
$12313$&$188,730,2933,12061,50423,213423$ \\ \hline
$12323, 12234, 12332, 12123, 12132, 12213, 12231, 12312$&  \\
  $ 12321, 12331,  12134$&$188,731,2950,12235,51822,223191$  \\
  \hline
$12233, 12133$&$188,732,2969,12452,53769,238379$  \\ \hline
$11223, 11232$ &$189,746,3094,13371,59873,276670$  \\\hline
$11234$ &$189,747,3109,13507,60837,282503$ \\ \hline
$12131$&$189,747,3109,13517,61061,285503$ \\ \hline
$11233$ &$189,747,3111,13550,61393,288157$ \\ \hline
$12223, 12232, 12322, 12333, 12311, 12113$&$189,747,3111,13551,61419,288543$ \\ \hline
$11231$ &$190,760,3222,14350,66715,322218$ \\ \hline
$11213$ &$190,760,3223,14366,66882,323663$ \\ \hline
$11123$ &$191,771,3310,14969,70831,348887$ \\ \hline
$12112$, $12122$, $12212$, $12221$ &$191,773,3336,15207,72697,362447$  \\ \hline
$12121$ &$191,773,3337,15224,72892,364317$ \\ \hline
$11212$ &$191,773,3337,15224,72893,364341$ \\ \hline
$12211$ &$191,774,3351,15361,74043,373270$ \\ \hline
$11221$ &$191,774,3353,15393,74395,376556$ \\ \hline
$11222$ &$191,774,3354,15409,74579,378365$ \\ \hline
$11122$ &$191,774,3355,15424,74738,379805$ \\ \hline
$11112, 11121, 11211,12111, 12222$ &
  $192,789,3495,16545,83142,441009$ \\ \hline
$11111$ & $196,827,3795,18755,99146,556711$ \\ \hline
\end{longtable}
\normalsize

\section*{Appendix B: Table of patterns of length six}
\footnotesize
\begin{longtable}{l|l}
\caption{Number of partitions in $P(n;\tau)$, where $\tau\in
P(6)$.}\\
$\tau$ & $\{p(n;\tau)\}_{n\geq7}$ \\
\hline
\endhead
\hline
\multicolumn{2}{r@{}}%
continue\\
\endfoot
\endlastfoot
$123415, 123425, 123435, 123445, 123451, 123452$& \\
$123453, 123454, 123455, 123456$& $855, 3845, 18002, 86472$\\ \hline
$123414$&    $856, 3867,    18286,   89291$ \\ \hline
$123413,123424$&    $856, 3867,    18288,   89348$ \\ \hline
$123134, 123143$& $856, 3867,    18289,   89375, 447219,
2277477$\\ \hline
$123241$&    $856, 3868,    18312, 89684, 450407$\\ \hline
$123314$& $856, 3868, 18312, 89684, 450408, 2305592,$\\
        &\hfil $11978961,62983208$\\ \hline
$123142$& $856, 3868, 18312, 89684, 450408, 2305592,$\\
        &\hfil $11978961,62983209$\\ \hline
$123124,123145,123214,123234,123243,123245$  & \\
$123324,123341, 123342, 123345,123412,123421$& $856,
3868, 18313, 89711, 450825, 2310453$\\
$123423,123431,123432,123434,123441,123442$& \\
$123443$& \\ \hline

$123144,123244,123344$&$856,3869,    18340,   90135, 455917$\\
\hline
$121342$&$857, 3888,    18555,   92027$  \\ \hline
$122314$&$857, 3889,    18578,   92339$  \\ \hline
$122341$&$857, 3889,    18578,   92341$  \\ \hline
$121324$&$857, 3889,    18579,   92369$  \\ \hline

$121334,122334,121343,122343$&$857, 3890,18605, 92767, 478726, 2544145$\\
\hline

$121345,122345$&$857,3891,18628,93074,481845,2570867$  \\ \hline

$123141$&$857,3891,18628,93082$  \\ \hline
$123242$&$857,3891,18628,93084$  \\ \hline
$121344,122344$&$857,3891,18630,93135,482921,2585332$ \\ \hline

$123114,123224,123334,123343,123411,123422$\\
$123433,123444$&$857,3891,18630,93136$ \\ \hline

$121234,122134$&$858,3908, 18801, 94448, 491234, 2628572$ \\
\hline

$123132$&$858, 3909, 18821, 94686$ \\ \hline
$123213$&$858, 3909, 18822, 94712$ \\ \hline
$122313$&$858, 3910, 18844, 95008$\\ \hline
$121332$&$858, 3910, 18845, 95037$ \\ \hline

$123123,123312,123321$&$858,3910,18846,95058$ \\ \hline

$123231$&$858 ,3910    ,18847   ,95086$\\ \hline $121323$&$858 ,3910
,18847   ,95087$\\ \hline $122331$&$858 ,3911    ,18871   ,95434$\\
\hline $121341$&$858 ,3911    ,18872   ,95455$\\ \hline
$121314$&$858 ,3911    ,18872   ,95460$\\ \hline
$122133,121233$&$858 ,3911    ,18872   ,95461,   523161$\\ \hline
$112342$&$858 ,3911
,18873   ,95485$\\ \hline $122342$&$858 ,3911    ,18873   ,95486$\\
\hline $122324$&$858,3911,18874,95511$\\ \hline
$112324$&$858,3911,18874,95513$\\ \hline

$112334,112343$&$858,3912,18897,95828,506812,2781704$\\ \hline
$112345$&$858,3912,18900,95904$\\ \hline
$112344$&$858,3912,18900,95909$\\ \hline

$121134,122234$&$859,3929, 19077, 97377,   518804,  2869604$ \\ \hline

$112234$&$859, 3930, 19096, 97599$ \\ \hline
$112323$&$859, 3930, 19100, 97700, 522417$ \\ \hline
$112332$&$859, 3930, 19100, 97700, 522415$ \\ \hline
$123131$&$859, 3931, 19115, 97828, 523161$ \\ \hline
$123313$&$859, 3931, 19115, 97828, 523144$ \\ \hline
$123133$&$859, 3931, 19115, 97831$ \\ \hline
$123113$&$859, 3931, 19116, 97852$ \\ \hline
$121313$&$859, 3931, 19117, 97872$ \\ \hline
$121132$&$859, 3931, 19117, 97882$ \\ \hline
$121312$&$859, 3931, 19118, 97898$ \\ \hline

$121223, 121232, 121322, 122123, 122132, 122213$&\\
$122231, 122312, 122321, 123112, 123122, 123212$&\\
$123221, 123223, 123233, 123323, 123331, 123332$&$859,3931,19119
97921,524460,2921730$\\ \hline

$121123$& $859,3931,19120,97945$\\ \hline

$123121,123232$&$859,3931,19120,97947,524870,2926845$\\\hline
$121213$&$859,3931,19120,97947,524870,2926847$\\ \hline
$122323$&$859,3931,19120,97947,524871$\\ \hline

$121231$&$859,3931,19120,97948$\\\hline
$121321$&$859,3931,19121,97972$\\\hline
$112341$&$859,3931,19122,97987$\\\hline
$112314$&$859,3931,19122,97992$\\\hline
$112233$&$859,3931,19123,98023$\\\hline
$122131$&$859,3932,19139,98173$\\\hline
$122113$&$859,3932,19141,98222$\\\hline
$121331$&$859,3932,19142,98242$\\\hline
$122311,123311,123211,123322$&$859,3932,19142,98246,528141,2958634$\\\hline
$122332$&$859,3932,19144,98296$\\\hline

$121333,122333$&$859,3932,19145,98321,529292,2972760$\\\hline
$121133,122233$&$859,3932,19146,98345,529646,2976939$\\\hline
$112312$&$860,3948,19308,99685$\\\hline
$112123$&$860,3948,19310,99730   540195$\\\hline
$112132$&$860,3948,19310,99730   540193$\\\hline
$112134$&$860,3949,19327,99908$\\\hline
$112313$&$860,3949,19327,99914$\\\hline
$112321$&$860,3949,19330,99990$\\\hline

$112223,112232,112322$&$860,3949,19330,99993,543077,3081145$\\\hline

$112213$&$860,3949,19331,100010$\\\hline
$112231$&$860,3949,19332,100031$\\\hline
$112331$&$860,3950,19350,100240$\\\hline
$112133$&$860,3950,19354,100332$\\\hline
$112333$&$860,3950,19354,100338$\\\hline
$121131$&$860,3954,19434,101338$\\\hline
$121311$&$860,3954,19434,101342$\\\hline

$121113,122232,122322,123111,123222,123333$&\\
$122223$&$860,3954,19434,101350,557570,3220754$\\\hline

$111223,111232$&$861,3964,19488,101434,555332,3181699$\\\hline

$111234$&$861,3966,19516,101662$\\\hline
$111233$&$861,3966,19523,101837$\\\hline
$112131$&$861,3970,19599,102778$\\\hline
$112113$&$861,3970,19599,102786$\\\hline
$112311$&$861,3970,19600,102802$\\\hline
$111213$&$862,3984,19731,103869$\\\hline
$111231$&$862,3984,19733,103905$\\\hline
$111123$&$863,3996,19837,104726$\\\hline
$121212$&$863,3999,19880,105134,  587479,  3449505$\\\hline
$122121$&$863,3999,19880,105134,  587479,  3449509$\\\hline
$121221$&$863,3999,19880,105135$\\\hline
$122112$&$863,3999,19881,105150$\\\hline
$112122$&$863,3999,19882,105176$\\\hline
$121122$&$863,3999,19883,105188$\\\hline
$112212$&$863,3999,19883,105192$\\\hline
$122211$&$863,3999,19885,105226$\\\hline
$112221$&$863,3999,19885,105233$\\\hline
$111222$&$863,3999,19889,105314$\\\hline

$122221,121222,122122,122212$&$863,4001,19917,105594$\\\hline

$112112$&$863, 4001, 19918,105614,  592676$\\\hline
$121112$&$863, 4001, 19918,105614,  592671$\\\hline
$121121$&$863, 4001, 19918,105618$\\\hline
$111212$&$863, 4001, 19919,105636,  592976,  3504921$\\\hline
$112121$&$863, 4001, 19919,105636,  592976,  3504918$\\\hline
$121211$&$863, 4001, 19919,105636,  592976,  3504916$\\\hline
$112222$&$863, 4002, 19938,105878$\\\hline
$122111$&$863, 4002, 19939,105886$\\\hline
$111221$&$863, 4002, 19939,105893$\\\hline
$112211$&$863, 4002, 19939,105895$\\\hline
$111122$&$863, 4002, 19939,105901$\\\hline

$111112,111121,111211,112111,121111,122222$&$864,4020,20150,107964$\\\hline
$111111$&$869,4075,20645,112124$\\\hline

\end{longtable}
\normalsize

\section*{Appendix C: Table of patterns of length seven}
{
\footnotesize
\begin{longtable}{l|l}
\caption{Number of partitions in $P(n;\tau)$, where $\tau\in
P(7)$.}\\
$\tau$ & $\{p(n;\tau)\}_{n\geq8}$ \\
\hline
\endhead
\hline
\multicolumn{2}{r@{}}%
continue\\
\endfoot
\endlastfoot
$1234516,1234526,1234536,1234546$&\\
$1234556,1234561,1234562,1234563$&\\
$1234564,1234565,1234566,1234567$&$4111,20648,109299,601492$\\\hline

$1234515$& $4112,20678,109817,608258$\\\hline

$1234514,1234525$& $4112,20678,109817,608261$\\\hline

$1234154$&$4112,20678,109818,608300,3478443$\\\hline

$1234145$&$4112,20678,109818,608300,3478444$\\\hline

$1234513,1234524,1234535$&$4112,20678,109819,608338,3479249$\\\hline

$1234153$&$4112,20678,109819,608338,3479251$\\\hline

$1234135,1234245,1234254$&$4112,20678,109820,608375$\\\hline

$1234251$&$4112,20679,109852,608957,3487954,20485468$\\\hline

$1234415$&$4112,20679,109852,608957,3487954,20485475,$\\
         &\hfil $122666770, 745713106$\\\hline
$1234152$&$4112,20679,109852,608957,3487954,20485475,$\\
         &\hfil $122666770, 745713111$\\\hline
$1234351,1234352$&$4112,20679,109852,608959,3488036,20487341$\\\hline

$1234315,1234425$&$4112,20679,109852,608959,3488037,20487386,$\\
                 &\hfil $122699078,746161492$\\ \hline
$1234253$&$4112,20679,109852,608959,3488037,20487386,$\\
         &\hfil $122699078,746161493$\\ \hline

$1234125,1234156,1234215,1234235$&\\
$1234256,1234325,1234345,1234354$&\\
$1234356,1234435,1234451,1234452$&\\
$1234453,1234456,1234512,1234521$&\\
$1234523,1234531,1234532,1234534$&\\
$1234541,1234542,1234543,1234545$&\\
$1234551,1234552,1234553,1234554$&$4112,20679,109853,608996,3488806$\\\hline

$1234155,1234255,1234355,1234455$&$4112,20680,109889,609735$\\\hline

$1231453$&$4113,20707,110309,614684$\\\hline

$1231452$&$4113,20707,110311,614747$\\\hline

$1232453$&$4113,20707,110311,614752$\\\hline

$1231435$&$4113,20707,110313,614824$\\\hline

$1232415$&$4113,20708,110342,615304,3558058$\\\hline

$1232451$&$4113,20708,110342,615304,3558067$\\\hline

$1231425$&$4113,20708,110342,615306$\\\hline

$1233415$&$4113,20708,110343,615337$\\\hline

$1233425$&$4113,20708,110343,615339$\\\hline

$1233451,1233452$&$4113,20708,110343,615341$\\\hline

$1232435$&$4113,20708,110344,615379$\\\hline

$1231445,1231454,1232445,1232454$&\\
$1233445,1233454$&$4113,20709,110379,616082$\\\hline

$1231456,1232456,1233456$&$4113,20710,110411,616664$\\\hline

$1234151$&$4113,20710,110411,616672,3578613$\\\hline

$1234252$&$4113,20710,110411,616672,3578615$\\\hline

$1234353$&$4113,20710,110411,616674$\\\hline

$1231455,1232455,1233455$&$4113,20710,110413,616745$\\\hline

$1234115,1234225,1234335,1234445$&\\
$1234454,1234511,1234522,1234533$&\\
$1234544,1234555$&$4113,20710,110413,616746$\\\hline

$1231345$&$4114,20734,110743,620127$\\\hline

$1233145$&$4114,20735,110772,620616$\\\hline

$1231245,1232145,1232345,1233245$&$4114,20735,110773,620653$\\\hline

$1234214$&$4114,20736,110799,621037$\\\hline

$1234143$&$4114,20736,110800,621064$\\\hline

$1234142$&$4114,20736,110800,621065$\\\hline

$1232413$&$4114,20736,110800,621066$\\\hline

$1231342$&$4114,20736,110800,621070$\\\hline

$1234314$&$4114,20736,110802,621134$\\\hline

$1234132,1234243$&$4114,20736,110802,621136$\\\hline

$1232143$&$4114,20736,110802,621137$\\\hline

$1231432$&$4114,20736,110802,621138$\\\hline

$1234124$&$4114,20736,110802,621145$\\\hline

$1234134$&$4114,20736,110803,621172,3622245$\\\hline

$1234213,1234324$&$4114,20736,110803,621172,3622246$\\\hline

$1231324$&$4114,20736,110803,621173$\\\hline

$1232134$&$4114,20736,110804,621207$\\\hline

$1234241$&$4114,20736,110804,621214$\\\hline

$1231243$&$4114,20736,110805,621243,3623689$\\\hline

$1213452$&$4114,20736,110805,621243,3623710$\\\hline

$1231424$&$4114,20736,110805,621246$\\\hline

$1231434$&$4114,20736,110807,621319$\\\hline

$1213425$&$4114,20736,110808,621351$\\\hline

$1232414$&$4114,20737,110833,621694$\\\hline

$1234413$&$4114,20737,110833,621698$\\\hline

$1233414$&$4114,20737,110834,621729$\\\hline

$1233142$&$4114,20737,110834,621730$\\\hline

$1233424$&$4114,20737,110834,621733$\\\hline

$1231443$&$4114,20737,110834,621737$\\\hline

$1233241$&$4114,20737,110835,621766$\\\hline

$1231423$&$4114,20737,110835,621767$\\\hline

$1233124$&$4114,20737,110835,621768,3630754$\\\hline

$1231442$&$4114,20737,110835,621768,3630761$\\\hline

$1232443$&$4114,20737,110835,621772$\\\hline

$1233412,1233421,1234123,1234234$&\\
$1234312,1234321,1234412,1234421$&\\
$1234423,1234431,1234432$&$4114,20737,110836,621803,3631456,21211085$\\\hline

$1231234,1233214$&$4114,20737,110836,621803,3631456,21933850$\\\hline

$1232314$&$4114,20737,110836,621804$\\\hline

$1232341$&$4114,20737,110837,621841,3632280$\\\hline

$1234231,1234341,1234342$&$4114,20737,110837,621841,3632281$\\\hline

$1232431$&$4114,20737,110837,621842,3632324$\\\hline

$1232434$&$4114,20737,110837,621842,3632325$\\\hline

$1223451$&$4114,20737,110838,621872$\\\hline

$1223415$&$4114,20737,110838,621873$\\\hline

$1231344$&$4114,20737,110841,621987$\\\hline

$1232441$&$4114,20738,110868,622399$\\\hline

$1233441,1233442$&$4114,20738,110870,622474$\\\hline

$1231451$&$4114,20738,110871,622504$\\\hline

$1232452$&$4114,20738,110871,622505$\\\hline

$1233144$&$4114,20738,110871,622507$\\\hline

$1231415$&$4114,20738,110871,622508$\\\hline

$1232425$&$4114,20738,110871,622510$\\\hline

$1231244,1232144,1232344,1233244$&$4114,20738,110871,622511$\\\hline

$1213453,1223453$&$4114,20738,110872,622545$\\\hline

$1233453$&$4114,20738,110872,622546$\\\hline

$1233435$&$4114,20738,110873,622581$\\\hline

$1213435,1223435$&$4114,20738,110873,622583$\\\hline

$1213445,1213454,1223445,1223454$&$4114,20739,110905,623173$\\\hline

$1213456,1223456$&$4114,20739,110908,623279$\\\hline

$1213455,1223455$&$4114,20739,110908,623284$\\\hline

$1213423$&$4115,20762,111212,626275$\\\hline

$1223413$&$4115,20763,111238,626659$\\\hline

$1223143$&$4115,20763,111239,626702$\\\hline

$1213424$&$4115,20763,111239,626706$\\\hline

$1213432$&$4115,20763,111240,626729$\\\hline

$1213342$&$4115,20763,111240,626730$\\\hline

$1223134$&$4115,20763,111240,626733$\\\hline

$1213324$&$4115,20763,111240,626737,3684077$\\\hline

$1213243$&$4115,20763,111243,626837,3686012$\\\hline

$1213234$&$4115,20763,111243,626837,3686013$\\\hline

$1223145$&$4115,20764,111269,627221$\\\hline

$1223414$&$4115,20764,111270,627257$\\\hline

$1213442$&$4115,20764,111270,627260$\\\hline

$1213245$&$4115,20764,111270,627261$\\\hline

$1223431$&$4115,20764,111272,627332$\\\hline

$1223314$&$4115,20764,111273,627365$\\\hline

$1223341$&$4115,20764,111274,627402$\\\hline

$1231145,1232245,1233345$&$4115,20764,111274,627407$\\\hline

$1213345,1223345$&$4115,20765,111302,627864 $\\\hline

$1223441$&$4115,20765,111303,627898$\\\hline

$1223144$&$4115,20765,111306,627998$\\\hline

$1213244$&$4115,20765,111306,628002$\\\hline

$1213443,1223443$&$4115,20765,111306,628005,3074775$\\\hline

$1213434,1223434$&$4115,20765,111306,628005,3074777$\\\hline

$1123452$&$4115,20765,111311,628168$\\\hline

$1123425$&$4115,20765,111311,628173$\\\hline

$1234141$&$4115,20766,111330,628320$\\\hline

$1234414$&$4115,20766,111330,628321$\\\hline

$1234144$&$4115,20766,111330,628322$\\\hline

$1231413$&$4115,20766,111330,628326$\\\hline

$1231314$&$4115,20766,111330,628327$\\\hline

$1234313,1234424,1233413$&$4115,20766,111330,628328,3705907$\\\hline

$1234131,1234242$&$4115,20766,111330,628328,3705924$\\\hline

$1233134,1233143$&$4115,20766,111330,628329,3705940,22709849$\\\hline

$1234133,1234244$&$4115,20766,111330,628331$\\\hline

$1231343$&$4115,20766,111330,628333$\\\hline

$1231334,1231433$&$4115,20766,111330,628335,3706167,22714756$\\\hline

$1231341$&$4115,20766,111331,628354$\\\hline

$1234114$&$4115,20766,111331,628355$\\\hline

$1231431$&$4115,20766,111331,628358$\\\hline

$1231143$&$4115,20766,111331,628361$\\\hline

$1234113,1234224$&$4115,20766,111331,628362$\\\hline

$1231242,1232142$&$4115,20766,111331,628364,3706656,22721090$\\\hline

$1231414$&$4115,20766,111332,628387$\\\hline

$1232424$&$4115,20766,111332,628392$\\\hline

$1232412,1232421$&$4115,20766,111332,628395,3707209,22728608$\\\hline

$1231142$&$4115,20766,111332,628397$\\\hline

$1231134$&$4115,20766,111332,628400$\\\hline

$1232243$&$4115,20766,111332,628402$\\\hline

$1231412$&$4115,20766,111333,628426$\\\hline

$1232423$&$4115,20766,111333,628428$\\\hline

$1231241$&$4115,20766,111333,628435$\\\hline

$1232241$&$4115,20766,111334,628454$\\\hline

$1233314$&$4115,20766,111334,628456$\\\hline

$1231422$&$4115,20766,111334,628457$\\\hline

$1231224,1232124,1232214,1232334$&\\
$1232343,1232433,1233234,1233243$&\\
$1233324,1233341,1233342,1233423$&\\
$1233431,1233432,1234112,1234122$&\\
$1234212,1234221,1234223,1234233$&\\
$1234323,1234331,1234332,1234334$&\\
$1234344,1234434,1234441,1234442$&\\
$1234443$&$4115,20766,111334,628461$\\\hline

$1231124,1232234$&$4115,20766,111335,628495$\\\hline

$1234121,1234232,1234343$&$4115,20766,111335,628497,3709217,22758862$\\\hline

$1231214,1232324$&$4115,20766,111335,628497,3709217,22758864$\\\hline

$1233434$&$4115,20766,111335,628497,3709218$\\\hline

$1232342$&$4115,20766,111335,628498$\\\hline

$1231421$&$4115,20766,111335,628499$\\\hline

$1213451$&$4115,20766,111336,628523,3709547$\\\hline

$1213415$&$4115,20766,111336,628523,3709569$\\\hline

$1232432$&$4115,20766,111336,628532$\\\hline

$1223452$&$4115,20766,111337,628557$\\\hline

$1223425$&$4115,20766,111337,628562$\\\hline

$1213344,1223344$&$4115,20766,111338,628603$\\\hline

$1123453$&$4115,20766,111342,628731$\\\hline

$1123435$&$4115,20766,111342,628735$\\\hline

$1123445,1123454$&$4115,20766,111345,628837$\\\hline

$1123456$&$4115,20766,111345,628846,3716242$\\\hline

$1123455$&$4115,20766,111345,628846,3716256$\\\hline

$1233141$&$4115,20767,111363,628953$\\\hline

$1232141$&$4115,20767,111363,628954$\\\hline

$1233242$&$4115,20767,111363,628958$\\\hline

$1233114$&$4115,20767,111365,629024$\\\hline

$1232114,1233224$&$4115,20767,111365,629027$\\\hline

$1231441$&$4115,20767,111366,629055$\\\hline

$1232411$&$4115,20767,111366,629057,3717013$\\\hline

$1232442$&$4115,20767,111366,629057,3717017$\\\hline

$1233411,1233422,1234211,1234311$&\\
$1234322,1234411,1234422,1234433$&$4115,20767,111366,629061$\\\hline

$1233443$&$4115,20767,111368,629131$\\\hline

$1231444,1232444,1233444$&$4115,20767,111369,629166$\\\hline

$1231144,1232244,1233344$&$4115,20767,111370,629200$\\\hline

$1212334,1212343,1221334,1221343$&$4116,20788,111626,631531$\\\hline

$1212345,1221345$&$4116,20790,111675,632216$\\\hline

$1212344,1221344$&$4116,20790,111682,632466$\\\hline

$1123424$&$4116,20790,111686,632609$\\\hline

$1213422$&$4116,20791,111706,632807$\\\hline

$1213242$&$4116,20791,111707,632831$\\\hline

$1213412$&$4116,20791,111708,632868$\\\hline

$1213142$&$4116,20791,111709,632895$\\\hline

$1213224$&$4116,20791,111709,632900$\\\hline

$1213421$&$4116,20791,111709,632907$\\\hline

$1223142$&$4116,20791,111710,632925,3754896$\\\hline

$1213413$&$4116,20791,111710,632925,3754918$\\\hline

$1223412$&$4116,20791,111710,632929$\\\hline

$1223421$&$4116,20791,111710,632932$\\\hline

$1213124$&$4116,20791,111710,632933$\\\hline

$1223241$&$4116,20791,111711,632962$\\\hline

$1223214$&$4116,20791,111711,632964,3755736$\\\hline

$1223124$&$4116,20791,111711,632964,3755740$\\\hline

$1213143$&$4116,20791,111711,632968$\\\hline

$1213134$&$4116,20791,111711,632969$\\\hline

$1223423$&$4116,20791,111712,633000$\\\hline

$1123342$&$4116,20791,111712,633005$\\\hline

$1123423$&$4116,20791,111712,633007$\\\hline

$1123432$&$4116,20791,111712,633008$\\\hline

$1123324$&$4116,20791,111712,633009,3756838$\\\hline

$1123243$&$4116,20791,111712,633009,3756839$\\\hline

$1123234$&$4116,20791,111712,633009,3756840$\\\hline

$1213241$&$4116,20791,111713,633029$\\\hline

$1213214$&$4116,20791,111713,633032$\\\hline

$1223243$&$4116,20791,111714,633065,3757738$\\\hline

$1223234$&$4116,20791,111714,633065,3757740$\\\hline

$1123442$&$4116,20791,111714,633070$\\\hline

$1123245$&$4116,20791,111715,633098$\\\hline

$1223141$&$4116,20792,111737,633365$\\\hline

$1213145$&$4116,20792,111738,633395$\\\hline

$1223114$&$4116,20792,111739,633435$\\\hline

$1213414$&$4116,20792,111740,633454$\\\hline

$1223245$&$4116,20792,111740,633458$\\\hline

$1223424$&$4116,20792,111740,633464$\\\hline

$1223411$&$4116,20792,111740,633468$\\\hline

$1213314$&$4116,20792,111742,633529$\\\hline

$1213431$&$4116,20792,111743,633561$\\\hline

$1213341$&$4116,20792,111743,633562$\\\hline

$1223432$&$4116,20792,111743,633570$\\\hline

$1213334,1213343,1213433,1223334$&\\
$1223343,1223433$&$4116,20792,111743,633573$\\\hline

$1123345$&$4116,20792,111744,633595$\\\hline

$1223324$&$4116,20792,111744,633600$\\\hline

$1123443$&$4116,20792,111744,633605,3765380$\\\hline

$1123434$&$4116,20792,111744,633605,3765382$\\\hline

$1223342$&$4116,20792,111745,633631$\\\hline

$1123244$&$4116,20792,111745,633640$\\\hline

$1123451$&$4116,20792,111749,633753,3767900$\\\hline

$1123415$&$4116,20792,111749,633753,3767916$\\\hline

$1213441$&$4116,20793,111772,634059$\\\hline

$1223442$&$4116,20793,111772,634065$\\\hline

$1213144$&$4116,20793,111774,634137$\\\hline

$1123344$&$4116,20793,111774,634138$\\\hline

$1223244$&$4116,20793,111776,634197$\\\hline

$1213444,1223444$&$4116,20793,111776,634203$\\\hline

$1231141$&$4116,20797,111892,636183,3800334$\\\hline

$1232242$&$4116,20797,111892,636183,3800350$\\\hline

$1231411$&$4116,20797,111892,636187,3800468$\\\hline

$1232422$&$4116,20797,111892,636187,3800476$\\\hline

$1231114,1232224,1233334,1233343$&\\
$1233433,1234111,1234222,1234333$&\\
$1234444$&$4116,20797,111892,636195$\\\hline

$1221342$&$4117,20814,112031,636507$\\\hline

$1212342$&$4117,20814,112033,636564$\\\hline

$1221324$&$4117,20814,112033,636567,3791466$\\\hline

$1212324$&$4117,20814,112033,636567,3791468$\\\hline

$1211342$&$4117,20814,112034,636593$\\\hline

$1212314$&$4117,20814,112034,636598$\\\hline

$1211324$&$4117,20814,112034,636610$\\\hline

$1222314$&$4117,20814,112035,636623$\\\hline

$1212341$&$4117,20814,112036,636655$\\\hline

$1222341$&$4117,20814,112037,636686$\\\hline

$1221314$&$4117,20815,112059,636978$\\\hline

$1221341$&$4117,20815,112061,637033$\\\hline

$1211334,1211343,1222334,1222343$&$4117,20815,112072,637389$\\\hline

$1122334,1122343$&$4117,20816,112092,637612$\\\hline

$1123412$&$4117,20816,112094,637657$\\\hline

$1123142$&$4117,20816,112095,637686$\\\hline

$1123124$&$4117,20816,112095,637694$\\\hline

$1211345,1222345$&$4117,20817,112118,637987$\\\hline

$1123242$&$4117,20817,112119,638047$\\\hline

$1122345$&$4117,20817,112120,638059$\\\hline

$1123413$&$4117,20817,112120,638067$\\\hline

$1123143$&$4117,20817,112121,638103$\\\hline

$1123134$&$4117,20817,112121,638105$\\\hline

$1123422$&$4117,20817,112121,638107$\\\hline

$1123224$&$4117,20817,112121,638109$\\\hline

$1122344$&$4117,20817,112122,638140$\\\hline

$1123421$&$4117,20817,112123,638166$\\\hline

$1123214$&$4117,20817,112123,638172$\\\hline

$1123241$&$4117,20817,112124,638199$\\\hline

$1211344,1222344$&$4117,20817,112125,638232$\\\hline

$1123145$&$4117,20817,112126,638255$\\\hline

$1123314$&$4117,20818,112150,638607$\\\hline

$1123431$&$4117,20818,112151,638635$\\\hline

$1123341$&$4117,20818,112151,638640$\\\hline

$1123334,1123343,1123433$&$4117,20818,112152,638673$\\\hline

$1123441$&$4117,20818,112155,638729$\\\hline

$1123444$&$4117,20818,112155,638779$\\\hline

$1123144$&$4117,20818,112156,638802$\\\hline

$1212234,1221234,1222134$&$4118,20835,112311,639591$\\\hline

$1212134$&$4118,20835,112312,639622$\\\hline

$1211234$&$4118,20835,112312,639624$\\\hline

$1221134$&$4118,20836,112338,640044$\\\hline

$1213141$&$4117,20821,112237,640107$\\\hline

$1223242$&$4117,20821,112237,640113,3839906$\\\hline

$1213114$&$4117,20821,112237,640113,3839908$\\\hline

$1223224$&$4117,20821,112237,640121$\\\hline

$1213411$&$4117,20821,112238,640143$\\\hline

$1223422$&$4117,20821,112238,640147$\\\hline

$1121342$&$4118,20839,112419,641372$\\\hline

$1121324$&$4118,20839,112423,641478$\\\hline

$1231332$&$4118,20840,112437,641524$\\\hline

$1233132$&$4118,20840,112437,641528,3845007$\\\hline

$1213132$&$4118,20840,112437,641528,3845023$\\\hline

$1231312$&$4118,20840,112437,641528,3845028$\\\hline

$1232313$&$4118,20840,112437,641529$\\\hline

$1232213$&$4118,20840,112437,641530$\\\hline

$1223213$&$4118,20840,112437,641537$\\\hline

$1232133$&$4118,20840,112438,641563$\\\hline

$1231132$&$4118,20840,112439,641589$\\\hline

$1231323$&$4118,20840,112439,641590$\\\hline

$1232113$&$4118,20840,112439,641591$\\\hline

$1231322$&$4118,20840,112439,641592,3846251$\\\hline

$1233213$&$4118,20840,112439,641592,3846236$\\\hline

$1232131$&$4118,20840,112439,641595$\\\hline

$1213213$&$4118,20840,112439,641601$\\\hline

$1231223$&$4118,20840,112440,641620$\\\hline

$1232123$&$4118,20840,112440,641622$\\\hline

$1213123$&$4118,20840,112440,641623,3846806$\\\hline

$1231232$&$4118,20840,112440,641623,3846810$\\\hline

$1223123$&$4118,20840,112440,641623,3846811$\\\hline

$1231213$&$4118,20840,112440,641623,3846814$\\\hline

$1231321$&$4118,20840,112440,641623,3846829$\\\hline

$1232132$&$4118,20840,112440,641624$\\\hline

$1223132$&$4118,20840,112440,641625$\\\hline

$1213323$&$4118,20840,112440,641634$\\\hline

$1231233$&$4118,20840,112441,641652,3847223$\\\hline

$1233123$&$4118,20840,112441,641652,3847318$\\\hline

$1213223$&$4118,20840,112441,641653$\\\hline

$1233231$&$4118,20840,112441,641654,3847386$\\\hline

$1213232$&$4118,20840,112441,641654,3847388$\\\hline

$1212313$&$4118,20840,112441,641654,3847397$\\\hline

$1232331$&$4118,20840,112441,641655,3847421,24108094$\\\hline

$1232231$&$4118,20840,112441,641655,3847421,24108095$\\\hline

$1231123$&$4118,20840,112442,641682,3847867$\\\hline

$1233312,1233321$&$4118,20840,112442,641682,3847869$\\\hline

$1232312$&$4118,20840,112442,641684,3847930$\\\hline

$1232321$&$4118,20840,112442,641684,3847931$\\\hline

$1221323$&$4118,20840,112442,641684,3847935$\\\hline

$1223231$&$4118,20840,112442,641684,3847937$\\\hline

$1212323$&$4118,20840,112442,641684,3847939$\\\hline

$1231231$&$4118,20840,112442,641686$\\\hline

$1211323$&$4118,20840,112442,641696$\\\hline

$1222313$&$4118,20840,112443,641708$\\\hline

$1213332$&$4118,20840,112443,641712$\\\hline

$1213231$&$4118,20840,112444,641749,3849207$\\\hline

$1213233$&$4118,20840,112444,641749,3849214$\\\hline

$1122342$&$4118,20840,112445,641778$\\\hline

$1122324$&$4118,20840,112446,641807$\\\hline

$1122314$&$4118,20840,112448,641860$\\\hline

$1122341$&$4118,20840,112448,641868$\\\hline

$1223131$&$4118,20841,112468,642087,3852874$\\\hline

$1221313$&$4118,20841,112468,642087,3852877$\\\hline

$1213312$&$4118,20841,112468,642092$\\\hline

$1223133$&$4118,20841,112468,642095$\\\hline

$1223313$&$4118,20841,112468,642098$\\\hline

$1213321$&$4118,20841,112468,642099$\\\hline

$1223113$&$4118,20841,112469,642117$\\\hline

$1213322$&$4118,20841,112469,642124$\\\hline

$1211332$&$4118,20841,112469,642126$\\\hline

$1233122,1233212,1233221$&$4118,20841,112469,642128,3853761,24182852$\\\hline

$1233112$&$4118,20841,112469,642128,3853761,24182855$\\\hline

$1233121$&$4118,20841,112469,642131$\\\hline

$1212332$&$4118,20841,112469,642132$\\\hline

$1221332$&$4118,20841,112470,642155$\\\hline

$1223312$&$4118,20841,112471,642189$\\\hline

$1212331$&$4118,20841,112471,642191$\\\hline

$1223321$&$4118,20841,112471,642192$\\\hline

$1223331$&$4118,20841,112472,642222$\\\hline

$1233211$&$4118,20841,112472,642224$\\\hline

$1232311$&$4118,20841,112472,642228$\\\hline

$1222331$&$4118,20841,112473,642250$\\\hline

$1212233,1221233,1222133$&$4118,20841,112473,642253$\\\hline

$1211233$&$4118,20841,112473,642257$\\\hline

$1212133$&$4118,20841,112474,642283$\\\hline

$1212333,1221333$&$4118,20841,112475,642320,3857455,24238859$\\\hline

$1121334,1121343$&$4118,20841,112480,642458$\\\hline

$1221331$&$4118,20842,112498,642629$\\\hline

$1223311$&$4118,20842,112500,642701,3862079,24288261$\\\hline

$1221133$&$4118,20842,112500,642701,3862079,24288297$\\\hline

$1121345$&$4118,20842,112506,642847$\\\hline

$1121344$&$4118,20842,112508,642928$\\\hline

$1211314$&$4118,20843,112531,643262$\\\hline

$1222324$&$4118,20843,112531,643274$\\\hline

$1211341$&$4118,20843,112533,643325$\\\hline

$1222342$&$4118,20843,112533,643330$\\\hline

$1123141$&$4118,20845,112588,644220$\\\hline

$1123114$&$4118,20845,112588,644226$\\\hline

$1123411$&$4118,20845,112588,644230$\\\hline

$1121234$&$4119,20860,112699,644440$\\\hline

$1122234$&$4119,20861,112724,644829$\\\hline

$1122134$&$4119,20861,112725,644859$\\\hline

$1123132$&$4119,20862,112751,645259$\\\hline

$1123213$&$4119,20862,112752,645297$\\\hline

$1123123$&$4119,20862,112753,645319$\\\hline

$1121323$&$4119,20862,112757,645441$\\\hline

$1112342$&$4119,20862,112758,645467$\\\hline

$1112324$&$4119,20862,112760,645522$\\\hline

$1123312$&$4119,20863,112778,645710$\\\hline

$1121332$&$4119,20863,112778,645721$\\\hline

$1123223$&$4119,20863,112779,645747$\\\hline

$1123232$&$4119,20863,112780,645776,3892647$\\\hline

$1122323$&$4119,20863,112780,645776,3892648$\\\hline

$1122313$&$4119,20863,112781,645796$\\\hline

$1123321$&$4119,20863,112781,645802$\\\hline

$1123231$&$4119,20863,112781,645806$\\\hline

$1123332$&$4119,20863,112781,645807,3893233,24584926$\\\hline

$1123323$&$4119,20863,112781,645807,3893233,24584935$\\\hline

$1123233$&$4119,20863,112781,645808$\\\hline

$1211134,1222234$&$4119,20863,112781,645851$\\\hline

$1121233$&$4119,20863,112783,645862$\\\hline

$1123322$&$4119,20864,112805,646169$\\\hline

$1112334,1112343$&$4119,20864,112807,646204$\\\hline

$1122332$&$4119,20864,112807,646223$\\\hline

$1122233$&$4119,20864,112809,646282$\\\hline

$1122331$&$4119,20864,112810,646306$\\\hline

$1122133$&$4119,20864,112810,646307$\\\hline

$1122333$&$4119,20864,112811,646343$\\\hline

$1112345$&$4119,20864,112814,646384$\\\hline

$1112344$&$4119,20864,112814,646403$\\\hline

$1231313$&$4119,20866,112848,646718$\\\hline

$1231331$&$4119,20866,112850,646778$\\\hline

$1233131$&$4119,20866,112850,646779$\\\hline

$1221312$&$4119,20866,112850,646794$\\\hline

$1221231$&$4119,20866,112850,646795$\\\hline

$1212132$&$4119,20866,112850,646796$\\\hline

$1212231$&$4119,20866,112850,646797$\\\hline

$1212312$&$4119,20866,112851,646823$\\\hline

$1231212,1232323$&$4119,20866,112851,646824,3905223,24711200,163188860$\\\hline

$1212123$&$4119,20866,112851,646824,3905223,24711200,163188863$\\\hline

$1223121$&$4119,20866,112851,646824,3905223,24711200,163188865$\\\hline

$1232121,1233232$&$4119,20866,112851,646824,3905223,24711204$\\\hline

$1221213$&$4119,20866,112851,646824,3905223$\\\hline

$1213221$&$4119,20866,112851,646824,3905224$\\\hline

$1221321$&$4119,20866,112851,646824,3905230$\\\hline

$1213212$&$4119,20866,112851,646824,3905236$\\\hline

$1213122$&$4119,20866,112851,646824,3905254$\\\hline

$1231221,1232332$&$4119,20866,112851,646825$\\\hline

$1212213$&$4119,20866,112851,646826,3905284$\\\hline

$1212321$&$4119,20866,112851,646826,3905292$\\\hline

$1233113$&$4119,20866,112852,646837$\\\hline

$1221132$&$4119,20866,112852,646842,3905432$\\\hline

$1231133$&$4119,20866,112852,646842,3905461$\\\hline

$1213313$&$4119,20866,112852,646845$\\\hline

$1223112$&$4119,20866,112852,646847$\\\hline

$1213133$&$4119,20866,112852,646848$\\\hline

$1232112,1233223$&$4119,20866,112852,646850$\\\hline

$1221123$&$4119,20866,112852,646851$\\\hline

$1211232$&$4119,20866,112852,646856$\\\hline

$1211322$&$4119,20866,112853,646874$\\\hline

$1223233$&$4119,20866,112853,646886$\\\hline

$1222131$&$4119,20866,112854,646902$\\\hline

$1231122,1232233$&$4119,20866,112854,646908$\\\hline

$1211223$&$4119,20866,112854,646909$\\\hline

$1223323$&$4119,20866,112854,646912$\\\hline

$1222113$&$4119,20866,112856,646963$\\\hline

$1222311,1223211,1232211,1233311$&\\
$1233322$&$4119,20866,112856,646966$\\\hline

$1213331$&$4119,20866,112856,646967$\\\hline

$1223332$&$4119,20866,112856,646973$\\\hline

$1211333,1222333$&$4119,20866,112860,647094$\\\hline

$1121314$&$4119,20867,112883,647408$\\\hline

$1121341$&$4119,20867,112884,647438$\\\hline

$1233133$&$4119,20868,112906,647727$\\\hline

$1233313$&$4119,20868,112906,647729,3917135$\\\hline

$1231333$&$4119,20868,112906,647729,3917150$\\\hline

$1211312$&$4119,20868,112906,647739$\\\hline

$1212223,1212232,1212322,1213222$&\\
$1221223,1221232,1221322,1222123$&\\
$1222132,1222213,1222231,1222312$&\\
$1222321,1223122,1223212,1223221$&\\
$1231222,1232122,1232212,1232221$&\\
$1232333,1233233,1233323,1233331$&\\
$1233332$&$4119,20868,112906,647744,3917573$\\\hline

$1231131$&$4119,20868,112907,647757$\\\hline

$1213113$&$4119,20868,112907,647763$\\\hline

$1231113$&$4119,20868,112907,647764$\\\hline

$1213112$&$4119,20868,112907,647767$\\\hline

$1211132$&$4119,20868,112907,647769$\\\hline

$1231112,1232223$&$4119,20868,112907,647774,3918115$\\\hline

$1223223$&$4119,20868,112907,647774,3918120$\\\hline

$1211321$&$4119,20868,112907,647776$\\\hline

$1231121,1232232$&$4119,20868,112907,647778$\\\hline

$1211213$&$4119,20868,112907,647779$\\\hline

$1213131$&$4119,20868,112908,647788,3918150$\\\hline

$1231311$&$4119,20868,112908,647788,3918180$\\\hline

$1211313$&$4119,20868,112908,647796$\\\hline

$1212131$&$4119,20868,112908,647797$\\\hline

$1211231$&$4119,20868,112908,647801$\\\hline

$1213121$&$4119,20868,112908,647802$\\\hline

$1212113$&$4119,20868,112908,647804$\\\hline

$1211123$&$4119,20868,112908,647806,3918705$\\\hline

$1231211,1232322$&$4119,20868,112908,647806,3918717,24875277$\\\hline

$1223232$&$4119,20868,112908,647806,3918717,24875279$\\\hline

$1222323$&$4119,20868,112908,647806,3918717,24875282$\\\hline

$1212311$&$4119,20868,112909,647832$\\\hline

$1213211$&$4119,20868,112909,647836$\\\hline

$1221131$&$4119,20869,112936,648262$\\\hline

$1221311$&$4119,20869,112936,648263$\\\hline

$1221113$&$4119,20869,112936,648270$\\\hline

$1213333,1223333$&$4119,20869,112936,648283$\\\hline

$1213311$&$4119,20869,112937,648299$\\\hline

$1223111,1233111,1232111,1233222$&$4119,20869,112937,648301,3925257$\\\hline

$1211331$&$4119,20869,112937,648301,3925333$\\\hline

$1222332$&$4119,20869,112937,648308$\\\hline

$1223322$&$4119,20869,112937,648310$\\\hline

$1211133,1222233$&$4119,20869,112937,648316$\\\hline

$1112234$&$4120,20883,113036,648480$\\\hline

$1112323$&$4120,20883,113048,648817$\\\hline

$1112332$&$4120,20883,113048,648818$\\\hline

$1112233$&$4120,20884,113074,649231$\\\hline

$1123122$&$4120,20886,113110,649636$\\\hline

$1123212$&$4120,20886,113110,649640$\\\hline

$1121322$&$4120,20886,113110,649650$\\\hline

$1121232$&$4120,20886,113111,649672$\\\hline

$1121223$&$4120,20886,113111,649673$\\\hline

$1122312$&$4120,20886,113112,649692$\\\hline

$1122132$&$4120,20886,113112,649694,3934278$\\\hline

$1122123$&$4120,20886,113112,649694,3934279$\\\hline

$1123221$&$4120,20886,113114,649746$\\\hline

$1122213$&$4120,20886,113114,649751$\\\hline

$1122231$&$4120,20886,113114,649752,3935305$\\\hline

$1122321$&$4120,20886,113114,649752,3935315$\\\hline

$1121134$&$4120,20887,113134,650020$\\\hline

$1123133$&$4120,20887,113137,650074,3938938$\\\hline

$1123313$&$4120,20887,113137,650074,3938943$\\\hline

$1123331$&$4120,20887,113141,650187$\\\hline

$1121333$&$4120,20887,113145,650310$\\\hline

$1112314$&$4120,20888,113162,650462$\\\hline

$1112341$&$4120,20888,113162,650471$\\\hline

$1123112$&$4120,20888,113163,650505$\\\hline

$1121132$&$4120,20888,113163,650509$\\\hline

$1121312$&$4120,20888,113164,650531$\\\hline

$1123121$&$4120,20888,113164,650537$\\\hline

$1121213$&$4120,20888,113164,650539$\\\hline

$1121123$&$4120,20888,113164,650540$\\\hline

$1121231$&$4120,20888,113165,650564$\\\hline

$1121321$&$4120,20888,113165,650568$\\\hline

$1123113$&$4120,20889,113188,650884$\\\hline

$1122223,1122232,1122322,1123222$&$4120,20889,113188,650904$\\\hline

$1123131$&$4120,20889,113189,650907$\\\hline

$1121313$&$4120,20889,113189,650913$\\\hline

$1122131$&$4120,20889,113191,650978$\\\hline

$1122113$&$4120,20889,113191,650985$\\\hline

$1123211$&$4120,20889,113192,651006$\\\hline

$1122311$&$4120,20889,113192,651015$\\\hline

$1121331$&$4120,20890,113216,651360$\\\hline

$1123311$&$4120,20890,113216,651362$\\\hline

$1121133$&$4120,20890,113216,651375$\\\hline

$1123333$&$4120,20890,113216,651379$\\\hline

$1112223,1112232,1112322$&$4121,20904,113336,652122$\\\hline

$1112312$&$4121,20906,113380,652756$\\\hline

$1112132$&$4121,20906,113381,652784,3967043,25325818$\\\hline

$1112123$&$4121,20906,113381,652784,3967043,25325832$\\\hline

$1112333$&$4121,20906,113385,652864$\\\hline

$1112134$&$4121,20908,113412,653041$\\\hline

$1112213$&$4121,20907,113406,653173$\\\hline

$1112321$&$4121,20907,113407,653198$\\\hline

$1112231$&$4121,20907,113407,653204$\\\hline

$1112313$&$4121,20908,113425,653385$\\\hline

$1211131$&$4120,20895,113346,653419,3982042$\\\hline

$1211311$&$4120,20895,113346,653419,3982043$\\\hline

$1213111$&$4120,20895,113346,653419,3982063$\\\hline

$1211113,1222223,1222232,1222322$&\\
$1223222,1231111,1232222,1233333$&$4120,20895,113346,653419,3982093$\\\hline

$1112331$&$4121,20909,113450,653784$\\\hline

$1112133$&$4121,20909,113450,653791$\\\hline

$1111223,1111232$&$4122,20923,113586,654964$\\\hline

$1111234$&$4122,20926,113627,655233$\\\hline

$1111233$&$4122,20926,113643,655673$\\\hline

$1121131$&$4121,20914,113574,655691,4002872$\\\hline

$1121113$&$4121,20914,113574,655691,4002915$\\\hline

$1121311$&$4121,20914,113574,655692$\\\hline

$1123111$&$4121,20914,113575,655720$\\\hline

$1112113$&$4122,20931,113763,657478$\\\hline

$1112131$&$4122,20931,113763,657480$\\\hline

$1112311$&$4122,20931,113765,657533$\\\hline

$1111213$&$4123,20946,113917,658862$\\\hline

$1111231$&$4123,20946,113920,658937$\\\hline

$1111123$&$4124,20959,114044,660003$\\\hline

$1221221$&$4124,20966,114159,661276$\\\hline

$1212221$&$4124,20966,114159,661277$\\\hline

$1212122$&$4124,20966,114159,661279,4052523$\\\hline

$1221212$&$4124,20966,114159,661279,4052524,26151249,176986852$\\\hline

$1222121$&$4124,20966,114159,661279,4052524,26151249,176986866$\\\hline

$1212212$&$4124,20966,114159,661279,4052526$\\\hline

$1212112$&$4124,20966,114159,661280$\\\hline

$1211212$&$4124,20966,114159,661282$\\\hline

$1121212$&$4124,20966,114159,661284$\\\hline

$1212121$&$4124,20966,114159,661286$\\\hline

$1221121$&$4124,20966,114161,661324$\\\hline

$1222112$&$4124,20966,114161,661326$\\\hline

$1221112$&$4124,20966,114161,661328$\\\hline

$1221122$&$4124,20966,114161,661329,4053313$\\\hline

$1221211$&$4124,20966,114161,661329,4053324$\\\hline

$1212211$&$4124,20966,114161,661330$\\\hline

$1121221$&$4124,20966,114161,661331,4053377,26162919$\\\hline

$1211221$&$4124,20966,114161,661331,4053377,26162926$\\\hline

$1122212$&$4124,20966,114161,661332,4053398$\\\hline

$1122122$&$4124,20966,114161,661332,4053400$\\\hline

$1211222$&$4124,20966,114161,661333,4053404$\\\hline

$1121122$&$4124,20966,114161,661333,4053421$\\\hline

$1122112$&$4124,20966,114161,661333,4053424$\\\hline

$1211122$&$4124,20966,114161,661334,4053429$\\\hline

$1121222$&$4124,20966,114161,661333,4053431$\\\hline

$1122121$&$4124,20966,114161,661334,4053447$\\\hline

$1112122$&$4124,20966,114161,661337$\\\hline

$1222211$&$4124,20966,114161,661339$\\\hline

$1112212$&$4124,20966,114161,661340,4053579$\\\hline

$1122221$&$4124,20966,114161,661340,4053593$\\\hline

$1122211$&$4124,20966,114165,661427$\\\hline

$1222111$&$4124,20966,114165,661429$\\\hline

$1112221$&$4124,20966,114165,661432$\\\hline

$1112222$&$4124,20966,114165,661439$\\\hline

$1111222$&$4124,20966,114165,661444$\\\hline

$1212222,1221222,1222122,1222212$&\\
$1222221$&$4124,20968,114209,662046$\\\hline

$1121112$&$4124,20968,114210,662074,4062705$\\\hline

$1211112$&$4124,20968,114210,662074,4062709,26270000$\\\hline

$1211121$&$4124,20968,114210,662074,4062709,26270389$\\\hline

$1112112$&$4124,20968,114210,662075,4062730$\\\hline

$1121121$&$4124,20968,114210,662075,4062742$\\\hline

$1211211$&$4124,20968,114210,662075,4062747$\\\hline

$1112121$&$4124,20968,114211,662099,4063105,26275156$\\\hline

$1212111$&$4124,20968,114211,662099,4063105,26275158,178364080$\\\hline

$1111212$&$4124,20968,114211,662099,4063105,26275158,178364101$\\\hline

$1121211$&$4124,20968,114211,662099,4063106$\\\hline

$1122222$&$4124,20969,114233,662412$\\\hline

$1221111$&$4124,20969,114234,662438,4067231$\\\hline

$1111221$&$4124,20969,114234,662438,4067256$\\\hline

$1112211$&$4124,20969,114234,662438,4067265$\\\hline

$1122111$&$4124,20969,114234,662438,4067266$\\\hline

$1111122$&$4124,20969,114234,662438,4067287$\\\hline

$1111112,1111121,1111211,1112111$&\\
$1121111,1211111,1222222$&$4125,20990,114516,665604$\\\hline

$1111111$&$4131,21065,115274,672673$\\\hline

\end{longtable}
}

\end{document}